\documentclass{tac} 

\PassOptionsToPackage{hyphens}{url}\usepackage[colorlinks=true]{hyperref}
\hypersetup{urlcolor=[rgb]{0.1,0.1,0.4}}
\usepackage{etex}
\usepackage[all]{xy}
\usepackage{pb-diagram}
\usepackage{amsfonts, amssymb, amsmath}
\usepackage{tikz}
\usepackage{epsfig}
\usepackage{latexsym,url}
\usepackage{txfonts}
\usepackage[tiling]{pst-fill}
\usepackage{pstcol}
\usepackage{pdfpages}

\usetikzlibrary{matrix,arrows}
\usetikzlibrary{backgrounds}
\usetikzlibrary{decorations.pathmorphing}
\usetikzlibrary{decorations.pathreplacing}
\usetikzlibrary{decorations.shapes}
\usetikzlibrary{decorations.markings}
\usetikzlibrary{calc}

\tikzset{l/.style={font=\fontsize{8}{8}\selectfont}}
\tikzset{m/.style={font=\fontsize{6}{6}\selectfont}}
\tikzset{
electron/.style={postaction=decorate, decoration={markings, mark=at position 0.55 with {\arrow{triangle 45}}}}, photon/.style={decorate, decoration={coil, aspect=0}}
}
\definecolor{salmon}{rgb}{1,0.47,0.425}

\newtheorem{thm}{Theorem}

\newtheorem{defn}{Definition}
\newtheorem{lemma}{Lemma}

\newcommand{\Set}{\mathrm{Set}}

\newcommand{\Cat}{{\rm Cat}}

\newcommand{\maps}{\colon}

\newcommand{\lhom}{\multimap}

\newcommand{\tensor}{\otimes}
\newcommand{\x}{\times}

\newcommand{\id}{{\rm i}}

\newcommand{\op}{{\rm op}}

\renewcommand{\text}{\mbox}

\newcommand{\Span}{\mbox{Span}}

\newcommand{\nCob}{n\mbox{Cob}}

\newcommand{\C}{ {\mathcal{K}}  }

\newcommand{\acirc}{ {\stackrel{\mbox{\circle{2}}}{a}} }
\newcommand{\lcirc}{ {\stackrel{\mbox{\circle{2}}}{l}} }
\newcommand{\rcirc}{ {\stackrel{\mbox{\circle{2}}}{r}} }

\newtheorem{cor}[thm]{Corollary}

\author{Michael Stay}

\thanks{Thanks to John Baez and Mike Shulman for explanations of concepts and helpful
comments; to Alex Hoffnung for his patience in explaining weak pullbacks;
to Richard Garner and an anonymous reader for very careful readings and catching some errors,
and to the members of the n-Category Caf\'e and MathOverflow for 
expository background material and references.
}

\title{Compact closed bicategories}
\copyrightyear{2016}
\eaddress{metaweta@gmail.com}
\date{\today}

\begin{document}
\maketitle

\begin{abstract}
A compact closed bicategory is a symmetric monoidal bicategory where every object is equipped with a weak dual. The unit and counit satisfy the usual ``zig-zag'' identities of a compact closed category only up to natural isomorphism, and the isomorphism is subject to a coherence law.  We give several examples of compact closed bicategories, then review previous work.  In particular, Day and Street defined compact closed bicategories indirectly via Gray monoids and then appealed to a coherence theorem to extend the concept to bicategories; we restate the definition directly.

We prove that given a 2-category $T$ with finite products and weak pullbacks, the bicategory of objects of $C$, spans, and isomorphism classes of maps of spans is compact closed.  As corollaries, the bicategory of spans of sets and certain bicategories of ``resistor networks'' are compact closed.
\end{abstract}

\section{Introduction}

When moving from set theory to category theory and higher to $n$-category theory, we typically weaken equations at one level to natural isomorphisms at the next.  These isomorphisms are then subject to new coherence equations.  For example, multiplication in a monoid is associative, but the tensor product in a monoidal category is, in general, only associative up to a natural isomorphism.  This ``associator'' natural isomorphism has to satisfy some extra equations that are trivial in the case of the monoid.  In a similar way, when we move from compact closed categories to compact closed bicategories, the ``zig-zag'' equations governing the dual get weakened to natural isomorphisms and we need to introduce some new coherence laws.

In Section \ref{examples}, we will give several examples of important mathematical structures and how they arise in relation to compact closed bicategories.  Following the examples, in Section \ref{previous} we give the history of compact closed bicategories and related work.  Next, in Section \ref{definition} we give the complete definition, which to our knowledge has not appeared elsewhere; we try to motivate each piece of the definition so that the reader could afterwards reconstruct the definition without the aid of this paper.  In Section \ref{Hoffnung}, we prove that a construction by Hoffnung is an instance of a compact closed bicategory, and obtain few others as corollaries.

\section{Examples}
\label{examples}

In order to get across some of the flavor of these different compact closed bicategories, we will describe the bicategories as well as some weak monoids and monads in them.  That these bicategories are compact closed is mostly folklore; we prove that a few of them are compact closed as corollaries of the main theorem at the end of this paper.  The weak monoids and monads play no role in the rest of the paper, but we have found that comparing and contrasting them helps when trying to develop intuition about the bicategories.  The monoids are variations on the notion of an associative algebra, while the monads are variations on the notion of a category.

\begin{itemize}
  \item A {\bf span} from $A$ to $B$ in a category $T$ is an object $C$
in $T$ together with an ordered pair of morphisms $(f\maps C \to A, g\maps C\to B).$
  \begin{center}
    \begin{tikzpicture}[scale=2]
      \node (A) at (0,0) {$A$};
      \node (B) at (2,0) {$B$};
      \node (C) at (1,1) {$C$}
        edge [->] node [above left,l] {$f$} (A)
        edge [->] node [above right,l] {$g$} (B);
    \end{tikzpicture}
  \end{center}
  If $T$ is a category with pullbacks, we can compose spans:
  \begin{center}
    \begin{tikzpicture}[scale=2]
      \node (A) at (0,0) {$A$};
      \node (B) at (2,0) {$B$};
      \node (C) at (4,0) {$C$};
      \node (D) at (1,1) {$D$}
        edge [->] node [above left,l] {$f$} (A)
        edge [->] node [above right,l] {$g$} (B);
      \node (E) at (3,1) {$E$}
        edge [->] node [above left,l] {$h$} (B)
        edge [->] node [above right,l] {$j$} (C);
      \node (DghE) at (2,2) {$D_{gh}E$}
        edge [->] node [above left,l] {$\pi_1$} (D)
        edge [->] node [above right,l] {$\pi_2$} (E);
      \draw (1.6, 1.4)--(2,1)--(2.4, 1.4);
    \end{tikzpicture}
  \end{center}
  A {\bf map of spans} between two spans $A \stackrel{f}{\leftarrow} C \stackrel{g}{\rightarrow} B$ and $A \stackrel{f'}{\leftarrow} C' \stackrel{g'}{\rightarrow} B$ is a morphism $h\maps C \to C'$ making the following diagram commute:
  \begin{center}
    \begin{tikzpicture}[scale=2]
      \node (A) at (0,0) {$A$};
      \node (B) at (2,0) {$B$};
      \node (C) at (1,1) {$C$}
        edge [->] node [above left,l] {$f$} (A)
        edge [->] node [above right,l] {$g$} (B);
      \node (C') at (1, -1) {$C'$}
        edge [->] node [below left,l] {$f'$} (A)
        edge [->] node [below right,l] {$g'$} (B)
        edge [<-] node [right,l] {$h$} (C);
    \end{tikzpicture}
  \end{center}
Since the pullback is associative only up to a natural isomorphism, the same is true of the composite of two spans, so this construction does not give a 2-category; however, we do get a bicategory $\Span(T)$ of objects of $T$, spans in $T$, and maps of spans.

  If $T$ is a category with finite products as well as pullbacks, then the bicategory $\Span(T)$ is a compact closed bicategory where the tensor product is given by the product in $T.$  A weak monoid object in $\Span(T)$ is a categorification of the notion of an associative algebra.  For example, one weak monoid in Span(Set) is equivalent to the category of polynomial functors from Set to itself; such functors can be ``added'' using disjoint union, ``multiplied'' using the cartesian product, and ``scaled'' by sets \cite{GK}.  A monad in Span($T$) is a category internal to $T$ \cite{Ben67}.
  
  \item Sets, relations, and implications form the compact closed bicategory Rel, where the tensor product is given by the product in Set.  A weak monoid object $M$ in Rel is a quantale on the powerset of $M$ \cite{ShulCafe}, while a monad in Rel is a preorder.

  \item A {\bf 2-rig} is a cocomplete monoidal category where the tensor product distributes over the colimits \cite{HDA3}, though for the purpose of constructing a compact closed bicategory we only need the tensor product to distribute over finite coproducts.  Given a symmetric 2-rig $R$, Mat($R$) is the bicategory of finitely-generated free $R$-modules, where the tensor product is the usual tensor product for matrices.  We expect it to be compact closed.  A weak monoid object in Mat($R$) is a categorified finite-dimensional associative algebra over $R.$  A monad in Mat($R$) is a finite $R$-enriched category.  A {\bf finite 2-rig} $S$ is only finitely cocomplete, but we expect Mat($S$) is still compact closed. Kapranov and Voevodsky \cite{KV94} described a bicategory equivalent to Mat(FinVect) and called its objects ``2-vector spaces''.

  \item For a category $C$, let $\widehat{C} = \Set^{C^{\op}}$ be the category of presheaves on $C$.  The 2-category Cocont has
    \begin{itemize}
      \item small categories as objects;
      \item cocontinuous functors $f\maps \widehat{C} \to \widehat{D}$ between the categories of presheaves on the source and target as morphisms;
      \item natural transformations as 2-morphisms.
    \end{itemize}
    We can think of cocontinuous functors as being ``Set-linear transformations'', since they preserve sums.  Day and Street \cite{DS97}
proved that Cocont is a compact closed 2-category, {\em i.e.} a compact closed bicategory where the associator and unitors for composition are equalities.

  \item Recall that a profunctor $F\maps C \not\to R$ is a functor $F\maps C^{\op} \times R \to \Set$; we can think of profunctors as being rather like matrices, where the set $F(c,r)$ is the ``matrix element'' at row $r$ and column $c$.  Composition of profunctors is given by taking the coend of the inner coordinates, just as matrix multiplication done by summing over the inner index.  Small categories, profunctors, and natural transformations form the compact closed bicategory Prof, where the tensor product is the product in Cat.  A weak monoid object in Prof is a promonoidal category \cite{Day, Houston}. A symmetric monoidal monad in Prof is a Freyd category, also known as an ``Arrow'' in the functional programming community \cite{Asada, JHH}.

    Cattani and Winskel \cite{CW04} showed that Cocont and Prof are equivalent as bicategories.  Though they do not explicitly state it, the equivalence they construct is symmetric monoidal; since symmetric monoidal equivalences preserve the dual, Cocont and Prof are equivalent as compact closed bicategories.
  \item So far the examples have been rather algebraic in flavor, but there are topological examples, too.  The category $\nCob$ is the compact closed category whose
    \begin{itemize}
      \item objects are $(n-1)$-dimensional manifolds and
      \item morphisms are diffeomorphism classes of collared $n$-dimensional cobordisms between them,
    \end{itemize}
    where the tensor product is disjoint union.  Atiyah \cite{AtiyahTQFT} introduced the category informally in his paper defining topological quantum field theories.

    Morton \cite{Morton} defined the bicategory $\nCob_2$ whose
    \begin{itemize}
      \item objects are $(n-2)$-dimensional manifolds,
      \item morphisms are collared $(n-1)$-dimensional cobordisms, or ``manifolds with boundary'', and
      \item 2-morphisms are diffeomorphism classes of collared $n$-dimensional maps of cobordisms, or ``manifolds with corners''.
    \end{itemize}
    The collars are necessary to preserve the smoothness when composing 1- and 2-morphisms.  Schommer-Pries proved a purely algebraic characterization of $2\mbox{Cob}_2$, essentially proving the ``Baez-Dolan cobordism hypothesis'' for the $n=2$ case \cite{HDATQFT}.  We expect that $\nCob_2$ is compact closed.
    
    \item In a letter to the author, John Baez defined two interesting compact closed bicategories.
    
    A {\bf directed multigraph} is a finite set $E$ of edges and a finite set $V$ of vertices equipped with functions $s, t\maps E \to V$
mapping each edge to its source and target.  A {\bf resistor network} is a directed multigraph equipped with a function $r$ assigning a resistance in $(0, \infty)$ to each edge:
    \begin{center}
      \begin{tikzpicture}[scale=2]
        \node (R) at (0,1) {$(0, \infty)$};
        \node (E1) at (1,1) {$E$}
          edge [->] node [below] {$r$} (R);
        \node (V1) at (2,1) {$V$};
        \draw [->] (E1.15) to node [above] {$s$} (V1.165);
        \draw [->] (E1.345) to node [below] {$t$} (V1.195);
      \end{tikzpicture}
    \end{center}
    There are various choices one could make for a morphism of such networks; Baez defined a morphism of resistor networks to be a pair of functions $\epsilon, \upsilon$ making the following diagrams commute:

    \begin{center}
      \begin{tikzpicture}[scale=2]
        \node (R) at (0,.5) {$(0, \infty)$};
        \node (E2) at (1, 0) {$E'$}
          edge [->] node [below] {$r'$} (R);
        \node (E1) at (1,1) {$E$}
          edge [->] node [above] {$r$} (R)
          edge [->] node [right] {$\epsilon$} (E2);
      \end{tikzpicture}
      $\quad$
      \begin{tikzpicture}[scale=2]
        \node (E2) at (1, 0) {$E'$};
        \node (E1) at (1,1) {$E$}
          edge [->] node [right] {$\epsilon$} (E2);
        \node (V2) at (2, 0) {$V'$}
          edge [<-] node [below] {$s$} (E2);
        \node (V1) at (2,1) {$V$}
          edge [<-] node [above] {$s$} (E1)
          edge [->] node [right] {$\upsilon$} (V2);
      \end{tikzpicture}
      $\quad$
      \begin{tikzpicture}[scale=2]
        \node (E2) at (1, 0) {$E'$};
        \node (E1) at (1,1) {$E$}
          edge [->] node [right] {$\epsilon$} (E2);
        \node (V2) at (2, 0) {$V'$}
          edge [<-] node [below] {$t$} (E2);
        \node (V1) at (2,1) {$V$}
          edge [<-] node [above] {$t$} (E1)
          edge [->] node [right] {$\upsilon$} (V2);
      \end{tikzpicture}
    \end{center}

    Resistor networks and morphisms between them form a category ResNet; this category has finite limits and colimits.
    
    There is a compact closed bicategory Cospan(ResNet) with an important compact closed subbicategory Circ consisting of cospans whose feet are resistor networks with no edges.  A morphism in Circ is a {\bf circuit,} a resistor network with chosen sets of input and output vertices across which one can measure a voltage drop.
\end{itemize}

\section{Previous work}
\label{previous}

Compact closed categories were first defined by Kelly \cite{KellyCC}, and later studied in depth by Kelly and Laplaza \cite{KellyLaplaza}.

B\'enabou \cite{Ben67} defined bicategories and showed that small categories, distributors, and natural transformations form a bicategory Dist.  Distributors later became more widely known as ``profunctors'', so we will call that bicategory ``Prof'' instead.  Later, B\'enabou defined closed bicategories and showed that Prof is closed  \cite
{Ben73}.  He defined $V$-enriched profunctors when $V$ is a cocomplete monoidal or symmetric monoidal closed category, defined $V$-Prof and proved that any $V$-enriched functor, regarded as a $V$-profunctor, has a right adjoint.  More applications and details are in his lecture notes
\cite{Ben00}.

Kapranov and Voevodsky \cite{KV94} defined braided semistrict monoidal 2-categories, but their definition left out some necessary axioms.  Baez and Neuchl \cite{HDA1} gave an improved definition, but it was still missing a clause; Crans \cite{Crans} gave the complete definition.  See Baez and Langford \cite{HDA4} and Shulman \cite{Shulman} for details.

Gordon, Power, and Street \cite{GPS} defined fully weak tricategories; a monoidal bicategory is a one-object tricategory.

Another name for semistrict monoidal 2-categories is ``Gray monoids'', {\em i.e.} monoid objects in the 2-category Gray \cite{Gray}.  Day and Street \cite{DS97} defined compact closed Gray monoids, and appealed to the coherence theorem of Gordon, Power, and Street to extend compact closedness to arbitrary bicategories.  The semistrict approach is somewhat artificial when dealing with most ``naturally occurring'' bicategories, since the associator for composition of 1-morphisms is rarely the identity.  Cocont is a notable exception.

Katis, Sabadini and Walters \cite{KSW98} gave a precise account of the double-entry bookkeeping method {\em partita doppia} in terms of the compact closed bicategory Span(RGraph); in a later paper \cite{KW99}, they cite a handwritten note by McCrudden for the ``swallowtail'' coherence law we use in this paper.

Preller and Lambek \cite{PL07} generalized compact monoidal categories in a different direction.  They considered a compact monoidal category to be a one-object bicategory satisfying some criteria, and then extend that definition to multiple objects.  The resulting concept of ``compact bicategory'' is {\em not} what is being studied in this paper.

McCrudden \cite{McCrudden} gave the first fully general definitions of braided, sylleptic, and symmetric monoidal bicategories.  Schommer-Pries \cite{SPT} gave the correct notion of a monoidal transformation between monoidal functors between monoidal bicategories.

Carboni and Walters \cite{CW87} proved that $V$-Prof is a cartesian bicategory.  Later, they showed \cite{CW04} that Prof is equivalent to Cocont as a bicategory.  Together with Kelly and Wood \cite{CKWW08}, they proved that any cartesian bicategory is symmetric monoidal in the sense of McCrudden.

Gurski and Osorno \cite{GO} proved that every symmetric monoidal bicategory is equivalent to a semistrict one. Schommer-Pries \cite{SPT} strengthened their result by proving that every symmetric monoidal bicategory is equivalent to a ``quasistrict symmetric'' monoidal bicategory. Bartlett \cite{Bartlett} went a step further and showed every symmetric monoidal bicategory is equivalent to a ``stringent'' one.  He also used Schommer-Pries' results to develop a graphical calculus for symmetric monoidal bicategories.

\section{Compact closed bicategories}
\label{definition}

In this section, we lay out the definition of a compact closed bicategory.  First we give the definition of a bicategory, then start adding structure to it: we introduce the tensor product and monoidal unit; then we look at the different ways to move objects around each other, giving braided, sylleptic and symmetric monoidal bicategories. Next, we define closed monoidal bicategories by introducing a right pseudoadjoint to tensoring with an object; and finally we introduce duals for objects in a bicategory.

\begin{defn}
  A {\bf bicategory} $\C$ consists of
  \begin{enumerate}
    \item a collection of {\bf objects}
    \item for each pair of objects $A, B$ in $\C$, a category $\C(A,B)$;
the objects of $\C(A,B)$ are called {\bf 1-morphisms}, while the morphisms of $\C(A,B)$ are called {\bf 2-morphisms}.
    \item for each triple of objects $A, B, C$ in $\C$, a {\bf composition} functor 
        \[ \circ_{A,B,C}\maps \C(B, C) \times \C(A, B) \to \C(A, C).\]
      We will leave off the indices and write it as an infix operator.
    \item for each object $A$ in $\C$, an object $1_A$ in $\C(A,A)$ called the {\bf identity 1-morphism on $A$}.  We will often write this simply as $A$.
    \item for each quadruple of objects $A, B, C, D$, a natural isomorphism called the {\bf associator for composition}; if $(f,g,h)$ is an object of $\C(C,D) \times \C(B,C) \times \C(A,B),$ then 
      $$\acirc_{f,g,h} \maps (f \circ g) \circ h \to f \circ (g \circ h).$$
    \item for each pair of objects $A, B$ in $\C$, natural isomorphisms called {\bf left and right unitors for composition.}  If $f$ is an object of $\C(A,B)$, then
      \[\begin{array}{l}
        \lcirc_f\maps B \circ f \stackrel{\sim}{\to} f\\
        \rcirc_f\maps f \circ A \stackrel{\sim}{\to} f
      \end{array}\]
  \end{enumerate}
  such that $\acirc, \lcirc,$ and $\rcirc$ satisfy the following coherence laws:
  \begin{enumerate}
    \item for all $(f,g,h,j)$ in $\C(D,E) \times \C(C,D) \times \C(B,C)
\times \C(A,B),$ the following diagram, called the {\bf pentagon equation}, commutes:
      \[\begin{diagram}
        \node[2]{((f \circ g) \circ h) \circ j} \arrow{sse,t}{\acirc_{f\circ g, h, j}} \arrow{sw,t}{\acirc_{f,g,h} \circ j} \\
        \node{(f \circ (g \circ h)) \circ j} \arrow[2]{s,l}{\acirc_{f,g\circ h, j}}\\
        \node[3]{(f \circ g) \circ (h \circ j)} \arrow{ssw,b}{\acirc_{f,g,h \circ j}}\\
        \node{f \circ ((g \circ h) \circ j)} \arrow{se,b}{f \circ\acirc_{g,h,j}}\\
        \node[2]{f \circ (g \circ (h \circ j))}
      \end{diagram}\]
    \item for all $(f,g)$ in $\C(B,C)\times \C(A,B)$ the following diagram, called the {\bf triangle equation}, commutes:
      \[\begin{diagram}
        \node{(f \circ B) \circ g} 
          \arrow[2]{e,t}{\acirc} 
          \arrow{se,b} {\rcirc_f \circ g} \node[2]{f \circ (B \circ g)} 
          \arrow{sw,b}{f \circ \lcirc_g} \\
        \node[2]{f \circ g}
      \end{diagram}\]
  \end{enumerate}
\end{defn}

The associator $\acirc$ and unitors $\rcirc, \lcirc$ for composition of 1-morphisms are necessary, but when we are drawing commutative diagrams of 1-morphisms they are very hard to show; fortunately, by the coherence theorem for bicategories \cite{LeinsterBB}, any consistent choice is equivalent to any other, so we leave them out.

We refer the reader to Tom Leinster's excellent ``Basic bicategories''
\cite{LeinsterBB} for definitions of
\begin{itemize}
  \item morphisms of bicategories, which we call functors,
  \item transformations between functors, which we call pseudonatural transformations, and
  \item modifications between transformations.
\end{itemize}

\begin{defn}
  An {\bf equivalence} of objects $A, B$ in a bicategory is a pair of morphisms ${f\maps A \to B,}$ ${g\maps B \to A}$ together with invertible 2-morphisms ${e\maps g \circ f \stackrel{\sim}{\Rightarrow} 1_A}$ and ${i\maps f \circ g \stackrel{\sim}{\Rightarrow} 1_B}$.
\end{defn}
\begin{defn}
  An {\bf adjoint equivalence} is one in which the 2-morphisms $e$ and $i^{-1}$ exhibit that $g$ is left adjoint to $f$.
\end{defn}
For a given morphism $f,$ any two choices of data $(g, e, i)$ making
$f$ an adjoint equivalence are canonically isomorphic, so any choice is as good as any other.  When $f, g$ form an adjoint equivalence, we write $g = f^\bullet$.  Any equivalence can be improved to an adjoint equivalence.

We can often take a 2-morphism and ``reverse'' one of its edges.  Given objects $A,B,C,D$, morphisms $f\maps A\to C$, $g\maps C\to D$, $h\maps D\to B$, $j\maps A \to B$ such that $h$ is an adjoint equivalence, and a 2-morphism 
\begin{center}
  \begin{tikzpicture}[scale=2]
    \node (A) at (0,0) {$A$};
    \node (alpha) at (1,0) {$\Downarrow \alpha$};
    \node (C) at (.67, -.5) {$C$}
      edge [<-] node [below left, l] {$f$} (A);
    \node (D) at (1.33, -.5) {$D$}
      edge [<-] node [below, l] {$g$} (C);
    \node (B) at (2,0) {$B$}
      edge [<-, out=135, in=45] node [above,l] {$j$} (A)
      edge [<-] node [below right, l] {$h$} (D);
  \end{tikzpicture}
\end{center}
we can get a new 2-morphism 
\[ (\rcirc(g) \circ f)(e_h \circ g \circ f)(h^\bullet \circ
\alpha)\maps h^\bullet \circ j \Rightarrow g \circ f, \]
\begin{center}
  \begin{tikzpicture}[scale=2]
    \node (A) at (0,0) {$A$};
    \node (alpha) at (1,0) {$\Downarrow \alpha$};
    \node (C) at (.67, -.5) {$C$}
      edge [<-] node [below left, l] {$f$} (A);
    \node (D) at (1.33, -.5) {$D$}
      edge [<-] node [below, l] {$g$} (C);
    \node (B) at (2,0) {$B$}
      edge [<-, out=135, in=45] node [above,l] {$j$} (A)
      edge [<-] node [above left, l] {$h$} (D);
    \node (D2) at (2.67, -.5) {$D$}
      edge [<-] node [above right, l] {$h^\bullet$} (B)
      edge [<-] node [below, l] {$1_D$} (D)
      edge [<-, out=-135, in=-45] node [below] {$g$} (C);
    \node at (2, -.3) {$\Downarrow e_h$};
    \node at (1.65, -.75) {$\Downarrow \rcirc(g)$};
    \node at (3,0) {$=$};
  \end{tikzpicture}  
  \begin{tikzpicture}[scale=2]
    \node (A) at (0,0) {$A$};
    \node (alpha) at (1,0) {$\Downarrow \alpha_1$};
    \node (C) at (.67, -.5) {$C$}
      edge [<-] node [below left, l] {$f$} (A);
    \node (D) at (1.33, -.5) {$D$}
      edge [<-] node [below, l] {$g$} (C);
    \node (B) at (2,0) {$B$}
      edge [<-, out=135, in=45] node [above,l] {$j$} (A)
      edge [->] node [below right, l] {$h^\bullet$} (D);
    \node at (0,-1.2) {};
  \end{tikzpicture}  
\end{center}
where $e_h\maps h^\bullet \circ h \Rightarrow 1$ is the 2-morphism from the equivalence.  We denote such variations of a 2-morphism $\alpha$ by adding numeric subscripts as in $\alpha_1$; the number simply records the order in which we introduce them, not any information about the particular variation.

In the following definitions, I have given some plausible combinatorical reasoning justifying many of the parts of the definition, but except where noted, this is not part of the definition; its intent is merely to help organize the rather long and dry content.  I am not aware of any work on the combinatorics of cells in higher categories beyond that mentioned below by Stasheff, Kapranov and Voevodsky.  

Also, some of the illustrations of 2-morphisms and coherence laws below are quite large.  In order to preserve legibility, I use expressions like $(A B) C$ as a shorthand for functors like
      \[\tensor \circ (\tensor \times 1)\maps  \C^3 \to \C,\]
since parentheses suffice to show where the tensor product should be.

\begin{defn} A {\bf monoidal bicategory} $\C$ is a bicategory in which we can ``multiply'' objects.  Monoidal bicategories were originally defined as one-object ``tricategories'' \cite{GPS}; unpacking that definition, a monoidal bicategory consists of the following: \\
  
\begin{itemize}
\item A bicategory $\C$.
\item A {\bf tensor product} functor $\otimes\maps \C \times \C \rightarrow \C$.  This functor involves an invertible ``tensorator'' 2-morphism $(f \otimes g) \circ (f' \otimes g') \Rightarrow (f \circ f') \otimes (g \circ g')$ which we elide in most of the coherence equations below.  The coherence theorem for monoidal bicategories implies that any 2-morphism involving the tensorator is the same no matter how it is inserted \cite[Remark~3.1.6]{Gurski}, so like the associator for composition of 1-morphisms, we leave it out.

The {\bf Stasheff polytopes} \cite{Stasheff} are a series of geometric figures whose vertices enumerate the ways to parenthesize the tensor product of $n$ objects, so the number of vertices is given by the Catalan numbers; for each polytope, we have a corresponding $(n-2)$-morphism of the same shape with directed edges and faces:
  \begin{enumerate}
    \item The tensor product of one object $A$ is the one object $A$ itself.
    \item The tensor product of two objects $A$ and $B$ is the one object $(A B)$.
    \item There are two ways to parenthesize the product of three objects, so we have an {\bf associator} adjoint equivalence pseudonatural in $A,B,C$
        \[ a\maps (A B) C \rightarrow A (B C) \]
      for moving parentheses from the left pair to the right pair.  The fact that $a$ is pseudonatural in $A,B,C$ means that given $f\maps A\to D, g\maps B\to E,$ and $g\maps C\to F,$ there is an invertible modification from $f(gh) \circ a_{ABC}$ to $a_{DEF} \circ (fg)h$; this invertible modification appears three times in the ``associahedron'' below on the green faces.
    \item There are five ways to parenthesize the product of four objects, so we have a {\bf pentagonator} invertible modification $\pi$ relating the two different ways of moving parentheses from being clustered at the left to being clustered at the right.  (Mnemonic: Pink Pentagonator.)
      
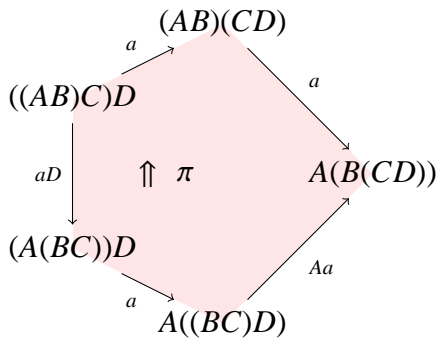
\begin{center}
        \begin{tikzpicture}
          \filldraw[white,fill=red,fill opacity=0.1](0,3)--(2,4)--(4,2)--(2,0)--(0,1)--cycle;
          \node (ABCD1) at (0,3) {$((A B) C) D$};
          \node (ABCD2) at (2,4) {$(A B) (C D)$}
            edge [<-] node [l, above left] {$a$} (ABCD1);
          \node (ABCD3) at (4,2) {$A (B (C D))$}
            edge [<-] node [l, above right] {$a$} (ABCD2);
          \node (ABCD4) at (2,0) {$A ((B C) D)$}
            edge [->] node [l, below right] {$A a$} (ABCD3);
          \node (ABCD5) at (0,1) {$(A (B C)) D$}
            edge [->] node [l, below left] {$a$} (ABCD4)
            edge [<-] node [l, left] {$a D$} (ABCD1);
          \node at (1,2) {\tikz\node [rotate=90] {$\Rightarrow$};};
          \node at (1.5,2) {$\pi$};
        \end{tikzpicture}
      \end{center}
    \item There are fourteen ways to parenthesize the product of five objects, so we have an {\bf associahedron} equation of modifications with fourteen vertices relating the various ways of getting from the parentheses clustered at the left to clustered at the right.
      
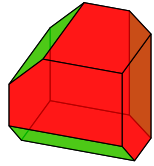
\begin{center}
        \begin{tikzpicture}[line join=round]
          \begin{scope}[scale=.2]
            \filldraw[fill=red,fill opacity=0.7](-4.306,-3.532)--(-2.391,-.901)--(-2.391,3.949)--(-5.127,.19)--(-5.127,-2.581)--cycle;
            \filldraw[fill=red,fill opacity=0.7](-4.306,-3.532)--(-2.391,-.901)--(2.872,-1.858)--(4.306,-3.396)--(3.212,-4.9)--cycle;
            \filldraw[fill=red,fill opacity=0.7](2.872,-1.858)--(2.872,5.07)--(-.135,5.617)--(-2.391,3.949)--(-2.391,-.901)--cycle;
            \filldraw[fill=green,fill opacity=0.7](4.306,-3.396)--(4.306,3.532)--(2.872,5.07)--(2.872,-1.858)--cycle;
            \filldraw[fill=green,fill opacity=0.7](-2.872,1.858)--(-.135,5.617)--(-2.391,3.949)--(-5.127,.19)--cycle;
            \filldraw[fill=red,fill opacity=0.7](3.212,-4.9)--(4.306,-3.396)--(4.306,3.532)--(2.391,.901)--(2.391,-3.949)--cycle;
            \filldraw[fill=green,fill opacity=0.7](-4.306,-3.532)--(3.212,-4.9)--(2.391,-3.949)--(-5.127,-2.581)--cycle;
            \filldraw[fill=red,fill opacity=0.7](-2.872,1.858)--(2.391,.901)--(4.306,3.532)--(2.872,5.07)--(-.135,5.617)--cycle;
            \filldraw[fill=red,fill opacity=0.7](-5.127,-2.581)--(-5.127,.19)--(-2.872,1.858)--(2.391,.901)--(2.391,-3.949)--cycle;
          \end{scope}
        \end{tikzpicture}
      \end{center}
    The associahedron is a cube with three of its edges bevelled. It holds in the bicategory $\C$, where the unmarked 2-morphisms are instances of the pseudonaturality invertible modification for the associator.  (Mnemonic for the rectangular invertible modifications: GReen conGRuences.)
      \begin{center}
        \scalebox{0.9}{
        \begin{tikzpicture}[line join=round]
          \filldraw[white,fill=green,fill opacity=0.1](-2.872,1.858)--(-.135,5.617)--(-2.391,3.949)--(-5.127,.19)--cycle;
          \filldraw[white,fill=red,fill opacity=0.1](3.212,-4.9)--(4.306,-3.396)--(4.306,3.532)--(2.391,.901)--(2.391,-3.949)--cycle;
          \filldraw[white,fill=green,fill opacity=0.1](-4.306,-3.532)--(3.212,-4.9)--(2.391,-3.949)--(-5.127,-2.581)--cycle;
          \filldraw[white,fill=red,fill opacity=0.1](-2.872,1.858)--(2.391,.901)--(4.306,3.532)--(2.872,5.07)--(-.135,5.617)--cycle;
          \filldraw[white,fill=red,fill opacity=0.1](-5.127,-2.581)--(-5.127,.19)--(-2.872,1.858)--(2.391,.901)--(2.391,-3.949)--cycle;
          \begin{scope}[font=\fontsize{8}{8}\selectfont]
            \node (A) at (-2.391,3.949) {$(A(B(CD)))E$};
            \node (B) at (-5.127,.19) {$(A((BC)D))E$}
              edge [->] node [l, above left] {$(Aa)E$} (A);
            \node (C) at (-5.127,-2.581) {$((A(BC))D)E$}
              edge [->] node [l, left] {$aE$} (B);
            \node (D) at (-4.306,-3.532) {$(((AB)C)D)E$}
              edge [->] node [l, left] {$(aD)E$} (C);
            \node (E) at (3.212,-4.9) {$((AB)C)(DE)$}
              edge [<-] node [l, below] {$a$} (D);
            \node (F) at (4.306,-3.396) {$(AB)(C(DE))$}
              edge [<-] node [l, below right] {$a$} (E);
            \node (G) at (4.306,3.532) {$A(B(C(DE)))$}
              edge [<-] node [l, right] {$a$} (F);
            \node (H) at (2.872,5.07) {$A(B((CD)E))$}
              edge [->] node [l, above right] {$A(Ba)$} (G);
            \node (I) at (-.135,5.617) {$A((B(CD))E)$}
              edge [->] node [l, above] {$Aa$} (H)
              edge [<-] node [l, above left] {$a$} (A);
            \node (J) at (-2.872,1.858) {$A(((BC)D)E)$}
              edge [->] node [l, below right] {$A(aE)$} (I)
              edge [<-] node [l, below right] {$a$} (B);
            \node (K) at (2.391,-3.949) {$(A(BC))(DE)$}
              edge [<-] node [l, left] {$a(DE)$} (E)
              edge [<-] node [l, above] {$a$} (C);
            \node (L) at (2.391,.901) {$A((BC)(DE))$}
              edge [<-] node [l, left] {$a$} (K)
              edge [<-] node [l, above] {$Aa$} (J)
              edge [->] node [l, above left] {$Aa$} (G);
            \node at (-1,-1) {\tikz\node [rotate=-45] {$\Rightarrow$};};
            \node at (-.5,-1) {$\pi$};
            \node at (1,3) {\tikz\node [rotate=-45] {$\Rightarrow$};};
            \node at (1.5,3) {$A \pi$};
            \node at (3,-1.5) {\tikz\node [rotate=-45] {$\Rightarrow$};};
            \node at (3.5,-1.5) {$\pi$};
            \node at (-1,-3.7) {$\cong$};
            \node at (-2.5,3) {$\cong$};
          \end{scope}
        \end{tikzpicture}}
        \\
        \scalebox{3}{
        \begin{tikzpicture}
          \node at (0,0) [rotate=90] {=};
        \end{tikzpicture}
        }
        \\
        \scalebox{0.9}{
        \begin{tikzpicture}[line join=round]
          \filldraw[white,fill=red,fill opacity=0.1](-4.306,-3.532)--(-2.391,-.901)--(-2.391,3.949)--(-5.127,.19)--(-5.127,-2.581)--cycle;
          \filldraw[white,fill=red,fill opacity=0.1](-4.306,-3.532)--(-2.391,-.901)--(2.872,-1.858)--(4.306,-3.396)--(3.212,-4.9)--cycle;
          \filldraw[white,fill=red,fill opacity=0.1](2.872,-1.858)--(2.872,5.07)--(-.135,5.617)--(-2.391,3.949)--(-2.391,-.901)--cycle;
          \filldraw[white,fill=green,fill opacity=0.1](4.306,-3.396)--(4.306,3.532)--(2.872,5.07)--(2.872,-1.858)--cycle;
          \begin{scope}[font=\fontsize{8}{8}\selectfont]
            \node (A) at (-2.391,3.949) {$(A(B(CD)))E$};
            \node (B) at (-5.127,.19) {$(A((BC)D))E$}
              edge [->] node [l, above left] {$(Aa)E$} (A);
            \node (C) at (-5.127,-2.581) {$((A(BC))D)E$}
              edge [->] node [l, left] {$aE$} (B);
            \node (D) at (-4.306,-3.532) {$(((AB)C)D)E$}
              edge [->] node [l, left] {$(aD)E$} (C);
            \node (E) at (3.212,-4.9) {$((AB)C)(DE)$}
              edge [<-] node [l, below] {$a$} (D);
            \node (F) at (4.306,-3.396) {$(AB)(C(DE))$}
              edge [<-] node [l, below right] {$a$} (E);
            \node (G) at (4.306,3.532) {$A(B(C(DE)))$}
              edge [<-] node [l, right] {$a$} (F);
            \node (H) at (2.872,5.07) {$A(B((CD)E))$}
              edge [->] node [l, above right] {$A(Ba)$} (G);
            \node (I) at (-.135,5.617) {$A((B(CD))E)$}
              edge [->] node [l, above] {$Aa$} (H)
              edge [<-] node [l, above left] {$a$} (A);
            \node (M) at (-2.391,-.901) {$((AB)(CD))E$}
              edge [<-] node [l, right] {$aE$} (D)
              edge [->] node [l, right] {$aE$} (A);
            \node (N) at (2.872,-1.858) {$(AB)((CD)E)$}
              edge [<-] node [l, above] {$a$} (M)
              edge [->] node [l, left] {$a$} (H)
              edge [->] node [l, left] {$(AB)a$} (F);
            \node at (-4,-.5) {$\Rightarrow \pi E$};
            \node at (0,-3) {\tikz\node [rotate=-90] {$\Rightarrow$};};
            \node at (0.5,-3) {$\pi$};
            \node at (0,2) {\tikz\node [rotate=-45] {$\Rightarrow$};};
            \node at (0.5,2) {$\pi$};
            \node at (3.5,1) {$\cong$};
          \end{scope}
        \end{tikzpicture}}
      \end{center}
  
  \end{enumerate}

\item Just as in any monoid there is an identity element 1, in every monoidal bicategory there is a {\bf monoidal unit} object $I$. Associated to the monoidal unit are a series of morphisms that express how to ``cancel'' the unit in a product. Each morphism of dimension $n>0$ has two Stasheff polytopes of dimension $n-1$ as ``subcells'', one for parenthesizing $n+1$ objects and the other for parenthesizing the $n$ objects left over after cancellation.  There are $n+1$ ways to insert $I$ into $n$ objects, so there are $n+1$ morphisms of dimension $n$.
  \begin{enumerate}
    \item There is one monoidal unit object $I$.
    \item There are two {\bf unitor} adjoint equivalences $l$ and $r$ that are pseudonatural in $A$.  The Stasheff polytopes for two objects and for one object are both points, so the unitors are line segments joining them.
      \[l\maps I A \rightarrow A, \qquad r\maps A I \rightarrow A.\]
    \item There are three {\bf 2-unitor} invertible modifications $\lambda, \mu,$ and $\rho$. The Stasheff polytope for three objects is a line segment and the Stasheff polytope for two objects is a point, so these modifications are triangles.  (Mnemonic: Umber Unitor.)
      \begin{center}
        \scalebox{0.9}
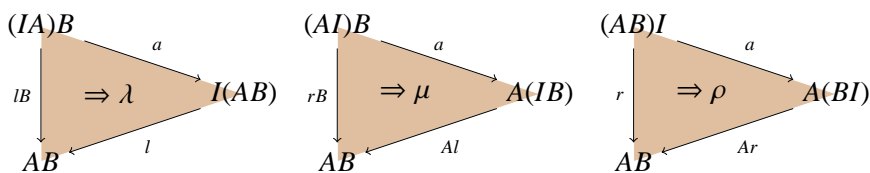
{
        \begin{tikzpicture}
          \filldraw[white,fill=brown,fill opacity=0.5](0,2)--(3,1)--(0,0)--cycle;
          \node (IAB1) at (0,2) {$(I A) B$};
          \node (IAB2) at (3,1) {$I (A  B)$}
            edge [<-] node [l, above right] {$a$} (IAB1);
          \node (AB) at (0,0) {$A  B$}
            edge [<-] node [l, below right] {$l$} (IAB2)
            edge [<-] node [l, left] {$l  B$} (IAB1);
          \node at (1,1) {$\Rightarrow \lambda$};
        \end{tikzpicture}
        \begin{tikzpicture}
          \filldraw[white,fill=brown,fill opacity=0.5](0,2)--(3,1)--(0,0)--cycle;
          \node (AB) at (0,0) {$A  B$};
          \node (AIB2) at (3,1) {$A (I B)$}
            edge [->] node [l, below right] {$A  l$} (AB);
          \node (AIB1) at (0,2) {$(A I)  B$}
            edge [->] node [l, above] {$a$} (AIB2)
            edge [->] node [l, left] {$r  B$} (AB);
          \node at (1,1) {$\Rightarrow \mu$};
        \end{tikzpicture}
        \begin{tikzpicture}
          \filldraw[white,fill=brown,fill opacity=0.5](0,2)--(3,1)--(0,0)--cycle;
          \node (ABI1) at (0,2) {$(A  B)  I$};
          \node (ABI2) at (3,1) {$A  (B  I)$}
            edge [<-] node [l, above right] {$a$} (ABI1);
          \node (AB) at (0,0) {$A  B$}
            edge [<-] node [l, below right] {$A  r$} (ABI2)
            edge [<-] node [l, left] {$r$} (ABI1);
          \node at (1,1) {$\Rightarrow \rho$};
        \end{tikzpicture}
        }
      \end{center}
    \item There are four equations of modifications.  The Stasheff polytope for four objects is a pentagon and the Stasheff polytope for three objects is a line segment, so these equations are irregular prisms with seven vertices.


    \begin{center}
      \scalebox{0.9}
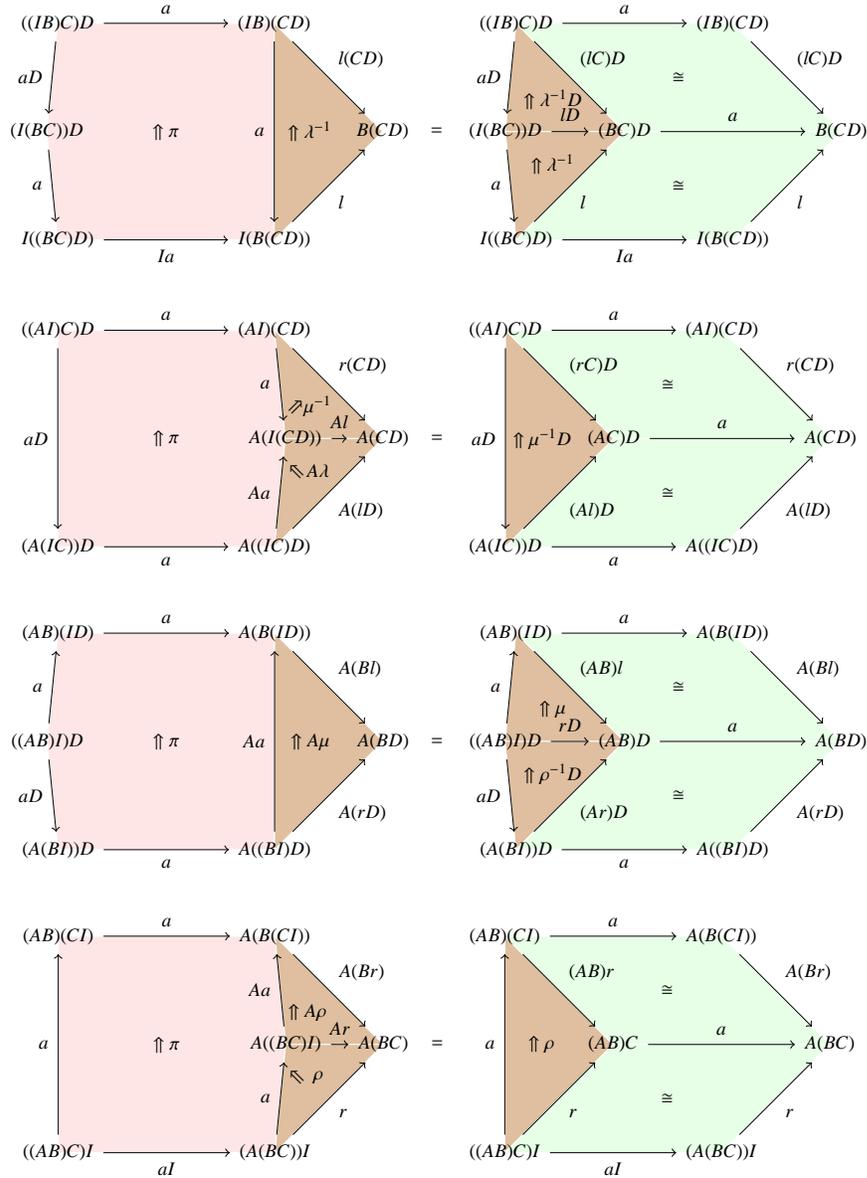
{
      \begin{tikzpicture}[scale=1.6]
        \filldraw[white,fill=red,fill opacity=0.1](-1,1)--(1,1)--(1,-1)--(-1,-1)--(-1.1,0)--cycle;
        \filldraw[white,fill=brown,fill opacity=0.5](1,1)--(2,0)--(1,-1)--cycle;
        \begin{scope}[font=\fontsize{8}{8}\selectfont]
          \node (A) at (-1.1,0) {$(I(BC))D$};
          \node (B) at (-1,1) {$((IB)C)D$}
            edge [->] node [l, left] {$aD$} (A);
          \node (C) at (1,1) {$(IB)(CD)$}
            edge [<-] node [l, above] {$a$} (B);
          \node (D) at (1,-1) {$I(B(CD))$}
            edge [<-] node [l, left] {$a$} (C);
          \node (E) at (-1,-1) {$I((BC)D)$}
            edge [->] node [l, below] {$Ia$} (D)
            edge [<-] node [l, left] {$a$} (A);
          \node (F) at (2,0) {$B(CD)$}
            edge [<-] node [l, above right] {$l(CD)$} (C)
            edge [<-] node [l, below right] {$l$} (D);
          \node at (0,0) {$\Uparrow \pi$};
          \node at (1.33,0) {$\Uparrow \lambda^{-1}$};
          \node at (2.5, 0) {=};
        \end{scope}
      \end{tikzpicture}
      \begin{tikzpicture}[scale=1.6]
        \filldraw[white,fill=brown,fill opacity=0.5](-1,1)--(0,0)--(-1.1,0)--cycle;
        \filldraw[white,fill=brown,fill opacity=0.5](-1,-1)--(0,0)--(-1.1,0)--cycle;
        \filldraw[white,fill=green,fill opacity=0.1](-1,1)--(0,0)--(2,0)--(1,1)--cycle;
        \filldraw[white,fill=green,fill opacity=0.1](-1,-1)--(0,0)--(2,0)--(1,-1)--cycle;
        \begin{scope}[font=\fontsize{8}{8}\selectfont]
          \node (A) at (-1.1,0) {$(I(BC))D$};
          \node (B) at (-1,1) {$((IB)C)D$}
            edge [->] node [l, left] {$aD$} (A);
          \node (C) at (1,1) {$(IB)(CD)$}
            edge [<-] node [l, above] {$a$} (B);
          \node (D) at (1,-1) {$I(B(CD))$};
          \node (E) at (-1,-1) {$I((BC)D)$}
            edge [->] node [l, below] {$Ia$} (D)
            edge [<-] node [l, left] {$a$} (A);
          \node (F) at (2,0) {$B(CD)$}
            edge [<-] node [l, above right] {$(lC)D$} (C)
            edge [<-] node [l, below right] {$l$} (D);
          \node (G) at (0,0) {$(BC)D$}
            edge [<-] node [l, above] {$lD$} (A)
            edge [<-] node [l, above right] {$(lC)D$} (B)
            edge [<-] node [l, below right] {$l$} (E)
            edge [->] node [l, above] {$a$} (F);
          \node at (-.66,.3) {$\Uparrow \lambda^{-1} D$};
          \node at (-.66,-.3) {$\Uparrow \lambda^{-1}$};
          \node at (.5,.5) {$\cong$};
          \node at (.5,-.5) {$\cong$};
        \end{scope}
      \end{tikzpicture}
      }
    \end{center}
    \begin{center}
      \scalebox{0.9}{
      \begin{tikzpicture}[scale=1.6]
        \filldraw[white,fill=red,fill opacity=0.1](-1,-1)--(-1,1)--(1,1)--(1.1,0)--(1,-1)--cycle;
        \filldraw[white,fill=brown,fill opacity=0.5](1,1)--(1.1,0)--(2,0)--cycle;
        \filldraw[white,fill=brown,fill opacity=0.5](1.1,0)--(1,-1)--(2,0)--cycle;
        \begin{scope}[font=\fontsize{8}{8}\selectfont]
          \node (A) at (-1,1) {$((AI)C)D$};
          \node (B) at (1,1) {$(AI)(CD)$}
            edge [<-] node [l,above] {$a$} (A);
          \node (C) at (1.1,0) {$A(I(CD))$}
            edge [<-] node [l,left] {$a$} (B);
          \node (D) at (1,-1) {$A((IC)D)$}
            edge [->] node [l,left] {$Aa$} (C);
          \node (E) at (-1,-1) {$(A(IC))D$}
            edge [->] node [l,below] {$a$} (D)
            edge [<-] node [l,left] {$aD$} (A);
          \node (F) at (2,0) {$A(CD)$}
            edge [<-] node [l,above right] {$r(CD)$} (B)
            edge [<-] node [l,above] {$Al$} (C)
            edge [<-] node [l,below right] {$A(lD)$} (D);
          \node at (0,0) {$\Uparrow \pi$};
          \node at (1.2,.3) {\tikz\node [rotate=45] {$\Rightarrow$};};
          \node at (1.4,.3) {$\mu^{-1}$};
          \node at (1.2,-.3) {\tikz\node [rotate=135] {$\Rightarrow$};};
          \node at (1.4,-.3) {$A\lambda$};
          \node at (2.5,0) {=};
        \end{scope}
      \end{tikzpicture}
      \begin{tikzpicture}[scale=1.6]
        \filldraw[white,fill=brown,fill opacity=0.5](-1,1)--(0,0)--(-1,-1)--cycle;
        \filldraw[white,fill=green,fill opacity=0.1](-1,1)--(0,0)--(2,0)--(1,1)--cycle;
        \filldraw[white,fill=green,fill opacity=0.1](0,0)--(2,0)--(1,-1)--(-1,-1)--cycle;
        \begin{scope}[font=\fontsize{8}{8}\selectfont]
          \node (A) at (-1,1) {$((AI)C)D$};
          \node (B) at (1,1) {$(AI)(CD)$}
            edge [<-] node [l,above] {$a$} (A);
          \node (D) at (1,-1) {$A((IC)D)$};
          \node (E) at (-1,-1) {$(A(IC))D$}
            edge [->] node [l,below] {$a$} (D)
            edge [<-] node [l,left] {$aD$} (A);
          \node (F) at (2,0) {$A(CD)$}
            edge [<-] node [l,above right] {$r(CD)$} (B)
            edge [<-] node [l,below right] {$A(lD)$} (D);
          \node (G) at (0,0) {$(AC)D$}
            edge [<-] node [l,above right] {$(rC)D$} (A)
            edge [<-] node [l,below right] {$(Al)D$} (E)
            edge [->] node [l, above] {$a$} (F);
          \node at (-.66,0) {$\Uparrow \mu^{-1}D$};
          \node at (.5,.5) {$\cong$};
          \node at (.5,-.5) {$\cong$};
        \end{scope}
      \end{tikzpicture}
      }
    \end{center}
    \begin{center}
      \scalebox{0.9}{
      \begin{tikzpicture}[scale=1.6]
        \filldraw[white,fill=red,fill opacity=0.1](-1,1)--(1,1)--(1,-1)--(-1,-1)--(-1.1,0)--cycle;
        \filldraw[white,fill=brown,fill opacity=0.5](1,1)--(2,0)--(1,-1)--cycle;
        \begin{scope}[font=\fontsize{8}{8}\selectfont]
          \node (A) at (-1.1,0) {$((AB)I)D$};
          \node (B) at (-1,1) {$(AB)(ID)$}
            edge [<-] node [l, left] {$a$} (A);
          \node (C) at (1,1) {$A(B(ID))$}
            edge [<-] node [l, above] {$a$} (B);
          \node (D) at (1,-1) {$A((BI)D)$}
            edge [->] node [l, left] {$Aa$} (C);
          \node (E) at (-1,-1) {$(A(BI))D$}
            edge [->] node [l, below] {$a$} (D)
            edge [<-] node [l, left] {$aD$} (A);
          \node (F) at (2,0) {$A(BD)$}
            edge [<-] node [l, above right] {$A(Bl)$} (C)
            edge [<-] node [l, below right] {$A(rD)$} (D);
          \node at (0,0) {$\Uparrow \pi$};
          \node at (1.33,0) {$\Uparrow A\mu$};
          \node at (2.5, 0) {=};
        \end{scope}
      \end{tikzpicture}
      \begin{tikzpicture}[scale=1.6]
        \filldraw[white,fill=brown,fill opacity=0.5](-1,1)--(0,0)--(-1.1,0)--cycle;
        \filldraw[white,fill=brown,fill opacity=0.5](-1,-1)--(0,0)--(-1.1,0)--cycle;
        \filldraw[white,fill=green,fill opacity=0.1](-1,1)--(0,0)--(2,0)--(1,1)--cycle;
        \filldraw[white,fill=green,fill opacity=0.1](-1,-1)--(0,0)--(2,0)--(1,-1)--cycle;
        \begin{scope}[font=\fontsize{8}{8}\selectfont]
          \node (A) at (-1.1,0) {$((AB)I)D$};
          \node (B) at (-1,1) {$(AB)(ID)$}
            edge [<-] node [l, left] {$a$} (A);
          \node (C) at (1,1) {$A(B(ID))$}
            edge [<-] node [l, above] {$a$} (B);
          \node (D) at (1,-1) {$A((BI)D)$};
          \node (E) at (-1,-1) {$(A(BI))D$}
            edge [->] node [l, below] {$a$} (D)
            edge [<-] node [l, left] {$aD$} (A);
          \node (F) at (2,0) {$A(BD)$}
            edge [<-] node [l, above right] {$A(Bl)$} (C)
            edge [<-] node [l, below right] {$A(rD)$} (D);
          \node (G) at (0,0) {$(AB)D$}
            edge [<-] node [l, above] {$rD$} (A)
            edge [<-] node [l, above right] {$(AB)l$} (B)
            edge [<-] node [l, below right] {$(Ar)D$} (E)
            edge [->] node [l, above] {$a$} (F);
          \node at (-.66,.3) {$\Uparrow \mu$};
          \node at (-.66,-.3) {$\Uparrow \rho^{-1}D$};
          \node at (.5,.5) {$\cong$};
          \node at (.5,-.5) {$\cong$};
        \end{scope}
      \end{tikzpicture}
      }
    \end{center}
    \begin{center}
      \scalebox{0.9}{
      \begin{tikzpicture}[scale=1.6]
        \filldraw[white,fill=red,fill opacity=0.1](-1,-1)--(-1,1)--(1,1)--(1.1,0)--(1,-1)--cycle;
        \filldraw[white,fill=brown,fill opacity=0.5](1,1)--(1.1,0)--(2,0)--cycle;
        \filldraw[white,fill=brown,fill opacity=0.5](1.1,0)--(1,-1)--(2,0)--cycle;
        \begin{scope}[font=\fontsize{8}{8}\selectfont]
          \node (A) at (-1,1) {$(AB)(CI)$};
          \node (B) at (1,1) {$A(B(CI))$}
            edge [<-] node [l,above] {$a$} (A);
          \node (C) at (1.1,0) {$A((BC)I)$}
            edge [->] node [l,left] {$Aa$} (B);
          \node (D) at (1,-1) {$(A(BC))I$}
            edge [->] node [l,left] {$a$} (C);
          \node (E) at (-1,-1) {$((AB)C)I$}
            edge [->] node [l,below] {$aI$} (D)
            edge [->] node [l,left] {$a$} (A);
          \node (F) at (2,0) {$A(BC)$}
            edge [<-] node [l,above right] {$A(Br)$} (B)
            edge [<-] node [l,above] {$Ar$} (C)
            edge [<-] node [l,below right] {$r$} (D);
          \node at (0,0) {$\Uparrow \pi$};
          \node at (1.3,.3) {$\Uparrow A\rho$};
          \node at (1.2,-.3) {\tikz\node [rotate=135] {$\Rightarrow$};};
          \node at (1.4,-.3) {$\rho$};
          \node at (2.5,0) {=};
        \end{scope}
      \end{tikzpicture}
      \begin{tikzpicture}[scale=1.6]
        \filldraw[white,fill=brown,fill opacity=0.5](-1,1)--(0,0)--(-1,-1)--cycle;
        \filldraw[white,fill=green,fill opacity=0.1](-1,1)--(0,0)--(2,0)--(1,1)--cycle;
        \filldraw[white,fill=green,fill opacity=0.1](0,0)--(2,0)--(1,-1)--(-1,-1)--cycle;
        \begin{scope}[font=\fontsize{8}{8}\selectfont]
          \node (A) at (-1,1) {$(AB)(CI)$};
          \node (B) at (1,1) {$A(B(CI))$}
            edge [<-] node [l,above] {$a$} (A);
          \node (D) at (1,-1) {$(A(BC))I$};
          \node (E) at (-1,-1) {$((AB)C)I$}
            edge [->] node [l,below] {$aI$} (D)
            edge [->] node [l,left] {$a$} (A);
          \node (F) at (2,0) {$A(BC)$}
            edge [<-] node [l,above right] {$A(Br)$} (B)
            edge [<-] node [l,below right] {$r$} (D);
          \node (G) at (0,0) {$(AB)C$}
            edge [<-] node [l,above right] {$(AB)r$} (A)
            edge [<-] node [l,below right] {$r$} (E)
            edge [->] node [l, above] {$a$} (F);
          \node at (-.66,0) {$\Uparrow \rho$};
          \node at (.5,.5) {$\cong$};
          \node at (.5,-.5) {$\cong$};
        \end{scope}
      \end{tikzpicture}
      }
    \end{center}
  \end{enumerate}
\end{itemize}
\end{defn}

\begin{defn}
  A {\bf braided} monoidal bicategory $\C$ is a monoidal bicategory in which objects can be moved past each other.  A braided monoidal bicategory consists of the following:
  \begin{itemize}
    \item A monoidal bicategory $\C$;
    \item A series of morphisms for ``shuffling''.  
    \begin{defn}
      A {\bf shuffle} of a list $\mathcal{A} = (A_1, \ldots, A_n)$ into a list $\mathcal{B} = (B_1, \ldots, B_k)$ inserts each element of
$\mathcal{A}$ into $\mathcal{B}$ such that if $0 < i < j < n+1$ then
$A_i$ appears to the left of $A_j$.
    \end{defn}

    An ``$(n,k)$-shuffle polytope'' is an $n$-dimensional polytope whose vertices are all the different shuffles of an $n$-element list into a
$k$-element list; there are ${n+k \choose k}$ ways to do this.  General shuffle polytopes were defined by Kapronov and Voevodsky \cite{KV94}. As with the Stasheff polytopes, we have morphisms of the same shape as $(n, k)$-shuffle polytopes with directed edges and faces.
    \begin{itemize}
      \item $(n=1,k=1)$: ${{1+1} \choose 1} = 2,$ so this polytope has two vertices, $(A,B)$ and $(B,A)$.  It has a single edge, which we call a ``braiding'', which encodes how $A$ moves past $B$.  It is an adjoint equivalence pseudonatural in $A, B$.
        \[ b\maps AB \to BA \]
      \item $(n=1,k=2)$ and $(n=2,k=1)$: ${{1+2} \choose 1} = {{2+1}
\choose 1} = 3,$ so whenever the associator is the identity---{\em e.g.} in a braided strictly monoidal bicategory---these polytopes are triangles, invertible modifications whose edges are the directed (1,1) polytope, the braiding.  There are two triangles because the braiding in a braided monoidal bicategory is not necessarily symmetric; when it happens to be symmetric, one can be derived from the other.
        
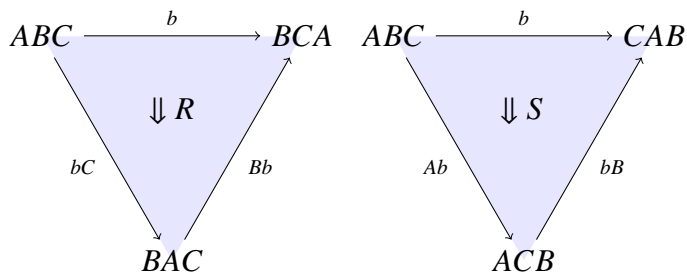
\begin{center}
          \begin{tikzpicture}
            \filldraw[white,fill=blue,fill opacity=0.1](30:2cm)--(150:2cm)--(270:2cm)--cycle;
            \node (A) at (30:2cm) {$BCA$};
            \node (C) at (150:2cm) {$ABC$}
              edge [->] node [l, above] {$b$} (A);
            \node (E) at (270:2cm) {$BAC$}
              edge [<-] node [l, below left] {$bC$} (C)
              edge [->] node [l, below right] {$Bb$} (A);
            \node at (0,0) {$\Downarrow R$};
          \end{tikzpicture}
          \begin{tikzpicture}
            \filldraw[white,fill=blue,fill opacity=0.1](30:2cm)--(150:2cm)--(270:2cm)--cycle;
            \node (A) at (30:2cm) {$CAB$};
            \node (C) at (150:2cm) {$ABC$}
              edge [->] node [l, above] {$b$} (A);
            \node (E) at (270:2cm) {$ACB$}
              edge [<-] node [l, below left] {$Ab$} (C)
              edge [->] node [l, below right] {$bB$} (A);
            \node at (0,0) {$\Downarrow S$};
          \end{tikzpicture}
        \end{center}
        When the associator is not the identity, the triangles' vertices get replaced with associators, effectively truncating them, and we are left with hexagon invertible modifications. (Mnemonic: Blue Braiding.)
        
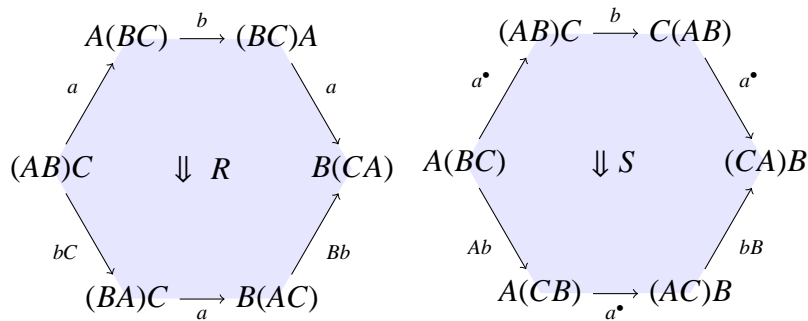
\begin{center}
          \begin{tikzpicture}
            \filldraw[white,fill=blue,fill opacity=0.1](0:2cm)--(60:2cm)--(120:2cm)--(180:2cm)--(240:2cm)--(300:2cm)--cycle;
            \node (A) at (  0:2cm) {$B(CA)$};
            \node (B) at ( 60:2cm) {$(BC)A$}
              edge [->] node [l, above right] {$a$} (A);
            \node (C) at (120:2cm) {$A(BC)$}
              edge [->] node [l, above] {$b$} (B);
            \node (D) at (180:2cm) {$(AB)C$}
              edge [->] node [l, above left] {$a$} (C);
            \node (E) at (240:2cm) {$(BA)C$}
              edge [<-] node [l, below left] {$bC$} (D);
            \node (F) at (300:2cm) {$B(AC)$}
              edge [<-] node [l, below] {$a$} (E)
              edge [->] node [l, below right] {$Bb$} (A);
            \node at (-0.25,0) {\tikz\node [rotate=-90] {$\Rightarrow$};};
            \node at (0.25,0) {$R$};
          \end{tikzpicture}
          \begin{tikzpicture}
            \filldraw[white,fill=blue,fill opacity=0.1](0:2cm)--(60:2cm)--(120:2cm)--(180:2cm)--(240:2cm)--(300:2cm)--cycle;
            \node (A) at (  0:2cm) {$(CA)B$};
            \node (B) at ( 60:2cm) {$C(AB)$}
              edge [->] node [l, above right] {$a^\bullet$} (A);
            \node (C) at (120:2cm) {$(AB)C$}
              edge [->] node [l, above] {$b$} (B);
            \node (D) at (180:2cm) {$A(BC)$}
              edge [->] node [l, above left] {$a^\bullet$} (C);
            \node (E) at (240:2cm) {$A(CB)$}
              edge [<-] node [l, below left] {$Ab$} (D);
            \node (F) at (300:2cm) {$(AC)B$}
              edge [<-] node [l, below] {$a^\bullet$} (E)
              edge [->] node [l, below right] {$bB$} (A);
            \node at (0,0) {$\Downarrow S$};
          \end{tikzpicture}
        \end{center}
      \item $(n=3,k=1)$ and $(n=1,k=3)$: ${{3+1} \choose 1} = {{1+3} \choose 1} = 4,$ so in a braided strictly monoidal bicategory, these polytopes are tetrahedra whose faces are the $(2,1)$ polytope.  As with $R$ and $S$ above, there are two polytopes because the braiding is not necessarily symmetric.
        
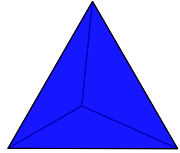
\begin{center}
          \begin{tikzpicture}[line join=round,scale=.2]
            \filldraw[fill=blue,fill opacity=0.7](-.551,-.318)--(.172,6.631)--(-5.449,-3.146)--cycle;
            \filldraw[fill=blue,fill opacity=0.7](-.551,-.318)--(5.828,-3.167)--(-5.449,-3.146)--cycle;
            \filldraw[fill=blue,fill opacity=0.7](-.551,-.318)--(.172,6.631)--(5.828,-3.167)--cycle;
            \filldraw[fill=blue,fill opacity=0.7](-5.449,-3.146)--(.172,6.631)--(5.828,-3.167)--cycle;
          \end{tikzpicture}
        \end{center}
        Again, when the associator is not the identity, the vertices get truncated, this time being replaced by pentagonators; as a side-effect, four of the six edges are also beveled.
        \begin{center}
          \begin{tikzpicture}[line join=round,scale=.2]
            \filldraw[fill=green,fill opacity=0.7](3.11,-1.226)--(1.675,-.585)--(1.112,-1.562)--(2.548,-2.203)--cycle;
            \filldraw[fill=blue,fill opacity=0.7](3.281,.425)--(2.009,2.629)--(-.144,3.591)--(-.235,2.722)--(1.675,-.585)--(3.11,-1.226)--cycle;
            \filldraw[fill=blue,fill opacity=0.7](-3.234,-1.214)--(-2.132,-.578)--(-.235,2.722)--(-.144,3.591)--(-1.798,2.636)--(-3.062,.436)--cycle;
            \filldraw[fill=red,fill opacity=0.7](-1.566,-1.557)--(-2.132,-.578)--(-.235,2.722)--(1.675,-.585)--(1.112,-1.562)--cycle;
            \filldraw[fill=blue,fill opacity=0.7](2.548,-2.203)--(1.112,-1.562)--(-1.566,-1.557)--(-2.668,-2.194)--(-.515,-3.155)--(.894,-3.158)--cycle;
            \filldraw[fill=green,fill opacity=0.7](-3.234,-1.214)--(-2.132,-.578)--(-1.566,-1.557)--(-2.668,-2.194)--cycle;
            \filldraw[fill=red,fill opacity=0.7](1.519,2.346)--(-1.16,2.351)--(-1.798,2.636)--(-.144,3.591)--(2.009,2.629)--cycle;
            \filldraw[fill=red,fill opacity=0.7](3.281,.425)--(2.791,.142)--(.894,-3.158)--(2.548,-2.203)--(3.11,-1.226)--cycle;
            \filldraw[fill=red,fill opacity=0.7](-3.234,-1.214)--(-2.668,-2.194)--(-.515,-3.155)--(-2.425,.152)--(-3.062,.436)--cycle;
            \filldraw[fill=green,fill opacity=0.7](3.281,.425)--(2.009,2.629)--(1.519,2.346)--(2.791,.142)--cycle;
            \filldraw[fill=green,fill opacity=0.7](-3.062,.436)--(-1.798,2.636)--(-1.16,2.351)--(-2.425,.152)--cycle;
            \filldraw[fill=blue,fill opacity=0.7](2.791,.142)--(1.519,2.346)--(-1.16,2.351)--(-2.425,.152)--(-.515,-3.155)--(.894,-3.158)--cycle;
          \end{tikzpicture}
        \end{center}
        This equation governs shuffling one object $A$ into three objects $B, C, D$: \\
        \begin{center}
          \begin{tikzpicture}[line join=round,scale=1.25]
            \filldraw[white,fill=red,fill opacity=0.1](1.519,2.346)--(-1.16,2.351)--(-1.798,2.636)--(-.144,3.591)--(2.009,2.629)--cycle;
            \filldraw[white,fill=red,fill opacity=0.1](3.281,.425)--(2.791,.142)--(.894,-3.158)--(2.548,-2.203)--(3.11,-1.226)--cycle;
            \filldraw[white,fill=red,fill opacity=0.1](-3.234,-1.214)--(-2.668,-2.194)--(-.515,-3.155)--(-2.425,.152)--(-3.062,.436)--cycle;
            \filldraw[white,fill=green,fill opacity=0.1](3.281,.425)--(2.009,2.629)--(1.519,2.346)--(2.791,.142)--cycle;
            \filldraw[white,fill=green,fill opacity=0.1](-3.062,.436)--(-1.798,2.636)--(-1.16,2.351)--(-2.425,.152)--cycle;
            \filldraw[white,fill=blue,fill opacity=0.1](2.791,.142)--(1.519,2.346)--(-1.16,2.351)--(-2.425,.152)--(-.515,-3.155)--(.894,-3.158)--cycle;
            \begin{scope}[font=\fontsize{8}{8}\selectfont]
              \node (A) at (-.144,3.591) {};
              \node[anchor=south] at (A) {$(A(BC))D$};
              \node (B) at (2.009,2.629) {}
                edge [<-] node [l, above right] {$a$} (A);
              \node[anchor=west] at (B) {$A((BC)D)$};
              \node (C) at (3.281,.425) {}
                edge [<-] node [l, above right] {$b$} (B);
              \node[anchor=west] at (C) {$((BC)D)A$};
              \node (D) at (3.11,-1.226) {}
                edge [<-] node [l, right] {$a$} (C);
              \node[anchor=west] at (D) {$(BC)(DA)$};
              \node (E) at (2.548,-2.203) {}
                edge [<-] node [l, right] {$a$} (D);
              \node[anchor=north west] at (E) {$B(C(DA))$};
              \node (F) at (.894,-3.158) {}
                edge [->] node [l, below right] {$Ba$} (E);
              \node[anchor=north] at (F) {$B((CD)A)$};
              \node (G) at (-.515,-3.155) {}
                edge [->] node [l, above] {$Bb$} (F);
              \node[anchor=north] at (G) {$B(A(CD))$};
              \node (H) at (-2.668,-2.194) {}
                edge [->] node [l, below left] {$Ba$} (G);
              \node[anchor=north east] at (H) {$B((AC)D)$};
              \node (I) at (-3.234,-1.214) {}
                edge [->] node [l, left] {$a$} (H);
              \node[anchor=east] at (I) {$(B(AC))D$};
              \node (J) at (-3.062,.436) {}
                edge [->] node [l, left] {$aD$} (I);
              \node[anchor=east] at (J) {$((BA)C)D$};
              \node (K) at (-1.798,2.636) {}
                edge [->] node [l, above left] {$(bC)D$} (J)
                edge [->] node [l, above left] {$aD$} (A);
              \node[anchor=east] at (K) {$((AB)C)D$};
              \node (L) at (-1.16,2.351) {}
                edge [<-] node [l, below left] {$a$} (K);
              \node[anchor=north west] at (L) {$(AB)(CD)$};
              \node (M) at (1.519,2.346) {}
                edge [<-] node [l, above] {$a$} (L)
                edge [<-] node [l, below, sloped] {$Aa$} (B);
              \node[anchor=north east] at (M) {$A(B(CD))$};
              \node (N) at (2.791,.142) {}
                edge [<-] node [l, left] {$b$} (M)
                edge [->] node [l, left] {$a$} (F)
                edge [<-] node [l, above, sloped] {$aA$} (C);
              \node[anchor=east] at (N) {$(B(CD))A$};
              \node (O) at (-2.425,.152) {}
                edge [<-] node [l, below right] {$b(CD)$} (L)
                edge [<-] node [l, above] {$a$} (J)
                edge [->] node [l, right] {$a$} (G);
              \node[anchor=west] at (O) {$(BA)(CD)$};
              \node at (-0.25,3) {\tikz\node [rotate=-90] {$\Rightarrow$};};
              \node at (0,3) {$\pi$};
              \node at (-2,1.5) {$\cong$};
              \node at (2.3,1.5) {$\cong$};
              \node at (-2.25,-1.5) {$\Leftarrow \pi^{-1}$};
              \node at (2.5,-1.5) {$\Leftarrow \pi^{-1}$};
              \node at (-0,0) {\tikz\node [rotate=-135] {$\Rightarrow$};};
              \node at (0.25,0) {$R$};
            \end{scope}
          \end{tikzpicture}
          \\
          \scalebox{3}{
          \begin{tikzpicture}
            \node at (0,0) [rotate=90] {=};
          \end{tikzpicture}
          }
          \\
          \begin{tikzpicture}[line join=round,scale=1.25]
            \filldraw[white,fill=green,fill opacity=0.1](3.11,-1.226)--(1.675,-.585)--(1.112,-1.562)--(2.548,-2.203)--cycle;
            \filldraw[white,fill=blue,fill opacity=0.1](3.281,.425)--(2.009,2.629)--(-.144,3.591)--(-.235,2.722)--(1.675,-.585)--(3.11,-1.226)--cycle;
            \filldraw[white,fill=blue,fill opacity=0.1](-3.234,-1.214)--(-2.132,-.578)--(-.235,2.722)--(-.144,3.591)--(-1.798,2.636)--(-3.062,.436)--cycle;
            \filldraw[white,fill=red,fill opacity=0.1](-1.566,-1.557)--(-2.132,-.578)--(-.235,2.722)--(1.675,-.585)--(1.112,-1.562)--cycle;
            \filldraw[white,fill=blue,fill opacity=0.1](2.548,-2.203)--(1.112,-1.562)--(-1.566,-1.557)--(-2.668,-2.194)--(-.515,-3.155)--(.894,-3.158)--cycle;
            \filldraw[white,fill=green,fill opacity=0.1](-3.234,-1.214)--(-2.132,-.578)--(-1.566,-1.557)--(-2.668,-2.194)--cycle;
            \begin{scope}[font=\fontsize{8}{8}\selectfont]
              \node (A) at (-.144,3.591) {};
              \node[anchor=south] at (A) {$(A(BC))D$};
              \node (B) at (2.009,2.629) {}
                edge [<-] node [l, above right] {$a$} (A);
              \node[anchor=west] at (B) {$A((BC)D)$};
              \node (C) at (3.281,.425) {}
                edge [<-] node [l, above right] {$b$} (B);
              \node[anchor=west] at (C) {$((BC)D)A$};
              \node (D) at (3.11,-1.226) {}
                edge [<-] node [l, right] {$a$} (C);
              \node[anchor=west] at (D) {$(BC)(DA)$};
              \node (E) at (2.548,-2.203) {}
                edge [<-] node [l, right] {$a$} (D);
              \node[anchor=north west] at (E) {$B(C(DA))$};
              \node (F) at (.894,-3.158) {}
                edge [->] node [l, below right] {$Ba$} (E);
              \node[anchor=north] at (F) {$B((CD)A)$};
              \node (G) at (-.515,-3.155) {}
                edge [->] node [l, above] {$Bb$} (F);
              \node[anchor=north] at (G) {$B(A(CD))$};
              \node (H) at (-2.668,-2.194) {}
                edge [->] node [l, below left] {$Ba$} (G);
              \node[anchor=north east] at (H) {$B((AC)D)$};
              \node (I) at (-3.234,-1.214) {}
                edge [->] node [l, left] {$a$} (H);
              \node[anchor=east] at (I) {$(B(AC))D$};
              \node (J) at (-3.062,.436) {}
                edge [->] node [l, left] {$aD$} (I);
              \node[anchor=east] at (J) {$((BA)C)D$};
              \node (K) at (-1.798,2.636) {}
                edge [->] node [l, above left] {$(bC)D$} (J)
                edge [->] node [l, above left] {$aD$} (A);
              \node[anchor=east] at (K) {$((AB)C)D$};
              \node (P) at (-2.132,-.578) {}
                edge [<-] node [l, above, sloped] {$(Bb)D$} (I);
              \node[anchor=west] at (P) {$(B(CA))D$};
              \node (Q) at (-.235,2.722) {}
                edge [->] node [l, below right] {$aD$} (P)
                edge [<-] node [l, left] {$bD$} (A);
              \node[anchor=west] at (Q) {$((BC)A)D$};
              \node (R) at (1.675,-.585) {}
                edge [<-] node [l, below left] {$a$} (Q)
                edge [->] node [l, above, sloped] {$(BC)b$} (D);
              \node[anchor=east] at (R) {$(BC)(AD)$};
              \node (S) at (1.112,-1.562) {}
                edge [<-] node [l, below right] {$a$} (R)
                edge [->] node [l, above, sloped] {$B(Cb)$} (E);
              \node[anchor=south east] at (S) {$B(C(AD))$};
              \node (T) at (-1.566,-1.557) {}
                edge [->] node [l, below] {$Ba$} (S)
                edge [<-] node [l, below left] {$a$} (P)
                edge [<-] node [l, above, sloped] {$B(bD)$} (H);
              \node[anchor=south west] at (T) {$B((CA)D)$};
              \node at (-0.25,0) {\tikz\node [rotate=-135] {$\Rightarrow$};};
              \node at (0.25,0) {$\pi^{-1}$};
              \node at (2,1) {$\Leftarrow R$};
              \node at (-2.25,1) {\tikz\node [rotate=-135] {$\Rightarrow$};};
              \node at (-2,1) {$RD$};
              \node at (2.25,-1.25) {$\cong$};
              \node at (-2.5,-1.25) {$\cong$};
              \node at (-0.5,-2.5) {\tikz\node [rotate=-90] {$\Rightarrow$};};
              \node at (0,-2.5) {$BR^{-1}$};
            \end{scope}
          \end{tikzpicture}
        \end{center}
        This equation governs shuffling one object $D$ into three objects $A,B,C$: \\
        \begin{center}
          \begin{tikzpicture}[line join=round,scale=1.25]
            \filldraw[white,fill=red,fill opacity=0.1](1.519,2.346)--(-1.16,2.351)--(-1.798,2.636)--(-.144,3.591)--(2.009,2.629)--cycle;
            \filldraw[white,fill=red,fill opacity=0.1](3.281,.425)--(2.791,.142)--(.894,-3.158)--(2.548,-2.203)--(3.11,-1.226)--cycle;
            \filldraw[white,fill=red,fill opacity=0.1](-3.234,-1.214)--(-2.668,-2.194)--(-.515,-3.155)--(-2.425,.152)--(-3.062,.436)--cycle;
            \filldraw[white,fill=green,fill opacity=0.1](3.281,.425)--(2.009,2.629)--(1.519,2.346)--(2.791,.142)--cycle;
            \filldraw[white,fill=green,fill opacity=0.1](-3.062,.436)--(-1.798,2.636)--(-1.16,2.351)--(-2.425,.152)--cycle;
            \filldraw[white,fill=blue,fill opacity=0.1](2.791,.142)--(1.519,2.346)--(-1.16,2.351)--(-2.425,.152)--(-.515,-3.155)--(.894,-3.158)--cycle;
            \begin{scope}[font=\fontsize{8}{8}\selectfont]
              \node (A) at (-.144,3.591) {};
              \node[anchor=south] at (A) {$A((BC)D)$};
              \node (B) at (2.009,2.629) {}
                edge [<-] node [l, above right] {$a^\bullet$} (A);
              \node[anchor=west] at (B) {$(A(BC))D$};
              \node (C) at (3.281,.425) {}
                edge [<-] node [l, above right] {$b$} (B);
              \node[anchor=west] at (C) {$D(A(BC))$};
              \node (D) at (3.11,-1.226) {}
                edge [<-] node [l, right] {$a^\bullet$} (C);
              \node[anchor=west] at (D) {$(DA)(BC)$};
              \node (E) at (2.548,-2.203) {}
                edge [<-] node [l, right] {$a^\bullet$} (D);
              \node[anchor=north west] at (E) {$((DA)B)C$};
              \node (F) at (.894,-3.158) {}
                edge [->] node [l, below right] {$a^\bullet C$} (E);
              \node[anchor=north] at (F) {$(D(AB))C$};
              \node (G) at (-.515,-3.155) {}
                edge [->] node [l, above] {$bC$} (F);
              \node[anchor=north] at (G) {$((AB)D)C$};
              \node (H) at (-2.668,-2.194) {}
                edge [->] node [l, below left] {$a^\bullet C$} (G);
              \node[anchor=north east] at (H) {$(A(BD))C$};
              \node (I) at (-3.234,-1.214) {}
                edge [->] node [l, left] {$a^\bullet$} (H);
              \node[anchor=east] at (I) {$A((BD)C)$};
              \node (J) at (-3.062,.436) {}
                edge [->] node [l, left] {$Aa^\bullet$} (I);
              \node[anchor=east] at (J) {$A(B(DC))$};
              \node (K) at (-1.798,2.636) {}
                edge [->] node [l, above left] {$a(Bb)$} (J)
                edge [->] node [l, above left] {$Aa^\bullet$} (A);
              \node[anchor=east] at (K) {$A(B(CD))$};
              \node (L) at (-1.16,2.351) {}
                edge [<-] node [l, below left] {$a^\bullet$} (K);
              \node[anchor=north west] at (L) {$(AB)(CD)$};
              \node (M) at (1.519,2.346) {}
                edge [<-] node [l, above] {$a^\bullet$} (L)
                edge [<-] node [l, below, sloped] {$a^\bullet D$} (B);
              \node[anchor=north east] at (M) {$((AB)C)D$};
              \node (N) at (2.791,.142) {}
                edge [<-] node [l, left] {$b$} (M)
                edge [->] node [l, left] {$a^\bullet$} (F)
                edge [<-] node [l, above, sloped] {$Da^\bullet$} (C);
              \node[anchor=east] at (N) {$D((AB)C)$};
              \node (O) at (-2.425,.152) {}
                edge [<-] node [l, below right] {$(AB)b$} (L)
                edge [<-] node [l, above] {$a^\bullet$} (J)
                edge [->] node [l, right] {$a^\bullet$} (G);
              \node[anchor=west] at (O) {$(AB)(DC)$};
              \node at (-0.25,3) {\tikz\node [rotate=-90] {$\Rightarrow$};};
              \node at (0,3) {$\pi^\bullet$};
              \node at (-2,1.5) {$\cong$};
              \node at (2.3,1.5) {$\cong$};
              \node at (-2.25,-1.5) {$\Leftarrow \pi^{\bullet-1}$};
              \node at (2.5,-1.5) {$\Leftarrow \pi^{\bullet-1}$};
              \node at (-0,0) {\tikz\node [rotate=-135] {$\Rightarrow$};};
              \node at (0.25,0) {$S$};
            \end{scope}
          \end{tikzpicture}
          \\
          \scalebox{3}{
          \begin{tikzpicture}
            \node at (0,0) [rotate=90] {=};
          \end{tikzpicture}
          }
          \\
          \begin{tikzpicture}[line join=round,scale=1.25]
            \filldraw[white,fill=green,fill opacity=0.1](3.11,-1.226)--(1.675,-.585)--(1.112,-1.562)--(2.548,-2.203)--cycle;
            \filldraw[white,fill=blue,fill opacity=0.1](3.281,.425)--(2.009,2.629)--(-.144,3.591)--(-.235,2.722)--(1.675,-.585)--(3.11,-1.226)--cycle;
            \filldraw[white,fill=blue,fill opacity=0.1](-3.234,-1.214)--(-2.132,-.578)--(-.235,2.722)--(-.144,3.591)--(-1.798,2.636)--(-3.062,.436)--cycle;
            \filldraw[white,fill=red,fill opacity=0.1](-1.566,-1.557)--(-2.132,-.578)--(-.235,2.722)--(1.675,-.585)--(1.112,-1.562)--cycle;
            \filldraw[white,fill=blue,fill opacity=0.1](2.548,-2.203)--(1.112,-1.562)--(-1.566,-1.557)--(-2.668,-2.194)--(-.515,-3.155)--(.894,-3.158)--cycle;
            \filldraw[white,fill=green,fill opacity=0.1](-3.234,-1.214)--(-2.132,-.578)--(-1.566,-1.557)--(-2.668,-2.194)--cycle;
            \begin{scope}[font=\fontsize{8}{8}\selectfont]
              \node (A) at (-.144,3.591) {};
              \node[anchor=south] at (A) {$A((BC)D)$};
              \node (B) at (2.009,2.629) {}
                edge [<-] node [l, above right] {$a^\bullet$} (A);
              \node[anchor=west] at (B) {$(A(BC))D$};
              \node (C) at (3.281,.425) {}
                edge [<-] node [l, above right] {$b$} (B);
              \node[anchor=west] at (C) {$D(A(BC))$};
              \node (D) at (3.11,-1.226) {}
                edge [<-] node [l, right] {$a^\bullet$} (C);
              \node[anchor=west] at (D) {$(DA)(BC)$};
              \node (E) at (2.548,-2.203) {}
                edge [<-] node [l, right] {$a^\bullet$} (D);
              \node[anchor=north west] at (E) {$((DA)B)C$};
              \node (F) at (.894,-3.158) {}
                edge [->] node [l, below right] {$a^\bullet C$} (E);
              \node[anchor=north] at (F) {$(D(AB))C$};
              \node (G) at (-.515,-3.155) {}
                edge [->] node [l, above] {$bC$} (F);
              \node[anchor=north] at (G) {$((AB)D)C$};
              \node (H) at (-2.668,-2.194) {}
                edge [->] node [l, below left] {$a^\bullet C$} (G);
              \node[anchor=north east] at (H) {$(A(BD))C$};
              \node (I) at (-3.234,-1.214) {}
                edge [->] node [l, left] {$a^\bullet$} (H);
              \node[anchor=east] at (I) {$A((BD)C)$};
              \node (J) at (-3.062,.436) {}
                edge [->] node [l, left] {$Aa^\bullet$} (I);
              \node[anchor=east] at (J) {$A(B(DC))$};
              \node (K) at (-1.798,2.636) {}
                edge [->] node [l, above left] {$A(Bb)$} (J)
                edge [->] node [l, above left] {$Aa^\bullet$} (A);
              \node[anchor=east] at (K) {$A(B(CD))$};
              \node (P) at (-2.132,-.578) {}
                edge [<-] node [l, above, sloped] {$A(bC)$} (I);
              \node[anchor=west] at (P) {$A((DB)C)$};
              \node (Q) at (-.235,2.722) {}
                edge [->] node [l, below right] {$Aa^\bullet$} (P)
                edge [<-] node [l, left] {$Ab$} (A);
              \node[anchor=west] at (Q) {$A(D(BC))$};
              \node (R) at (1.675,-.585) {}
                edge [<-] node [l, below left] {$a^\bullet$} (Q)
                edge [->] node [l, above, sloped] {$b(BC)$} (D);
              \node[anchor=east] at (R) {$(AD)(BC)$};
              \node (S) at (1.112,-1.562) {}
                edge [<-] node [l, below right] {$a^\bullet$} (R)
                edge [->] node [l, above, sloped] {$(bB)C$} (E);
              \node[anchor=south east] at (S) {$((AD)B)C$};
              \node (T) at (-1.566,-1.557) {}
                edge [->] node [l, below] {$a^\bullet C$} (S)
                edge [<-] node [l, below left] {$a^\bullet$} (P)
                edge [<-] node [l, above, sloped] {$(Ab)C$} (H);
              \node[anchor=south west] at (T) {$(A(DB))C$};
              \node at (-0.25,0) {\tikz\node [rotate=-135] {$\Rightarrow$};};
              \node at (0.25,0) {$\pi^{\bullet-1}$};
              \node at (2,1) {$\Leftarrow S$};
              \node at (-2.25,1) {\tikz\node [rotate=-135] {$\Rightarrow$};};
              \node at (-2,1) {$AS$};
              \node at (2.25,-1.25) {$\cong$};
              \node at (-2.5,-1.25) {$\cong$};
              \node at (-0.5,-2.5) {\tikz\node [rotate=-90] {$\Rightarrow$};};
              \node at (0,-2.5) {$S^{-1}C$};
            \end{scope}
          \end{tikzpicture}
        \end{center}
        
      \item $(n=2,k=2)$: ${{2+2} \choose 2} = 6;$ in a braided strictly monoidal bicategory, this polytope is composed mostly of (2,1) triangles, but there is a pair of braidings that commute, so one face is a square.
        
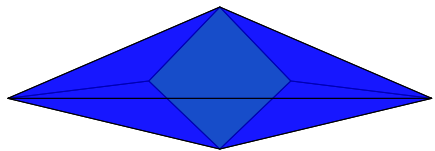
\begin{center}
          \begin{tikzpicture}[line join=round, scale=.5]
            \filldraw[fill=blue,fill opacity=0.7](-5.635,-.347)--(-1.888,.116)--(0,-1.695)--cycle;
            \filldraw[fill=blue,fill opacity=0.7](5.635,-.347)--(0,-1.695)--(1.888,.116)--cycle;
            \filldraw[fill=green,fill opacity=0.7](0,2.083)--(-1.888,.116)--(0,-1.695)--(1.888,.116)--cycle;
            \filldraw[fill=blue,fill opacity=0.7](5.635,-.347)--(-5.635,-.347)--(0,-1.695)--cycle;
            \filldraw[fill=blue,fill opacity=0.7](5.635,-.347)--(1.888,.116)--(0,2.083)--cycle;
            \filldraw[fill=blue,fill opacity=0.7](0,2.083)--(-1.888,.116)--(-5.635,-.347)--cycle;
            \filldraw[fill=blue,fill opacity=0.7](5.635,-.347)--(-5.635,-.347)--(0,2.083)--cycle;
          \end{tikzpicture}
        \end{center}
        When the associator is not the identity, the six vertices get truncated and six of the edges get beveled.
        
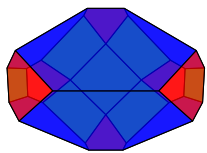
\begin{center}
          \begin{tikzpicture}[line join=round, scale=.5]
            \filldraw[fill=green,fill opacity=0.7](0,-.807)--(-.944,.116)--(-1.678,-.488)--(-.546,-1.564)--cycle;
            \filldraw[fill=green,fill opacity=0.7](.944,.116)--(0,-.807)--(.546,-1.564)--(1.678,-.488)--cycle;
            \filldraw[fill=green,fill opacity=0.7](.944,.116)--(0,-.807)--(-.944,.116)--(0,1.079)--cycle;
            \filldraw[fill=red,fill opacity=0.7](.546,-1.564)--(.453,-1.75)--(-.453,-1.75)--(-.546,-1.564)--(0,-.807)--cycle;
            \filldraw[fill=green,fill opacity=0.7](0,1.079)--(-.944,.116)--(-1.72,.642)--(-.588,1.829)--cycle;
            \filldraw[fill=green,fill opacity=0.7](.944,.116)--(0,1.079)--(.588,1.829)--(1.72,.642)--cycle;
            \filldraw[fill=red,fill opacity=0.7](1.678,-.488)--(2.075,-.32)--(2.109,.417)--(1.72,.642)--(.944,.116)--cycle;
            \filldraw[fill=red,fill opacity=0.7](-1.72,.642)--(-2.109,.417)--(-2.075,-.32)--(-1.678,-.488)--(-.944,.116)--cycle;
            \filldraw[fill=red,fill opacity=0.7](.588,1.829)--(.492,2.021)--(-.492,2.021)--(-.588,1.829)--(0,1.079)--cycle;
            \filldraw[fill=blue,fill opacity=0.7](1.678,-.488)--(2.075,-.32)--(2.573,-.48)--(2.308,-.977)--(.453,-1.75)--(.546,-1.564)--cycle;
            \filldraw[fill=blue,fill opacity=0.7](-.546,-1.564)--(-.453,-1.75)--(-2.308,-.977)--(-2.573,-.48)--(-2.075,-.32)--(-1.678,-.488)--cycle;
            \filldraw[fill=blue,fill opacity=0.7](-.453,-1.75)--(-2.308,-.977)--(-1.412,-.174)--(1.412,-.174)--(2.308,-.977)--(.453,-1.75)--cycle;
            \filldraw[fill=blue,fill opacity=0.7](1.72,.642)--(.588,1.829)--(.492,2.021)--(2.405,.899)--(2.625,.433)--(2.109,.417)--cycle;
            \filldraw[fill=blue,fill opacity=0.7](-.588,1.829)--(-.492,2.021)--(-2.405,.899)--(-2.625,.433)--(-2.109,.417)--(-1.72,.642)--cycle;
            \filldraw[fill=green,fill opacity=0.7](2.075,-.32)--(2.573,-.48)--(2.625,.433)--(2.109,.417)--cycle;
            \filldraw[fill=green,fill opacity=0.7](-2.075,-.32)--(-2.573,-.48)--(-2.625,.433)--(-2.109,.417)--cycle;
            \filldraw[fill=blue,fill opacity=0.7](.492,2.021)--(-.492,2.021)--(-2.405,.899)--(-1.412,-.174)--(1.412,-.174)--(2.405,.899)--cycle;
            \filldraw[fill=red,fill opacity=0.7](2.625,.433)--(2.405,.899)--(1.412,-.174)--(2.308,-.977)--(2.573,-.48)--cycle;
            \filldraw[fill=red,fill opacity=0.7](-2.625,.433)--(-2.405,.899)--(-1.412,-.174)--(-2.308,-.977)--(-2.573,-.48)--cycle;
          \end{tikzpicture}
        \end{center}
        This equation governs shuffling two objects $A, B$ into two objects $C, D$: \\
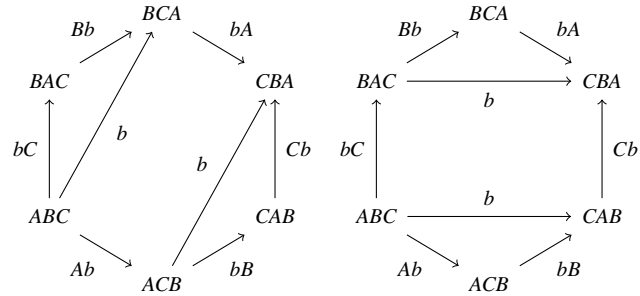
        \begin{center}
          \scalebox{0.8}{
          \begin{tikzpicture}[line join=round, scale=2.5, font=\fontsize{8}{8}\selectfont]
            \filldraw[white,fill=blue,fill opacity=0.1](-.453,-1.75)--(-2.308,-.977)--(-1.412,-.174)--(1.412,-.174)--(2.308,-.977)--(.453,-1.75)--cycle;
            \filldraw[white,fill=blue,fill opacity=0.1](.492,2.021)--(-.492,2.021)--(-2.405,.899)--(-1.412,-.174)--(1.412,-.174)--(2.405,.899)--cycle;
            \filldraw[white,fill=red,fill opacity=0.1](2.625,.433)--(2.405,.899)--(1.412,-.174)--(2.308,-.977)--(2.573,-.48)--cycle;
            \filldraw[white,fill=red,fill opacity=0.1](-2.625,.433)--(-2.405,.899)--(-1.412,-.174)--(-2.308,-.977)--(-2.573,-.48)--cycle;
            \node (A) at (-2.573,-.48) {};
            \node [left] at (A) {$(A(BC))D$};

            \node (B) at (-2.625,.433) {}
              edge [<-] node [l, left] {$a$} (A);
            \node [left] at (B) {$A((BC)D)$};

            \node (C) at (-2.405,.899) {}
              edge [<-] node [l, left] {$Aa$} (B);
            \node [left] at (C) {$A(B(CD))$};

            \node (D) at (-.492,2.021) {}
              edge [<-] node [l, above left] {$Ab$} (C);
            \node [above] at (D) {$A((CD)B)$};

            \node (E) at (.492,2.021) {}
              edge [<-] node [l, above] {$a^\bullet$} (D);
            \node [above] at (E) {$(A(CD))B$};

            \node (F) at (2.405,.899) {}
              edge [<-] node [l, above right] {$bB$} (E);
            \node [right] at (F) {$((CD)A)B$};

            \node (G) at (2.625,.433) {}
              edge [<-] node [l, right] {$aB$} (F);
            \node [right] at (G) {$(C(DA))B$};

            \node (H) at (2.573,-.48) {}
              edge [<-] node [l, right] {$a$} (G);
            \node [right] at (H) {$C((DA)B)$};

            \node (I) at (2.308,-.977) {}
              edge [->] node [l, right] {$Ca^\bullet$} (H);
            \node [right] at (I) {$C(D(AB))$};

            \node (J) at (.453,-1.75) {}
              edge [->] node [l, below right] {$Cb$} (I);
            \node [below] at (J) {$C((AB)D)$};

            \node (K) at (-.453,-1.75) {}
              edge [->] node [l, below] {$a$} (J);
            \node [below] at (K) {$(C(AB))D$};

            \node (L) at (-2.308,-.977) {}
              edge [->] node [l, below left] {$bD$} (K)
              edge [<-] node [l, left] {$a^\bullet D$} (A);
            \node [left] at (L) {$((AB)C)D$};

            \node (M) at (-1.412,-.174) {$(AB)(CD)$}
              edge [<-] node [l, below right] {$a$} (L)
              edge [<-] node [l, above right] {$a^\bullet$} (C);

            \node (N) at (1.412,-.174) {$(CD)(AB)$}
              edge [<-] node [l, above] {$b$} (M)
              edge [->] node [l, below left] {$a$} (I)
              edge [->] node [l, above left] {$a^\bullet$} (F);
            \begin{scope}[font=\fontsize{10}{10}\selectfont]
              \node at (-2.2,0) {$\Uparrow \pi_1$};
              \node at (0,1) {$\Uparrow S$};
              \node at (0,-1) {$\Uparrow R^{-1}$};
              \node at (2.2,0) {$\Uparrow \pi_2$};
            \end{scope}
          \end{tikzpicture}
          }
          \\
          \scalebox{3}{
          \begin{tikzpicture}
            \node at (0,0) [rotate=90] {=};
          \end{tikzpicture}
          }
          \\
          \scalebox{0.8}{
          \begin{tikzpicture}[line join=round, scale=2.5, font=\fontsize{8}{8}\selectfont]
            \filldraw[white,fill=green,fill opacity=0.1](0,-.807)--(-.944,.116)--(-1.678,-.488)--(-.546,-1.564)--cycle;
            \filldraw[white,fill=green,fill opacity=0.1](.944,.116)--(0,-.807)--(.546,-1.564)--(1.678,-.488)--cycle;
            \filldraw[white,fill=green,fill opacity=0.1](.944,.116)--(0,-.807)--(-.944,.116)--(0,1.079)--cycle;
            \filldraw[white,fill=red,fill opacity=0.1](.546,-1.564)--(.453,-1.75)--(-.453,-1.75)--(-.546,-1.564)--(0,-.807)--cycle;
            \filldraw[white,fill=green,fill opacity=0.1](0,1.079)--(-.944,.116)--(-1.72,.642)--(-.588,1.829)--cycle;
            \filldraw[white,fill=green,fill opacity=0.1](.944,.116)--(0,1.079)--(.588,1.829)--(1.72,.642)--cycle;
            \filldraw[white,fill=red,fill opacity=0.1](1.678,-.488)--(2.075,-.32)--(2.109,.417)--(1.72,.642)--(.944,.116)--cycle;
            \filldraw[white,fill=red,fill opacity=0.1](-1.72,.642)--(-2.109,.417)--(-2.075,-.32)--(-1.678,-.488)--(-.944,.116)--cycle;
            \filldraw[white,fill=red,fill opacity=0.1](.588,1.829)--(.492,2.021)--(-.492,2.021)--(-.588,1.829)--(0,1.079)--cycle;
            \filldraw[white,fill=blue,fill opacity=0.1](1.678,-.488)--(2.075,-.32)--(2.573,-.48)--(2.308,-.977)--(.453,-1.75)--(.546,-1.564)--cycle;
            \filldraw[white,fill=blue,fill opacity=0.1](-.546,-1.564)--(-.453,-1.75)--(-2.308,-.977)--(-2.573,-.48)--(-2.075,-.32)--(-1.678,-.488)--cycle;
            \filldraw[white,fill=blue,fill opacity=0.1](1.72,.642)--(.588,1.829)--(.492,2.021)--(2.405,.899)--(2.625,.433)--(2.109,.417)--cycle;
            \filldraw[white,fill=blue,fill opacity=0.1](-.588,1.829)--(-.492,2.021)--(-2.405,.899)--(-2.625,.433)--(-2.109,.417)--(-1.72,.642)--cycle;
            \filldraw[white,fill=green,fill opacity=0.1](2.075,-.32)--(2.573,-.48)--(2.625,.433)--(2.109,.417)--cycle;
            \filldraw[white,fill=green,fill opacity=0.1](-2.075,-.32)--(-2.573,-.48)--(-2.625,.433)--(-2.109,.417)--cycle;

            \node (A) at (-2.573,-.48) {};
            \node [left] at (A) {$(A(BC))D$};

            \node (B) at (-2.625,.433) {}
              edge [<-] node [l, left] {$a$} (A);
            \node [left] at (B) {$A((BC)D)$};

            \node (C) at (-2.405,.899) {}
              edge [<-] node [l, left] {$Aa$} (B);
            \node [left] at (C) {$A(B(CD))$};

            \node (D) at (-.492,2.021) {}
              edge [<-] node [l, above left] {$Ab$} (C);
            \node [above] at (D) {$A((CD)B)$};

            \node (E) at (.492,2.021) {}
              edge [<-] node [l, above] {$a^\bullet$} (D);
            \node [above] at (E) {$(A(CD))B$};

            \node (F) at (2.405,.899) {}
              edge [<-] node [l, above right] {$bB$} (E);
            \node [right] at (F) {$((CD)A)B$};

            \node (G) at (2.625,.433) {}
              edge [<-] node [l, right] {$aB$} (F);
            \node [right] at (G) {$(C(DA))B$};

            \node (H) at (2.573,-.48) {}
              edge [<-] node [l, right] {$a$} (G);
            \node [right] at (H) {$C((DA)B)$};

            \node (I) at (2.308,-.977) {}
              edge [->] node [l, right] {$Ca^\bullet$} (H);
            \node [right] at (I) {$C(D(AB))$};

            \node (J) at (.453,-1.75) {}
              edge [->] node [l, below right] {$Cb$} (I);
            \node [below] at (J) {$C((AB)D)$};

            \node (K) at (-.453,-1.75) {}
              edge [->] node [l, below] {$a$} (J);
            \node [below] at (K) {$(C(AB))D$};

            \node (L) at (-2.308,-.977) {}
              edge [->] node [l, below left] {$bD$} (K)
              edge [<-] node [l, left] {$a^\bullet D$} (A);
            \node [left] at (L) {$((AB)C)D$};

            \node (O) at (-2.075,-.32) {}
              edge [<-] node [l, below, sloped] {$(Ab)D$} (A);
            \node [above right] at (O) {$(A(CB))D$};

            \node (P) at (-2.109,.417) {}
              edge [<-] node [l, left] {$a$} (O)
              edge [<-] node [l, above, sloped] {$A(bD)$} (B);
            \node [below right] at (P) {$A((BC)D)$};

            \node (Q) at (-1.72,.642) {}
              edge [<-] node [l, above, sloped] {$Aa$} (P);
            \node [right] at (Q) {$A(C(BD))$};

            \node (R) at (-.944,.116) {$(AC)(BD)$}
              edge [<-] node [l, above right] {$a^\bullet$} (Q);

            \node (S) at (-1.678,-.488) {}
              edge [->] node [l, below right] {$a$} (R)
              edge [<-] node [l, below, sloped] {$a^\bullet D$} (O);
            \node [right] at (S) {$((AC)B)D$};

            \node (T) at (-.546,-1.564) {}
              edge [<-] node [l, above, sloped] {$(bB)D$} (S)
              edge [<-] node [right] {$a^\bullet D$} (K);
            \node at (T) {$((CA)B)D$};

            \node (U) at (0,-.807) {$(CA)(BD)$}
              edge [<-] node [l, above left] {$a$} (T)
              edge [<-] node [l, above, sloped] {$b(BD)$} (R);

            \node (V) at (.546,-1.564) {}
              edge [<-] node [l, above right] {$a$} (U)
              edge [->] node [left] {$Ca^\bullet$} (J);
            \node at (V) {$C(A(BD))$};

            \node (W) at (1.678,-.488) {}
              edge [<-] node [l, above, sloped] {$C(Ab)$} (V);
            \node [left] at (W) {$C(A(DB))$};

            \node (X) at (2.075,-.32) {}
              edge [<-] node [l, below, sloped] {$Ca^\bullet$} (W)
              edge [->] node [l, below, sloped] {$C(bB)$} (H);
            \node [above left] at (X) {$C((AD)B)$};

            \node (Y) at (2.109,.417) {}
              edge [->] node [l, right] {$a$} (X)
              edge [->] node [l, above, sloped] {$(Cb)B$} (G);
            \node [below left] at (Y) {$(C(AD))B$};

            \node (Z) at (1.72,.642) {}
              edge [->] node [l, above, sloped] {$aB$} (Y);
            \node [left] at (Z) {$((CA)D)B$};

            \node (A') at (.944,.116) {$(CA)(DB)$}
              edge [->] node [l, above left] {$a^\bullet$} (Z)
              edge [->] node [l, below left] {$a$} (W)
              edge [<-] node [l, above, sloped] {$(CA)b$} (U);

            \node (B') at (0,1.079) {$(AC)(DB)$}
              edge [->] node [l, below, sloped] {$b(DB)$} (A')
              edge [<-] node [l, below, sloped] {$(AC)b$} (R);

            \node (C') at (-.588,1.829) {}
              edge [->] node [l, below left] {$a^\bullet$} (B')
              edge [<-] node [l, below, sloped] {$A(Cb)$} (Q)
              edge [<-] node [right] {$Aa$} (D);
            \node at (C') {$A(C(DB))$};

            \node (D') at (.588,1.829) {}
              edge [<-] node [l, below right] {$a^\bullet$} (B')
              edge [->] node [l, below, sloped] {$(bD)B$} (Z)
              edge [->] node [left] {$aB$} (E);
            \node at (D') {$((AC)D)B$};

            \begin{scope} [font=\fontsize{10}{10}\selectfont]
              \node at (-2.4,0) {$\cong$};
              \node at (-1.7,0) {\tikz\node [rotate=90] {$\Rightarrow$};};
              \node at (-1.5,0) {$\pi_1$};
              \node at (-2,-.8) {\tikz\node [rotate=90] {$\Rightarrow$};};
              \node at (-1.8,-.8) {$SD$};
              \node at (-2,.9) {\tikz\node [rotate=90] {$\Rightarrow$};};
              \node at (-1.75,.9) {$AR^{-1}$};
              \node at (-0.1,1.6) {\tikz\node [rotate=90] {$\Rightarrow$};};
              \node at (0.1,1.6) {$\pi_4$};
              \node at (-.75,.75) {$\cong$};
              \node at (-.75,-.75) {$\cong$};
              \node at (0,0) {$\cong$};
              \node at (.75,-.75) {$\cong$};
              \node at (.75,.75) {$\cong$};
              \node at (-0.1,-1.5) {\tikz\node [rotate=90] {$\Rightarrow$};};
              \node at (0.1,-1.5) {$\pi_3$};
              \node at (1.8,.9) {\tikz\node [rotate=90] {$\Rightarrow$};};
              \node at (2.05,.9) {$R^{-1}B$};
              \node at (1.8,-.8) {\tikz\node [rotate=90] {$\Rightarrow$};};
              \node at (2, -.8) {$CS$};
              \node at (1.5,0) {\tikz\node [rotate=90] {$\Rightarrow$};};
              \node at (1.7,0) {$\pi_2$};
              \node at (2.4,0) {$\cong$};
            \end{scope}
          \end{tikzpicture}
          }
        \end{center}
    \end{itemize}
    \item The Breen polytope.  In a braided monoidal category, the Yang-Baxter equations hold; there are two fundamentally distinct proofs of this fact.
    \begin{center}
      \begin{tikzpicture}
        \begin{scope}[scale=.3, font=\fontsize{8}{8}\selectfont]
          \node (ACB) at (5,0) {ACB};
          \node (ABC) at (0,3) {ABC}
            edge [->] node [l, below left] {$Ab$} (ACB);
          \node (CAB) at (10,3) {CAB}
            edge [<-] node [l, below right] {$b B$} (ACB);
          \node (BAC) at (0,9) {BAC}
            edge [<-] node [l, left] {$b C$} (ABC);
          \node (CBA) at (10,9) {CBA}
            edge [<-] node [l, above left] {$b$} (ACB)
            edge [<-] node [l, right] {$Cb$} (CAB);
          \node (BCA) at (5,12) {BCA}
            edge [<-] node [l, above left] {$Bb$} (BAC)
            edge [<-] node [l, below right] {$b$} (ABC)
            edge [->] node [l, above right] {$b A$} (CBA);
        \end{scope}
      \end{tikzpicture}
      \begin{tikzpicture}
        \begin{scope}[scale=.3, font=\fontsize{8}{8}\selectfont]
          \node (ACB) at (5,0) {ACB};
          \node (ABC) at (0,3) {ABC}
            edge [->] node [l, below left] {$Ab$} (ACB);
          \node (CAB) at (10,3) {CAB}
            edge [<-] node [l, above] {$b$} (ABC)
            edge [<-] node [l, below right] {$b B$} (ACB);
          \node (BAC) at (0,9) {BAC}
            edge [<-] node [l, left] {$b C$} (ABC);
          \node (CBA) at (10,9) {CBA}
            edge [<-] node [l, below] {$b$} (BAC)
            edge [<-] node [l, right] {$Cb$} (CAB);
          \node (BCA) at (5,12) {BCA}
            edge [<-] node [l, above left] {$Bb$} (BAC)
            edge [->] node [l, above right] {$b A$} (CBA);
        \end{scope}
      \end{tikzpicture}
    \end{center}
    In a braided strictly monoidal bicategory, the two proofs become the front and back face of another coherence law governing the interaction of the (2,1)-shuffle polytopes; when the associator is nontrivial, the vertices get truncated.  That the coherence law is necessary was something of a surprise: Kapranov and Voevodsky did not include it in their definition of braided semistrict monoidal 2-categories; Breen
\cite{Breen} corrected the definition.  We therefore call the following coherence law the ``Breen polytope''.  In retrospect, we can see that this is the start of a more subtle collection of polytopes relevant to braided monoidal $n$-categories, which can be systematically obtained using Batanin's approach to weak $n$-categories \cite{Batanin}.
    \begin{center}
      \begin{tikzpicture}
        \begin{scope}[scale=.8, font=\fontsize{8}{8}\selectfont]
          \filldraw[white,fill=blue,fill opacity=0.1] (1,2.5)--(0,4)--(0,8)--(1,9.5)--(4,11.5)--(6,11.5)--cycle;
          \filldraw[white,fill=blue,fill opacity=0.1] (9,9.5)--(10,8)--(10,4)--(9,2.5)--(6,0.5)--(4,0.5)--cycle;
          \filldraw[white,fill=green,fill opacity=0.1] (6,11.5)--(9,9.5)--(4,0.5)--(1,2.5)--cycle;
          \node (ABC1) at (1, 2.5) {};
          \node [below left, l] at (ABC1) {$A(BC)$};
          \node (ABC2) at (0, 4) {}
            edge [<-] node [below left, l] {$a^\bullet$} (ABC1);
          \node [left, l] at (ABC2) {$(AB)C$};
          \node (BAC1) at (0, 8) {}
            edge [<-] node [left, l] {$bC$} (ABC2);
          \node [left, l] at (BAC1) {$(BA)C$};
          \node (BAC2) at (1, 9.5) {}
            edge [<-] node [left, l] {$a$} (BAC1);
          \node [above left, l] at (BAC2) {$B(AC)$};
          \node (BCA1) at (4, 11.5) {}
            edge [<-] node [above left, l] {$Bb$} (BAC2);
          \node [above, l] at (BCA1) {$B(CA)$};
          \node (BCA2) at (6, 11.5) {}
            edge [<-] node [above, l] {$a^\bullet$} (BCA1)
            edge [<-] node [below right] {$b$} (ABC1);
          \node [above, l] at (BCA2) {$(BC)A$};
          \node (CBA1) at (9, 9.5) {}
            edge [<-] node [above right, l] {$bA$} (BCA2);
          \node [above right, l] at (CBA1) {$(CB)A$};
          \node (CBA2) at (10, 8) {}
            edge [->] node [right, l] {$a^\bullet$} (CBA1);
          \node [right, l] at (CBA2) {$C(BA)$};
          \node (CAB1) at (10, 4) {}
            edge [->] node [right, l] {$Cb$} (CBA2);
          \node [right, l] at (CAB1) {$C(AB)$};
          \node (CAB2) at (9, 2.5) {}
            edge [->] node [right, l] {$a$} (CAB1);
          \node [below right, l] at (CAB2) {$(CA)B$};
          \node (ACB1) at (6, 0.5) {}
            edge [->] node [below right, l] {$bB$} (CAB2);
          \node [below, l] at (ACB1) {$(AC)B$};
          \node (ACB2) at (4, 0.5) {}
            edge [->] node [below, l] {$a^\bullet$} (ACB1)
            edge [->] node [below right, l] {$b$} (CBA1)
            edge [<-] node [below left, l] {$Ab$} (ABC1);
          \node [below, l] at (ACB2) {$A(CB)$};
          \node at (2, 7.5) {$\Rightarrow R_1^{-1}$};
          \node at (5, 6) {$\cong$};
          \node at (8, 4.5) {$\Rightarrow R_1$};
        \end{scope}
      \end{tikzpicture}
      \\
      \scalebox{3}{
      \begin{tikzpicture}
        \node at (0,0) [rotate=90] {=};
      \end{tikzpicture}
      }
      \\
      \begin{tikzpicture}
        \begin{scope}[scale=.8, font=\fontsize{8}{8}\selectfont]
          \filldraw[white,fill=blue,fill opacity=0.1] (0,8)--(1,9.5)--(4,11.5)--(6,11.5)--(9,9.5)--(10,8)--cycle;
          \filldraw[white,fill=blue,fill opacity=0.1] (10,4)--(9,2.5)--(6,0.5)--(4,0.5)--(1,2.5)--(0,4)--cycle;
          \filldraw[white,fill=green,fill opacity=0.1] (0,4)--(0,8)--(10,8)--(10,4)--cycle;
          \node (ABC1) at (1, 2.5) {};
          \node [below left, l] at (ABC1) {$A(BC)$};
          \node (ABC2) at (0, 4) {}
            edge [<-] node [below left, l] {$a^\bullet$} (ABC1);
          \node [left, l] at (ABC2) {$(AB)C$};
          \node (BAC1) at (0, 8) {}
            edge [<-] node [left, l] {$bC$} (ABC2);
          \node [left, l] at (BAC1) {$(BA)C$};
          \node (BAC2) at (1, 9.5) {}
            edge [<-] node [left, l] {$a$} (BAC1);
          \node [above left, l] at (BAC2) {$B(AC)$};
          \node (BCA1) at (4, 11.5) {}
            edge [<-] node [above left, l] {$Bb$} (BAC2);
          \node [above, l] at (BCA1) {$B(CA)$};
          \node (BCA2) at (6, 11.5) {}
            edge [<-] node [above, l] {$a^\bullet$} (BCA1);
          \node [above, l] at (BCA2) {$(BC)A$};
          \node (CBA1) at (9, 9.5) {}
            edge [<-] node [above right, l] {$bA$} (BCA2);
          \node [above right, l] at (CBA1) {$(CB)A$};
          \node (CBA2) at (10, 8) {}
            edge [->] node [right, l] {$a^\bullet$} (CBA1)
            edge [<-] node [below, l] {$b$} (BAC1);
          \node [right, l] at (CBA2) {$C(BA)$};
          \node (CAB1) at (10, 4) {}
            edge [->] node [right, l] {$Cb$} (CBA2)
            edge [<-] node [below, l] {$b$} (ABC2);
          \node [right, l] at (CAB1) {$C(AB)$};
          \node (CAB2) at (9, 2.5) {}
            edge [->] node [right, l] {$a$} (CAB1);
          \node [below right, l] at (CAB2) {$(CA)B$};
          \node (ACB1) at (6, 0.5) {}
            edge [->] node [below right, l] {$bB$} (CAB2);
          \node [below, l] at (ACB1) {$(AC)B$};
          \node (ACB2) at (4, 0.5) {}
            edge [->] node [below, l] {$a^\bullet$} (ACB1)
            edge [<-] node [below left, l] {$Ab$} (ABC1);
          \node [below, l] at (ACB2) {$A(CB)$};
          \node at (5, 10) {$\Downarrow S_1^{-1}$};
          \node at (5, 6) {$\cong$};
          \node at (5, 2) {$\Downarrow S_1$};
        \end{scope}
      \end{tikzpicture}
    \end{center}
  \end{itemize}
\end{defn}

\begin{defn}
  A {\bf sylleptic} monoidal bicategory $\C$ is a braided monoidal bicategory equipped with \\
  \begin{itemize} 
    \item an invertible modification called the syllepsis, (Mnemonic: Salmon Syllepsis)
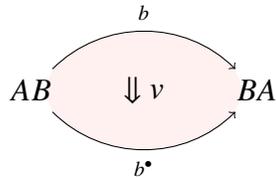
    \begin{center}
      \begin{tikzpicture}
        \filldraw[white,fill=salmon,fill opacity=0.1] (1.5,0) ellipse (37pt and 22pt);
        \node (AB) at (0,0) {$AB$};
        \node (BA) at (3,0) {$BA$}
          edge [<-, out=135, in=45] node [above,l] {$b$} (AB)
          edge [<-, out=-135, in=-45] node [below,l] {$b^\bullet$} (AB);
        \node at (1.5,0) {$\Downarrow v$};
      \end{tikzpicture}
    \end{center}
  \end{itemize}
  subject to the following axioms.
  \begin{itemize} 
    \item This equation governs the interaction of the syllepsis with the $(n=1,k=2)$ braiding:
    \begin{center}
      \scalebox{0.9}
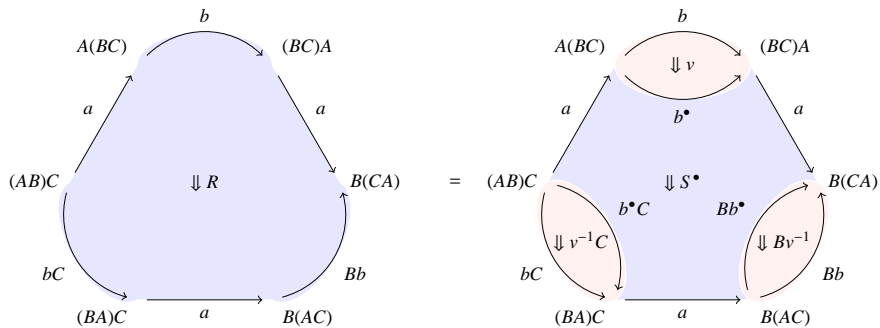
{
      \begin{tikzpicture}
        \begin{scope}[font=\fontsize{8}{8}\selectfont]
          \filldraw [white,fill=blue,fill opacity=0.1] (0:2cm)--(60:2cm) arc (0:180:1cm and 5mm)--(180:2cm) [rotate=120] arc(0:180:1cm and 5mm)--(180:2cm) [rotate=120] arc(0:180:1cm and 5mm)--cycle;
          \node (A) at (  0:2cm) {};
          \node at (A) [right] {$B(CA)$};
          \node (B) at ( 60:2cm) {}
            edge [->] node [l, above right] {$a$} (A);
          \node at (B) [above right] {$(BC)A$};
          \node (C) at (120:2cm) {}
            edge [->, out=45, in=135] node [l, above] {$b$} (B);
          \node at (C) [above left] {$A(BC)$};
          \node (D) at (180:2cm) {}
            edge [->] node [l, above left] {$a$} (C);
          \node at (D) [left] {$(AB)C$};
          \node (E) at (240:2cm) {}
            edge [<-, out=165, in=255] node [l, below left] {$bC$} (D);
          \node at (E) [below left] {$(BA)C$};
          \node (F) at (300:2cm) {}
            edge [<-] node [l, below] {$a$} (E)
            edge [->, out=15, in=285] node [l, below right] {$Bb$} (A);
          \node at (F) [below right] {$B(AC)$};
          \node at (0,0) {$\Downarrow R$};
          \node at (0:3.7cm) {=};
        \end{scope}
      \end{tikzpicture}
      \begin{tikzpicture}
        \begin{scope}[font=\fontsize{8}{8}\selectfont]
          \filldraw [white,fill=blue,fill opacity=0.1] (0:2cm)--(60:2cm) arc(0:-180:1cm and 5mm)--(180:2cm) [rotate=120] arc(0:-180:1cm and 5mm)--(180:2cm) [rotate=120] arc(0:-180:1cm and 5mm)--cycle;
          \filldraw [white,fill=salmon, fill opacity=0.1] (60:2cm) arc(0:360:1cm and 5mm);
          \filldraw [white,fill=salmon, fill opacity=0.1] (180:2cm) [rotate=120] arc(0:360:1cm and 5mm);
          \filldraw [white,fill=salmon, fill opacity=0.1] (300:2cm) [rotate=240] arc(0:360:1cm and 5mm);
          \node (A) at (  0:2cm) {};
          \node at (A) [right] {$B(CA)$};
          \node (B) at ( 60:2cm) {}
            edge [->] node [l, above right] {$a$} (A);
          \node at (B) [above right] {$(BC)A$};
          \node (C) at (120:2cm) {}
            edge [->, out=45, in=135] node [l, above] {$b$} (B)
            edge [->, out=-45, in=-135] node [l, below] {$b^\bullet$} (B);
          \node at (C) [above left] {$A(BC)$};
          \node (D) at (180:2cm) {}
            edge [->] node [l, above left] {$a$} (C);
          \node at (D) [left] {$(AB)C$};
          \node (E) at (240:2cm) {}
            edge [<-, out=165, in=255] node [l, below left] {$bC$} (D)
            edge [<-, out=75, in=-15] node [l, above right] {$b^\bullet C$} (D);
          \node at (E) [below left] {$(BA)C$};
          \node (F) at (300:2cm) {}
            edge [<-] node [l, below] {$a$} (E)
            edge [->, out=15, in=285] node [l, below right] {$Bb$} (A)
            edge [->, out=105, in=195] node [l, above left] {$Bb^\bullet$} (A);
          \node at (F) [below right] {$B(AC)$};
          \node at (0,0) {$\Downarrow S^\bullet$};
          \node at (90:1.732cm) {$\Downarrow v$};
          \node at (210:1.732cm) {$\Downarrow v^{-1}C$};
          \node at (330:1.732cm) {$\Downarrow Bv^{-1}$};
        \end{scope}
      \end{tikzpicture}
      }
    \end{center}
    \item This equation governs the interaction of the syllepsis with the $(n=2,k=1)$ braiding:
    \begin{center}
      \scalebox{0.9}
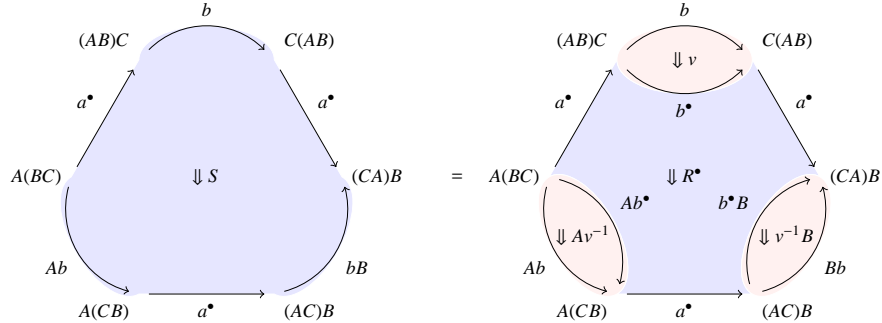
{
      \begin{tikzpicture}
        \begin{scope}[font=\fontsize{8}{8}\selectfont]
          \filldraw [white,fill=blue,fill opacity=0.1] (0:2cm)--(60:2cm) arc(0:180:1cm and 5mm)--(180:2cm) [rotate=120] arc(0:180:1cm and 5mm)--(180:2cm) [rotate=120] arc(0:180:1cm and 5mm)--cycle;
          \node (A) at (  0:2cm) {};
          \node at (A) [right] {$(CA)B$};
          \node (B) at ( 60:2cm) {}
            edge [->] node [l, above right] {$a^\bullet$} (A);
          \node at (B) [above right] {$C(AB)$};
          \node (C) at (120:2cm) {}
            edge [->, out=45, in=135] node [l, above] {$b$} (B);
          \node at (C) [above left] {$(AB)C$};
          \node (D) at (180:2cm) {}
            edge [->] node [l, above left] {$a^\bullet$} (C);
          \node at (D) [left] {$A(BC)$};
          \node (E) at (240:2cm) {}
            edge [<-, out=165, in=255] node [l, below left] {$Ab$} (D);
          \node at (E) [below left] {$A(CB)$};
          \node (F) at (300:2cm) {}
            edge [<-] node [l, below] {$a^\bullet$} (E)
            edge [->, out=15, in=285] node [l, below right] {$bB$} (A);
          \node at (F) [below right] {$(AC)B$};
          \node at (0,0) {$\Downarrow S$};
          \node at (0:3.7cm) {=};
        \end{scope}
      \end{tikzpicture}
      \begin{tikzpicture}
        \begin{scope}[font=\fontsize{8}{8}\selectfont]
          \filldraw [white,fill=blue,fill opacity=0.1] (0:2cm)--(60:2cm) arc(0:-180:1cm and 5mm)--(180:2cm) [rotate=120] arc(0:-180:1cm and 5mm)--(180:2cm) [rotate=120] arc(0:-180:1cm and 5mm)--cycle;
          \filldraw [white,fill=salmon, fill opacity=0.1] (60:2cm) arc(0:360:1cm and 5mm);
          \filldraw [white,fill=salmon, fill opacity=0.1] (180:2cm) [rotate=120] arc(0:360:1cm and 5mm);
          \filldraw [white,fill=salmon, fill opacity=0.1] (300:2cm) [rotate=240] arc(0:360:1cm and 5mm);
          \node (A) at (  0:2cm) {};
          \node at (A) [right] {$(CA)B$};
          \node (B) at ( 60:2cm) {}
            edge [->] node [l, above right] {$a^\bullet$} (A);
          \node at (B) [above right] {$C(AB)$};
          \node (C) at (120:2cm) {}
            edge [->, out=45, in=135] node [l, above] {$b$} (B)
            edge [->, out=-45, in=-135] node [l, below] {$b^\bullet$} (B);
          \node at (C) [above left] {$(AB)C$};
          \node (D) at (180:2cm) {}
            edge [->] node [l, above left] {$a^\bullet$} (C);
          \node at (D) [left] {$A(BC)$};
          \node (E) at (240:2cm) {}
            edge [<-, out=165, in=255] node [l, below left] {$Ab$} (D)
            edge [<-, out=75, in=-15] node [l, above right] {$Ab^\bullet$} (D);
          \node at (E) [below left] {$A(CB)$};
          \node (F) at (300:2cm) {}
            edge [<-] node [l, below] {$a^\bullet$} (E)
            edge [->, out=15, in=285] node [l, below right] {$Bb$} (A)
            edge [->, out=105, in=195] node [l, above left] {$b^\bullet B$} (A);
          \node at (F) [below right] {$(AC)B$};
          \node at (0,0) {$\Downarrow R^\bullet$};
          \node at (90:1.732cm) {$\Downarrow v$};
          \node at (210:1.732cm) {$\Downarrow Av^{-1}$};
          \node at (330:1.732cm) {$\Downarrow v^{-1}B$};
        \end{scope}
      \end{tikzpicture}
      }
    \end{center}
  \end{itemize}
\end{defn}

\begin{defn}
  A {\bf symmetric} monoidal bicategory is a sylleptic monoidal bicategory subject to the following axiom, where the unlabeled green cells are identities: \\
  \begin{itemize}
    \item for all objects $A$ and $B$ of $\C$, the following equation holds: \\
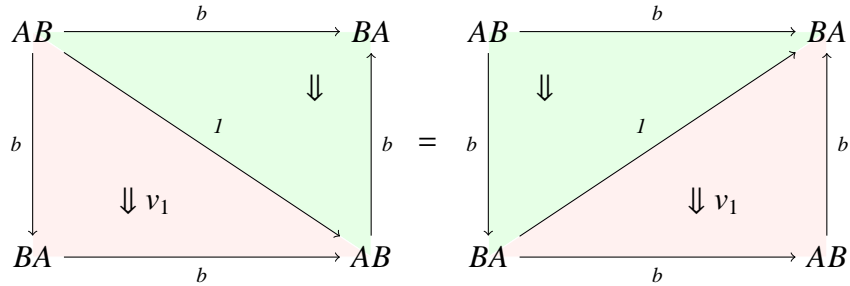
    \begin{center}
      \begin{tikzpicture}[scale=.75]
        \filldraw[white,fill=green,fill opacity=0.1](0,4)--(6,4)--(6,0)--cycle;
        \filldraw[white,fill=salmon,fill opacity=0.1](0,4)--(0,0)--(6,0)--cycle;
        \node (AB1) at (0,4) {$A  B$};
        \node (BA1) at (6,4) {$B  A$}
          edge [<-] node [l, above] {$b$} (AB1);
        \node (BA2) at (0,0) {$B  A$}
          edge [<-] node [l, left] {$b$} (AB1);
        \node (AB2) at (6,0) {$A  B$}
          edge [->] node [l, right] {$b$} (BA1)
          edge [<-] node [l, above right] {1} (AB1)
          edge [<-] node [l, below] {$b$} (BA2);
        \node (v) at (2,1) {$\Downarrow v_1$};
        \node (1) at (5,3) {$\Downarrow$};
        \node (=) at (7,2) {=};
      \end{tikzpicture}
      \begin{tikzpicture}[scale=.75]
        \filldraw[white,fill=green,fill opacity=0.1](0,4)--(6,4)--(0,0)--cycle;
        \filldraw[white,fill=salmon,fill opacity=0.1](0,0)--(6,0)--(6,4)--cycle;
        \node (AB1) at (0,4) {$A  B$};
        \node (BA1) at (6,4) {$B  A$}
          edge [<-] node [l, above] {$b$} (AB1);
        \node (BA2) at (0,0) {$B  A$}
          edge [->] node [l, above left] {1} (BA1)
          edge [<-] node [l, left] {$b$} (AB1);
        \node (AB2) at (6,0) {$A  B$}
          edge [->] node [l, right] {$b$} (BA1)
          edge [<-] node [l, below] {$b$} (BA2);
        \node (v) at (4,1) {$\Downarrow v_1$};
        \node (1) at (1,3) {$\Downarrow$};
      \end{tikzpicture}
    \end{center}
  \end{itemize}
\end{defn}  

\begin{defn}
  Given two bicategories $\mathcal{J}, \mathcal{K}$, two functors 
${L\maps\mathcal{J}\to \mathcal{K}}$ and ${R\maps\mathcal{K} \to \mathcal{J}}$ are
{\bf pseudoadjoint} if for all $A \in \mathcal{J}, B \in \mathcal{K}$
the categories Hom${}_\mathcal{K}(LA, B)$ and Hom${}_\mathcal{J}(A, RB)$
are adjoint equivalent pseudonaturally in $A$ and $B$.
\end{defn}

Symmetric monoidal closed bicategories satisfy the obvious weakening of the definition of symmetric monoidal closed categories:
\begin{defn}
  A symmetric monoidal {\bf closed} bicategory is one in which for every object
$A$, the functor $- \tensor A$ has a right pseudoadjoint $A \lhom -.$
\end{defn}

Similarly, compact closed bicategories weaken the notion of duality from compact closed categories.  In the following definition, we abstract the notion of pseudoadjointness from functors between bicategories to arbitrary objects of a bicategory.
\begin{defn}
  A {\bf compact closed} bicategory is a symmetric monoidal bicategory in which every object has a pseudoadjoint.
\end{defn}
This means that every object $A$ is equipped with a (weak) {\bf dual}, an object $A^*$ equipped with two 1-morphisms
  \def\cap[#1,#2,#3]{
    \coordinate (Astar) at ($ (#1,#2) + (1,0) $);
    \coordinate (top) at ($ (#1,#2) + (0.5,1) $);
    \coordinate (A) at (#1,#2);
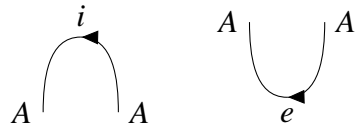
    \begin{scope}[decoration={markings,mark=at position 0.99 with {\arrow{triangle 45};}}]
      \draw[postaction={decorate}] (Astar) to [out=90,in=0] (top);
    \end{scope}
    \draw (top) to [out=180,in=90] (A);
    \node [left] at (A) {$#3$};
    \node [right] at (Astar) {$#3$};
    \node [above] at (top) {$i$};
  }
  \def\cup[#1,#2,#3]{
    \coordinate (Astar) at ($ (#1,#2) + (1,1) $);
    \coordinate (top) at ($ (#1,#2) + (0.5,0) $);
    \coordinate (A) at ($ (#1,#2) + (0,1) $);
    \begin{scope}[decoration={markings,mark=at position 0.99 with {\arrow{triangle 45};}}]
      \draw[postaction={decorate}] (Astar) to [out=-90,in=0] (top);
    \end{scope}
    \draw (top) to [out=180,in=-90] (A);
    \node [left] at (A) {$#3$};
    \node [right] at (Astar) {$#3$};
    \node [below] at (top) {$e$};
  }
  \def\id[#1,#2,#3,#4,#5]{
    \coordinate (A1) at (#1,#2);
    \coordinate (A2) at (#3,#4);
    \begin{scope}[decoration={markings,mark=at position 0.5 with {\arrow{triangle 45};}}]
      \draw[postaction={decorate}] (A1) to (A2);
      \node [right] at ($ (A1)!.5!(A2) $) {$#5$};
    \end{scope}
  }
  \[ i_A\maps  I \to A  A^* \quad \quad e_A\maps  A^*  A \to I \]
  \begin{center}
    \begin{tikzpicture}
      \cap[0,0,A];
    \end{tikzpicture}
    $\quad$
    \begin{tikzpicture}
      \cup[0,0,A];
    \end{tikzpicture}
  \end{center}
and two ``zig-zag'' 2-isomorphisms (Mnemonic: Yellow Yanking or Xanthic Zig-zag)
  \[ \zeta_A\maps  A \Rightarrow (A  e_A) \circ (i_A  A) \]
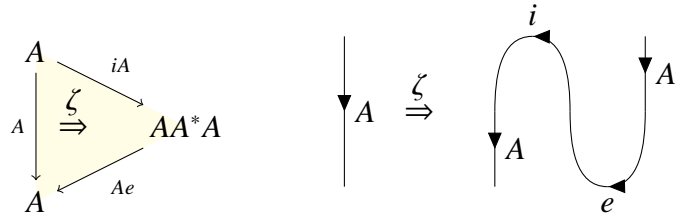
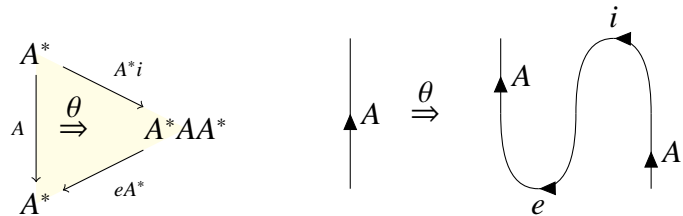
  \begin{center}
    \begin{tikzpicture}
      \node (A) at (0,0) {$A$};
      \node (B) at (2,1) {$AA^*A$}
        edge [->] node [below right, l] {$Ae$} (A);
      \node (C) at (0,2) {$A$}
        edge [->] node [left, l] {$A$} (A)
        edge [->] node [above right, l] {$iA$} (B);
      \node at (.5,1) [above] {$\zeta$};
      \node at (.5,1) {$\Rightarrow$};
      \begin{scope}[on background layer]
        \fill [white,fill=yellow,fill opacity=0.1] (A.center) to (B.center) to (C.center) to (A.center);
      \end{scope}
    \end{tikzpicture}
    $\quad\quad\quad$
    \begin{tikzpicture}
      \id[-2,2,-2,0,A];
      \node at (-1,1) {$\Rightarrow$};
      \node [above] at (-1,1) {$\zeta$};
      \cap[0,1,];
      \id[2,2,2,1,A];
      \id[0,1,0,0,A];
      \cup[1,0,];
    \end{tikzpicture}
  \end{center}
  \[ \theta_A\maps A^* \Rightarrow (e_A  A^*) \circ (A^*  i_A) \]
  \begin{center}
    \begin{tikzpicture}
      \node (A) at (0,0) {$A^*$};
      \node (B) at (2,1) {$A^*AA^*$}
        edge [->] node [below right, l] {$eA^*$} (A);
      \node (C) at (0,2) {$A^*$}
        edge [->] node [left, l] {$A$} (A)
        edge [->] node [above right, l] {$A^*i$} (B);
      \node at (.5,1) [above] {$\theta$};
      \node at (.5,1) {$\Rightarrow$};
      \begin{scope}[on background layer]
        \fill [white,fill=yellow,fill opacity=0.1] (A.center) to (B.center) to (C.center) to (A.center);
      \end{scope}
    \end{tikzpicture}
    $\quad\quad\quad$
    \begin{tikzpicture}
      \id[-2,0,-2,2,A];
      \node at (-1,1) {$\Rightarrow$};
      \node [above] at (-1,1) {$\theta$};
      \id[0,1,0,2,A];
      \cap[1,1,];
      \cup[0,0,];
      \id[2,0,2,1,A];
    \end{tikzpicture}
  \end{center}
such that the following ``swallowtail equation'' holds:
  
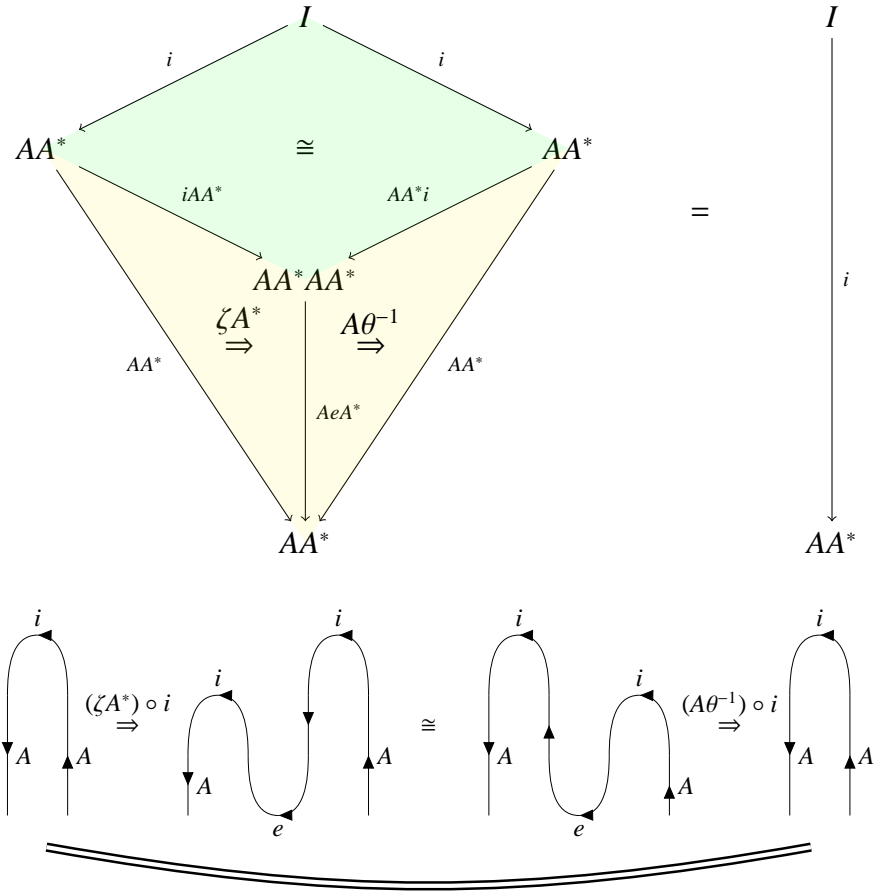
\begin{center}
    \begin{tikzpicture}[scale=1.75]
      \node (A) at (2,0) {$AA^*$};
      \node (B) at (2,2) {$AA^*AA^*$}
        edge [->] node [right, l] {$AeA^*$} (A);
      \node (C) at (0,3) {$AA^*$}
        edge [->] node [below left, l] {$AA^*$} (A)
        edge [->] node [above right, l] {$iAA^*$} (B);
      \node at (1.5,1.5) [above] {$\zeta A^*$};
      \node at (1.5,1.5) {$\Rightarrow$};
      \begin{scope}
        \fill [white,fill=yellow,fill opacity=0.1] (A.center) to (B.center) to (C.center) to (A.center);
      \end{scope}
      \node (D) at (4,3) {$AA^*$}
        edge [->] node [below right, l] {$AA^*$} (A)
        edge [->] node [above left, l] {$AA^*i$} (B);
      \begin{scope}
        \fill [white,fill=yellow,fill opacity=0.1] (A.center) to (B.center) to (D.center) to (A.center);
      \end{scope}
      \node at (2.5,1.5) [above] {$A\theta^{-1}$};
      \node at (2.5,1.5) {$\Rightarrow$};
      \node (E) at (2,4) {$I$}
        edge [->] node [above left, l] {$i$} (C)
        edge [->] node [above right, l] {$i$} (D);
      \node at (2,3) {$\cong$};
      \begin{scope}
        \fill [white,fill=green,fill opacity=0.1] (B.center) to (D.center) to (E.center) to (C.center) to (B.center);
      \end{scope}
      \node at (5,2.5) {$=$};
      \node (F) at (6,0) {$AA^*$};
      \node (G) at (6,4) {$I$}
        edge [->] node [right, l] {$i$} (F);
    \end{tikzpicture}
  \end{center}
  \begin{center}
    \scalebox{0.8}{
    \begin{tikzpicture}
      \id[0,2,0,0,A];
      \cap[0,2,];
      \id[1,0,1,2,A];

      \id[3,1,3,0,A];
      \cap[3,1,];
      \cup[4,0,];
      \id[5,2,5,1,];
      \cap[5,2,];
      \id[6,0,6,2,A];

      \id[8,2,8,0,A];
      \cap[8,2,];
      \id[9,1,9,2,];
      \cup[9,0,];
      \cap[10,1,];
      \id[11,0,11,1,A];

      \id[13,2,13,0,A];
      \cap[13,2,];
      \id[14,0,14,2,A];
      
      \node at (2,1.5) [above] {$(\zeta A^*)\circ i$};
      \node at (2,1.5) {$\Rightarrow$};
      \node at (7,1.5) {$\cong$};
      \node at (12,1.5) [above] {$(A\theta^{-1})\circ i$};
      \node at (12,1.5) {$\Rightarrow$};
      \node (B) at (13.5, -.5) {};
      \node (A) at (.5,-.5) {}
        edge [very thick,double distance=2pt, bend right=10] (B);
    \end{tikzpicture}
    }
  \end{center}
We have drawn the diagrams in a strictly monoidal compact closed bicategory for clarity; when the associator is not the identity, we truncate some corners:

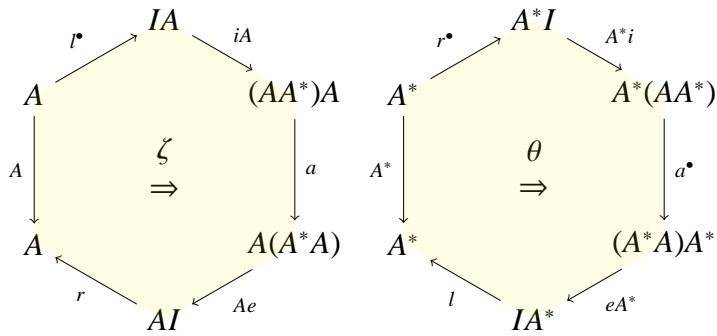
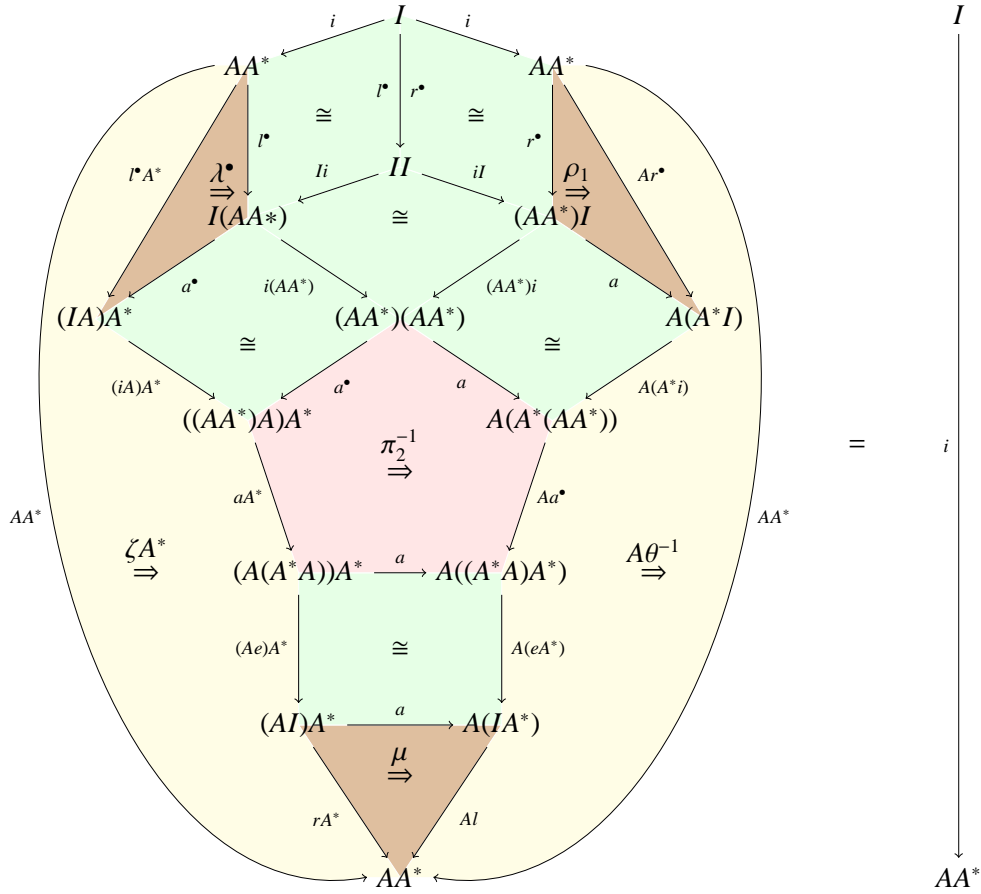
\begin{center}
  \begin{tikzpicture}
    \node (A) at (30:2cm) {$(AA^*)A$};
    \node (B) at (90:2cm) {$IA$}
      edge [->] node [above right, l] {$iA$} (A);
    \node (C) at (150:2cm) {$A$}
      edge [->] node [above left, l] {$l^\bullet$} (B);
    \node (D) at (210:2cm) {$A$}
      edge [<-] node [left, l] {$A$} (C);
    \node (E) at (270:2cm) {$AI$}
      edge [->] node [below left, l] {$r$} (D);
    \node (F) at (330:2cm) {$A(A^*A)$}
      edge [->] node [below right, l] {$Ae$} (E)
      edge [<-] node [right, l] {$a$} (A);
    \node at (90:.25cm) {$\zeta$};
    \node at (-90:.25cm) {$\Rightarrow$};
    \begin{scope}
      \fill [white,fill=yellow,fill opacity=0.1] (A.center) to (B.center) to (C.center) to (D.center) to (E.center) to (F.center) to (A.center);
    \end{scope}
  \end{tikzpicture}
  \begin{tikzpicture}
    \node (A) at (30:2cm) {$A^*(AA^*)$};
    \node (B) at (90:2cm) {$A^*I$}
      edge [->] node [above right, l] {$A^*i$} (A);
    \node (C) at (150:2cm) {$A^*$}
      edge [->] node [above left, l] {$r^\bullet$} (B);
    \node (D) at (210:2cm) {$A^*$}
      edge [<-] node [left, l] {$A^*$} (C);
    \node (E) at (270:2cm) {$IA^*$}
      edge [->] node [below left, l] {$l$} (D);
    \node (F) at (330:2cm) {$(A^*A)A^*$}
      edge [->] node [below right, l] {$eA^*$} (E)
      edge [<-] node [right, l] {$a^\bullet$} (A);
    \node at (90:.25cm) {$\theta$};
    \node at (-90:.25cm) {$\Rightarrow$};
    \begin{scope}
      \fill [white,fill=yellow,fill opacity=0.1] (A.center) to (B.center) to (C.center) to (D.center) to (E.center) to (F.center) to (A.center);
    \end{scope}
  \end{tikzpicture}
\end{center}

\begin{center}
  \scalebox{0.9}{
  \begin{tikzpicture}[scale=.75]
    \filldraw [white,fill=brown,fill opacity=0.5] (6,0)--(4,3)--(8,3)--cycle;
    \filldraw [white,fill=green,fill opacity=0.1] (4,6)--(8,6)--(8,3)--(4,3)--cycle;
    \filldraw [white,fill=red,fill opacity=0.1] (4,6)--(8,6)--(9,9)--(6,11)--(3,9)--cycle;
    \filldraw [white,fill=green,fill opacity=0.1] (3,9)--(0,11)--(3,13)--(6,11)--cycle;
    \filldraw [white,fill=green,fill opacity=0.1] (6,11)--(9,9)--(12,11)--(9,13)--cycle;
    \filldraw [white,fill=green,fill opacity=0.1] (6,11)--(3,13)--(6,14)--(9,13)--cycle;
    \filldraw [white,fill=green,fill opacity=0.1] (3,13)--(3,16)--(6,17)--(6,14)--cycle;
    \filldraw [white,fill=green,fill opacity=0.1] (6,17)--(6,14)--(9,13)--(9,16)--cycle;
    \filldraw [white,fill=brown,fill opacity=0.5] (0,11)--(3,13)--(3,16)--cycle;
    \filldraw [white,fill=brown,fill opacity=0.5] (12,11)--(9,13)--(9,16)--cycle;
    
    \node (A) at (6,0) {$AA^*$};
    \node (B) at (4,3) {$(AI)A^*$}
      edge [->] node [below left, l] {$rA^*$} (A);
    \node (C) at (8,3) {$A(IA^*)$}
      edge [<-] node [above, l] {$a$} (B)
      edge [->] node [below right, l] {$Al$} (A);
    \node at (6,2) [above] {$\mu$};
    \node at (6,2) {$\Rightarrow$};
    \node (D) at (4,6) {$(A(A^*A))A^*$}
      edge [->] node [left,l] {$(Ae)A^*$} (B);
    \node (E) at (8,6) {$A((A^*A)A^*)$}
      edge [<-] node [above, l] {$a$} (D)
      edge [->] node [right, l] {$A(eA^*)$} (C);
    \node at (6,4.5) {$\cong$};
    \node (F) at (3,9) {$((AA^*)A)A^*$}
      edge [->] node [left,l] {$aA^*$} (D);
    \node (G) at (9,9) {$A(A^*(AA^*))$}
      edge [->] node [right, l] {$Aa^\bullet$} (E);
    \node (H) at (0,11) {$(IA)A^*$}
      edge [->] node [below left, l] {$(iA)A^*$} (F);
    \node (I) at (6,11) {$(AA^*)(AA^*)$}
      edge [->] node [below right, l] {$a^\bullet$} (F)
      edge [->] node [below left, l] {$a$} (G);
    \node at (6,8) [above] {$\pi_2^{-1}$};
    \node at (6,8) {$\Rightarrow$};
    \node (J) at (12,11) {$A(A^*I)$}
      edge [->] node [below right, l] {$A(A^*i)$} (G);
    \node (K) at (3,13) {$I(AA*)$}
      edge [->] node [below right, l] {$a^\bullet$} (H)
      edge [->] node [below left, l] {$i(AA^*)$} (I);
    \node at (3,10.5) {$\cong$};
    \node (L) at (9,13) {$(AA^*)I$}
      edge [->] node [below right, l] {$(AA^*)i$} (I)
      edge [->] node [below left, l] {$a$} (J);
    \node at (9,10.5) {$\cong$};
    \node (M) at (3,16) {$AA^*$}
      edge [->] node [above left, l] {$l^\bullet A^*$} (H)
      edge [->] node [right, l] {$l^\bullet$} (K);
    \begin{scope}[on background layer]
      \fill[white,fill=yellow,fill opacity=0.1] (A.center) to (B.center) to (D.center) to (F.center) to (H.center) to (M.center) to (M.west) to [bend right=90] (A.west) to (A.center);
    \end{scope}
    \draw [->, black] (M.west) to [bend right=90] node [left, l] {$AA^*$} (A.west);
    \node at (2.5,13.5) [above] {$\lambda^\bullet$};
    \node at (2.5,13.5) {$\Rightarrow$};
    \node (N) at (9,16) {$AA^*$}
      edge [->] node [left, l] {$r^\bullet$} (L)
      edge [->] node [above right, l] {$Ar^\bullet$} (J);
    \begin{scope}[on background layer]
      \fill[white,fill=yellow,fill opacity=0.1] (A.center) to (C.center) to (E.center) to (G.center) to (J.center) to (N.center) to (N.east) to [bend left=90] (A.east) to (A.center);
    \end{scope}
    \draw [->, black] (N.east) to [bend left=90] node [right, l] {$AA^*$} (A.east);
    \node at (9.5,13.5) [above] {$\rho_1$};
    \node at (9.5,13.5) {$\Rightarrow$};
    \node (O) at (6,14) {$II$}
      edge [->] node [above left, l] {$Ii$} (K)
      edge [->] node [above right, l] {$iI$} (L);
    \node at (6,13) {$\cong$};
    \node (P) at (6,17) {$I$}
      edge [->] node [above, l] {$i$} (M)
      edge [->] node [above, l] {$i$} (N)
      edge [->] node [left, l] {$l^\bullet$} node [right, l] {$r^\bullet$} (O);
    \node at (4.5,15) {$\cong$};
    \node at (7.5,15) {$\cong$};
    \node at (1,6) [above] {$\zeta A^*$};
    \node at (1,6) {$\Rightarrow$};
    \node at (11,6) [above] {$A\theta^{-1}$};
    \node at (11,6) {$\Rightarrow$};
    \node at (15, 8.5) {$=$};

    \node (A') at (17,0) {$AA^*$};
    \node (P') at (17,17) {$I$}
      edge [->] node [left, l] {$i$} (A');
  \end{tikzpicture}
  }
\end{center}

\section{Bicategories of spans}
\label{Hoffnung}
At the start of this paper, we stated that spans of sets form a compact closed bicategory.  Street \cite{Street82} suggested weakening the notion of a map of spans to hold only up to 2-isomorphism, allowing to define spans in bicategories rather than mere categories; Hoffnung \cite
{Hoffnung} worked out the details.

A {\bf span} from $A$ to $B$ in a bicategory $T$ is a pair of morphisms with the same source: ${A \stackrel{f}{\leftarrow} C \stackrel{g}{\to} B.}$

A {\bf map of spans} $h$ between two spans ${A \stackrel{f}{\leftarrow} C \stackrel{g}{\to} B}$ and ${A \stackrel{f'}{\leftarrow} C' \stackrel{g'}{\to} B}$ is a triple $({h\maps C\to C',}$ ${\alpha\maps f \Rightarrow f'h,}$ ${\beta\maps g \Rightarrow g'h})$ such that $\alpha$ and $\beta$ are invertible.
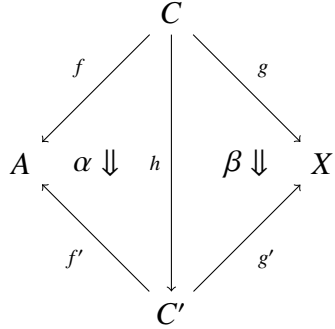
\begin{center}
  \begin{tikzpicture}[scale=2]
    \node (C) at (1,2) {$C$};
    \node (A) at (0,1) {$A$}
      edge [<-] node [l,above left] {$f$}(C);
    \node (B) at (2,1) {$X$}
      edge [<-] node [l,above right] {$g$}(C);
    \node (C') at (1,0) {$C'$}
      edge [->] node [l,below left] {$f'$}(A)
      edge [->] node [l,below right] {$g'$}(B)
      edge [<-] node [l,left] {$h$} (C);
    \node at (0.5,1) {$\alpha \Downarrow$};
    \node at (1.5,1) {$\beta \Downarrow$};
  \end{tikzpicture}
\end{center}
A {\bf map of maps of spans} is a 2-morphism $\gamma\maps h \Rightarrow h'$ such that $\alpha' = (f'\gamma)\cdot \alpha$ and $\beta' = (g'\gamma)\cdot \beta$.  Maps and maps of maps compose in the obvious ways.

Hoffnung showed that any 2-category $T$ with finite products and strict iso-comma objects (hereafter called ``weak pullbacks'') gives rise to a monoidal tricategory we will call $\Span_3(T)$ whose
\begin{itemize}
  \item objects are objects of $T$,
  \item morphisms are spans in $T$,
  \item 2-morphisms are maps of spans, and
  \item 3-morphisms are maps of maps of spans;
\end{itemize}
The tensor product of two spans $A \stackrel{f}{\leftarrow} C \stackrel{g}{\rightarrow} B$ 
and $A' \stackrel{f'}{\leftarrow} C' \stackrel{g'}{\rightarrow} B'$ is the span
\[A \times A' \stackrel{f \times f'}{\leftarrow} C \times C' \stackrel{g \times g'}{\rightarrow} B \times B'.\]

In this section, we will use $A$ to mean the object $A$, the identity 1-morphism on $A$, or the identity 2-morphism on the identity 1-morphism on $A,$ depending on the context.  Similarly, we will use $f$ to mean either the 1-morphism $f$ or the identity 2-morphism on $f$, depending on the context.  We will use juxtaposition to mean horizontal composition, {\em i.e.} composition of 1-morphisms: given $f\maps A \to B$ and $g\maps B \to C,$ we get $gf\maps A \to C.$  We also use juxtaposition to denote ``whiskering'': given a 2-morphism $K,$ we denote by $fK$ the the horizontal composite of $K$ and the identity 2-morphism on $f.$  We denote vertical composition of 2-morphisms $K\maps f\Rightarrow g$ and $L\maps g\Rightarrow h$ by $L\cdot K\maps f\Rightarrow h.$  We use $\pi_n$ to mean a projection out of a weak pullback, {\em not} a variant of the pentagonator 2-morphism.

We define composition of spans using the weak pullback in $T$.  The weak pullback of a cospan ${A \stackrel{f}{\to} C \stackrel{g}{\leftarrow} B}$ consists of an object $A_{f,g}B,$ 1-morphisms ${\pi_1\maps A_{f,g}B \to A}$ and ${\pi_2\maps A_{f,g}B \to B}$, and an invertible 2-morphism ${K\maps f\pi_1 \Rightarrow g\pi_2.}$  

The weak pullback satisfies two universal properties.  First, given any competitor $(X, {\pi'_1\maps X \to A,}$ ${\pi'_2\maps X \to B,}$ ${K'\maps f\pi'_1 \Rightarrow g\pi'_2})$ where $K'$ is invertible, there exists a unique 1-morphism ${\langle\pi'_1, \pi'_2\rangle\maps X \Rightarrow A_{f,g}B}$ such that ${\pi_1\langle\pi'_1, \pi'_2\rangle = \pi'_1,}$ ${\pi_2\langle\pi'_1, \pi'_2\rangle = \pi'_2,}$ and ${K\langle\pi'_1, \pi'_2\rangle = K'.}$
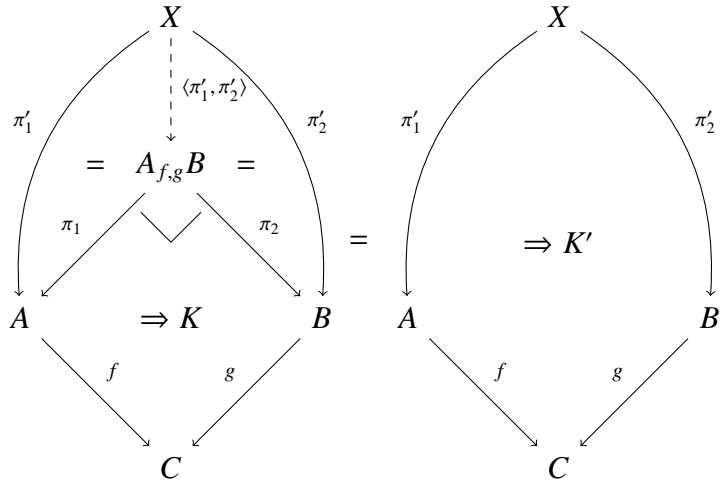
\begin{center}
  \begin{tikzpicture}[scale=2]
    \node (C) at (2,0) {$C$};
    \node (A) at (1,1) {$A$}
      edge [->] node [above right, l] {$f$} (C);
    \node (B) at (3,1) {$B$}
      edge [->] node [above left, l] {$g$} (C);
    \node (X) at (2,2) {$A_{f,g}B$}
      edge [->] node [above left, l] {$\pi_1$} (A)
      edge [->] node [above right, l] {$\pi_2$} (B);
    \node at (2,1) {$\Rightarrow K$};
    \draw (1.8,1.7)--(2,1.5)--(2.2,1.7);
    
    \node at (2, 3) {$X$}
      edge [->, dashed] node [right,l] {$\langle\pi'_1, \pi'_2\rangle$} (X)
      edge [->, bend right] node [above left, l] {$\pi'_1$} (A)
      edge [->, bend left] node [above right, l] {$\pi'_2$} (B);
    \node at (1.5, 2) {$=$};
    \node at (2.5, 2) {$=$};
    
    \node at (3.25, 1.5) {$=$};
  \end{tikzpicture}
  \begin{tikzpicture}[scale=2]
    \node (C) at (2,0) {$C$};
    \node (A) at (1,1) {$A$}
      edge [->] node [above right, l] {$f$} (C);
    \node (B) at (3,1) {$B$}
      edge [->] node [above left, l] {$g$} (C);    
    \node at (2, 3) {$X$}
      edge [->, bend right] node [above left, l] {$\pi'_1$} (A)
      edge [->, bend left] node [above right, l] {$\pi'_2$} (B);
    \node at (2, 1.5) {$\Rightarrow K'$};
  \end{tikzpicture}
\end{center}
Second, given any object $Y$, 1-morphisms ${j,k\maps Y \to A_{f,g}B,}$ and invertible 2-morphisms ${\omega\maps \pi_1j \Rightarrow \pi_1k}$ and ${\rho\maps \pi_2j \Rightarrow \pi_2k}$ such that
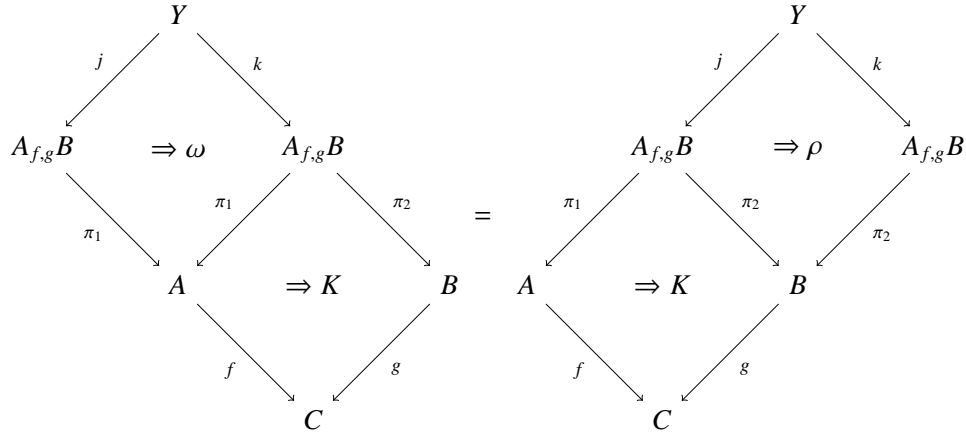
\begin{center}
  \scalebox{0.9}{
  \begin{tikzpicture}[scale=2]
    \node (C) at (1,0) {$C$};
    \node (A) at (0,1) {$A$}
      edge [->] node [below left,l] {$f$} (C);
    \node (B) at (2,1) {$B$}
      edge [->] node [below right,l] {$g$} (C);
    \node (P2) at (1,2) {$A_{f,g}B$}
      edge [->] node [above left,l] {$\pi_1$} (A)
      edge [->] node [above right,l] {$\pi_2$} (B);
    \node (P1) at (-1,2) {$A_{f,g}B$}
      edge [->] node [below left,l] {$\pi_1$} (A);
    \node (Y) at (0,3) {$Y$}
      edge [->] node [above left,l] {$j$} (P1)
      edge [->] node [above right,l] {$k$} (P2);
    \node at (1,1) {$\Rightarrow K$};
    \node at (0,2) {$\Rightarrow \omega$};
    \node at (2.25,1.5) {$=$};
  \end{tikzpicture}
  \begin{tikzpicture}[scale=2]
    \node (C) at (1,0) {$C$};
    \node (A) at (0,1) {$A$}
      edge [->] node [below left,l] {$f$} (C);
    \node (B) at (2,1) {$B$}
      edge [->] node [below right,l] {$g$} (C);
    \node (P1) at (1,2) {$A_{f,g}B$}
      edge [->] node [above left,l] {$\pi_1$} (A)
      edge [->] node [above right,l] {$\pi_2$} (B);
    \node (P2) at (3,2) {$A_{f,g}B$}
      edge [->] node [below right,l] {$\pi_2$} (B);
    \node (Y) at (2,3) {$Y$}
      edge [->] node [above left,l] {$j$} (P1)
      edge [->] node [above right,l] {$k$} (P2);
    \node at (1,1) {$\Rightarrow K$};
    \node at (2,2) {$\Rightarrow \rho$};
  \end{tikzpicture}
  }
\end{center}
there is a unique 2-morphism ${\gamma\maps j \Rightarrow k}$ such that ${\omega = \pi_1\gamma}$ and ${\rho = \pi_2\gamma.}$  

Here we show that the bicategory $\Span_2(T)$ whose
\begin{itemize}
  \item objects are objects of $T$,
  \item morphisms are spans in $T$, and
  \item 2-morphisms are 3-isomorphism classes of maps of spans.
\end{itemize}
forms a compact closed bicategory whenever $T$ is a 2-category with finite products and weak pullbacks.

Weak pullbacks are unique up to isomorphism \cite{nlab2pullbacks}.  The construction of $\Span_2(T)$ requires choosing specific weak pullbacks for each cospan \cite[3.2.1]{Hoffnung}; in our proof below, we choose especially nice pullbacks for the kinds of cospan that appear in the definition of a compact closed bicategory.  

This raises the question of whether some choices of weak pullback are fundamentally different than others.  As far as we know, nobody has proved that different choices give equivalent bicategories $\Span_2(T)$.  We conjecture that this is true, but for now we simply go ahead and take a particularly convenient choice.

Every weak pullback of a cospan comes equipped with two projections out of it. Now suppose that we compose four identity spans on $A,$ starting at the left; the resulting weak pullback is $((A_{A,A}A)_{\pi_2, A}A)_{\pi_2', A}A:$ 
\begin{center}
  \begin{equation}
    \label{fourfold}
    \begin{tikzpicture}
      \node (A1) at (0,0) {$A$};
      \node (A2) at (2,0)  {$A$};
      \node (A3) at (4,0)  {$A$};
      \node (A4) at (6,0)  {$A$};
      \node (A5) at (8,0)  {$A$};
      \node (AA1) at (1,1)  {$A$} edge [->] node [left,l] {$A$} (A1) edge [->] node [right,l] {$A$} (A2);
      \node (AA2) at (3,1)  {$A$} edge [->] node [left,l] {$A$} (A2) edge [->] node [right,l] {$A$} (A3);
      \node (AA3) at (5,1)  {$A$} edge [->] node [left,l] {$A$} (A3) edge [->] node [right,l] {$A$} (A4);
      \node (AA4) at (7,1)  {$A$} edge [->] node [left,l] {$A$} (A4) edge [->] node [right,l] {$A$} (A5);
      \node (AAA1) at (2,2) {$A_{A,A}A$} 
        edge [->] node [left,l] {$\pi_1$} (AA1)
        edge [->] node [right,l] {$\pi_2$} (AA2);
      \node (AAAA1) at (3,3) {$(A_{A,A}A)_{\pi_2,A}A$}
        edge [->] node [left,l] {$\pi'_1$} (AAA1)
        edge [->] node [right,l] {$\pi'_2$} (AA3);
      \node at (4,4) {$((A_{A,A}A)_{\pi_2,A}A)_{\pi'_2,A}A$}
        edge [->] node [left,l] {$\pi''_1$} (AAAA1)
        edge [->] node [right,l] {$\pi''_2$} (AA4);
      \node at (2,1) {$\stackrel{K}{\Rightarrow}$};
      \node at (3.5,1.5) {$\stackrel{K'}{\Rightarrow}$};
      \node at (5,2) {$\stackrel{K''}{\Rightarrow}$};
    \end{tikzpicture}
    .
  \end{equation}
\end{center}
This notation clearly becomes very cumbersome very quickly---particularly when dealing with the composite of many spans, as we will below.

We introduce a new notation $A^{\circ n}$ to mean the weak pullback in the composite of $n$ identity spans on $A$, beginning at the left; that is, $A^{\circ 1} = A, A^{\circ 2} = A_{A,A}A,$ and $A^{\circ n} = A^{\circ (n-1)}_{\pi_2, A}A,$ where $\pi_2\maps A^{\circ (n-1)} \to A$ is the second projection that $A^{\circ (n-1)}$ is equipped with.  

The construction $(-)^{\circ n}$ is an endofunctor on $T$; it takes an object $A$ to the object $A^{\circ n},$ a morphism $f\maps A\to B$ to the morphism $f^{\circ n}\maps A^{\circ n}\to B^{\circ n},$ and a 2-morphism $\alpha\maps f\Rightarrow g$ to the 2-morphism $\alpha^{\circ n}\maps f^{\circ n} \Rightarrow g^{\circ n}.$  For example, in the case $T=\Cat,$ the category $A^{\circ n}$ consists of length-$n$ chains of objects of $A$ equipped with isomorphisms between them:
\[ a_1 \stackrel{K_1}{\to} a_2 \stackrel{K_2}{\to} \cdots \stackrel{K_{n-1}}{\to} a_n. \]
Given a functor $f\maps A \to B$, the functor $f^{\circ n}\maps A^{\circ n}\to B^{\circ n}$ applies $f$ pointwise to each object and isomorphism in the chain:
\[ f(a_1) \stackrel{f(K_1)}{\to} f(a_2) \stackrel{f(K_2)}{\to} \cdots \stackrel{f(K_{n-1})}{\to} f(a_n). \]
Given a natural transformation $\alpha\maps f\Rightarrow g,$ the natural transformation $\alpha^{\circ n}\maps f^{\circ n}\Rightarrow g^{\circ n}$ assigns to each chain $a_1 \stackrel{K_1}{\to} a_2 \stackrel{K_2}{\to} \cdots \stackrel{K_{n-1}}{\to} a_n$ the list $(\alpha_{a_1},\ldots,\alpha_{a_n})$:
\[
\begin{array}{ccccccc}
  f(a_1) & \stackrel{f(K_1)}{\to} & f(a_2) & \stackrel{f(K_2)}{\to} & \cdots & \stackrel{f(K_{n-1})}{\to} & f(a_n) \\
  \alpha_{a_1}\downarrow & & \alpha_{a_2}\downarrow & & \cdots \downarrow & & \alpha_{a_n}\downarrow \\
  g(a_1) & \stackrel{g(K_1)}{\to} & g(a_2) & \stackrel{g(K_2)}{\to} & \cdots & \stackrel{g(K_{n-1})}{\to} & g(a_n)
\end{array}.
\]

Note that in diagram \ref{fourfold}, if we want to project from the apex onto the leftmost $A,$ we have to write $\pi_1\pi'_1\pi''_1;$ we are effectively forced to index the weak pullback using unary.  Going forward, we will write $\pi_1$ through $\pi_n$ for the $n$ projections out of $A^{\circ n}$ that result in an object of $A$.  In the case of $T=\Cat,$ for example,
\[ \pi_i(a_1 \stackrel{K_1}{\to} a_2 \stackrel{K_2}{\to} \cdots \stackrel{K_{n-1}}{\to} a_n) = a_i. \]
There is a dinatural transformation that assigns to each object $A$ of $T$ the morphism ${\pi_1\maps A^{\circ n} \to A,}$ and similarly for the other projections $\pi_i$; therefore we will use $\pi_i$ in a ``polymorphic'' way: we write both $\pi_1\maps A^{\circ n} \to A$ and $\pi_1\maps B^{\circ m} \to B,$ and will expect the reader to look at the source and target of such projections to determine exactly which morphism is being referred to.

\begin{lemma}
  \label{piling}
  Given isomorphisms $f\maps A \to C$ and $g\maps B \to C$, the weak pullback of the cospan $A^{\circ n} \stackrel{f\pi_n}{\to} C \stackrel{g\pi_1}{\leftarrow} B^{\circ m}$ is isomorphic to the weak pullback of the cospan $A^{\circ n} \stackrel{\pi_n}{\to} A \stackrel{\pi_1}{\leftarrow} A^{\circ m}.$
\end{lemma}
{\em Proof.}  
The weak pullbacks of the two cospans are
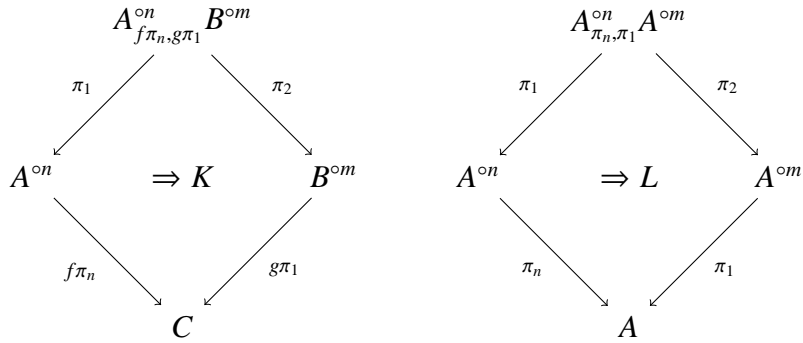
\begin{center}
  \begin{tikzpicture}[scale=2]
    \node(C) at (2,0) {$C$};
    \node(A) at (1,1) {$A^{\circ n}$}
      edge [->] node [below left,l] {$f\pi_n$} (C);
    \node(B) at (3,1) {$B^{\circ m}$}
      edge [->] node [below right,l] {$g\pi_1$} (C);
    \node(AfgB) at (2,2) {$A^{\circ n}_{f\pi_n,g\pi_1}B^{\circ m}$}
      edge [->] node [above left,l] {$\pi_1$} (A)
      edge [->] node [above right,l] {$\pi_2$} (B);
    \node at (2,1) {$\Rightarrow K$};
  \end{tikzpicture}
  $\quad \quad$
  \begin{tikzpicture}[scale=2]
    \node(C) at (2,0) {$A$};
    \node(A) at (1,1) {$A^{\circ n}$}
      edge [->] node [below left,l] {$\pi_n$} (C);
    \node(B) at (3,1) {$A^{\circ m}$}
      edge [->] node [below right,l] {$\pi_1$} (C);
    \node(AAAA) at (2,2) {$A^{\circ n}_{\pi_n,\pi_1}A^{\circ m}$}
      edge [->] node [above left,l] {$\pi_1$} (A)
      edge [->] node [above right,l] {$\pi_2$} (B);
    \node at (2,1) {$\Rightarrow L$};
  \end{tikzpicture}
  .
\end{center}
By the first universal property of weak pullbacks, there exist unique morphisms from $A^{\circ n}_{\pi_n,\pi_1}A^{\circ m}$ to $A^{\circ n}_{f\pi_n,g\pi_1}B^{\circ m}$ and back making the following diagrams commute.  The unique morphisms are evidently inverses.
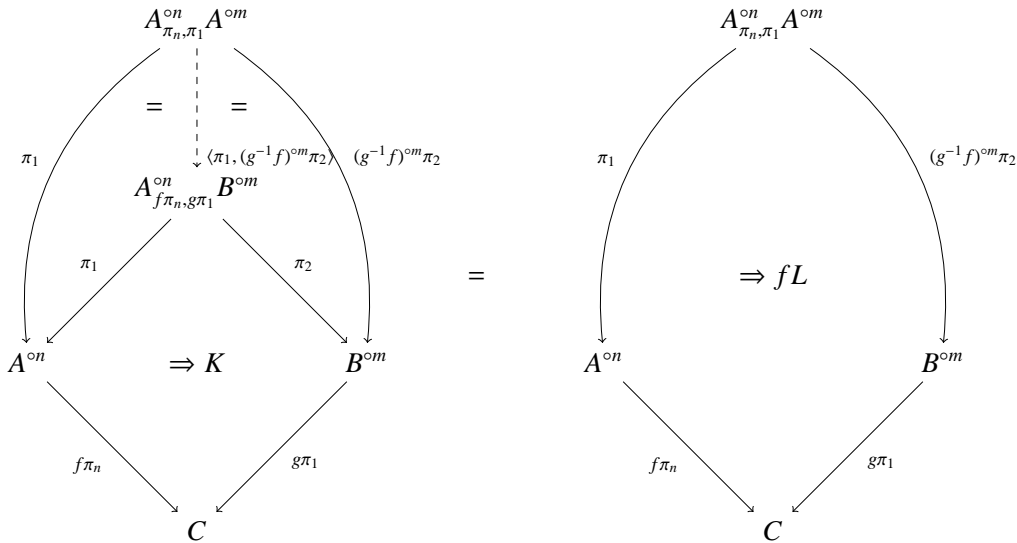
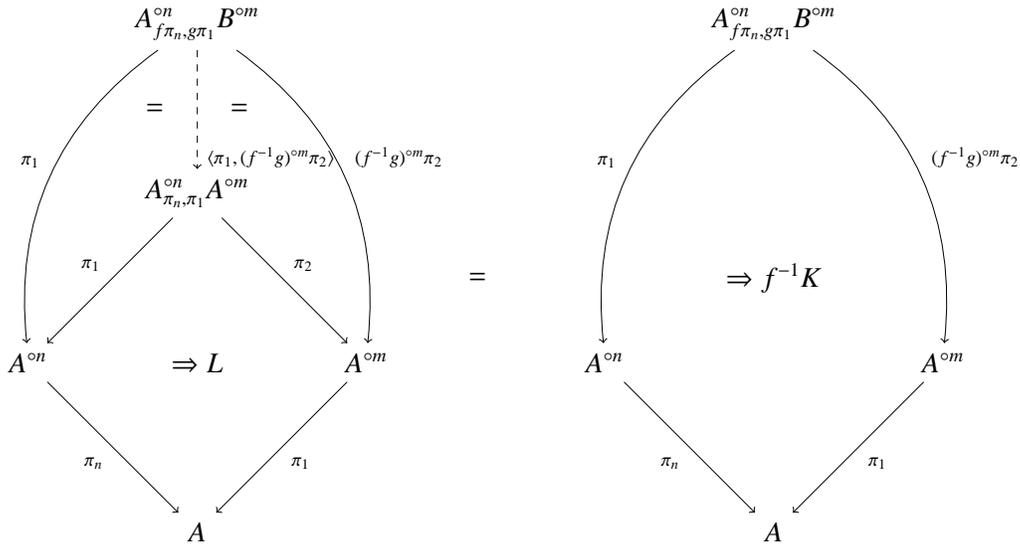
\begin{center}
  \scalebox{0.9}{
  \begin{tikzpicture}[scale=2.5]
    \node(C) at (2,0) {$C$};
    \node(A) at (1,1) {$A^{\circ n}$}
      edge [->] node [below left,l] {$f\pi_n$} (C);
    \node(B) at (3,1) {$B^{\circ m}$}
      edge [->] node [below right,l] {$g\pi_1$} (C);
    \node(AfgB) at (2,2) {$A^{\circ n}_{f\pi_n,g\pi_1}B^{\circ m}$}
      edge [->] node [above left,l] {$\pi_1$} (A)
      edge [->] node [above right,l] {$\pi_2$} (B);
    \node at (2,1) {$\Rightarrow K$};
    \node (AAAA) at (2,3) {$A^{\circ n}_{\pi_n,\pi_1}A^{\circ m}$}
      edge [->, dashed] node [right,l,pos=.9] {$\langle \pi_1, (g^{-1}f)^{\circ m}\pi_2\rangle$}(AfgB)
      edge [->, bend right] node [above left,l] {$\pi_1$} (A)
      edge [->, bend left] node [above right,l] {$(g^{-1}f)^{\circ m}\pi_2$} (B);
    \node at (1.75,2.5) {=};
    \node at (2.25,2.5) {=};
  \end{tikzpicture}
  \begin{tikzpicture}[scale=2.5]
    \node at (0.25,1.5) {$=$};
    \node(C) at (2,0) {$C$};
    \node(A) at (1,1) {$A^{\circ n}$}
      edge [->] node [below left,l] {$f\pi_n$} (C);
    \node(B) at (3,1) {$B^{\circ m}$}
      edge [->] node [below right,l] {$g\pi_1$} (C);
    \node at (2,1.5) {$\Rightarrow fL$};
    \node (AAAA) at (2,3) {$A^{\circ n}_{\pi_n,\pi_1}A^{\circ m}$}
      edge [->, bend right] node [above left,l] {$\pi_1$} (A)
      edge [->, bend left] node [above right,l] {$(g^{-1}f)^{\circ m}\pi_2$} (B);
  \end{tikzpicture}
  }
\end{center}

\begin{center}
  \scalebox{0.9}{
  \begin{tikzpicture}[scale=2.5]
    \node(C) at (2,0) {$A$};
    \node(A) at (1,1) {$A^{\circ n}$}
      edge [->] node [below left,l] {$\pi_n$} (C);
    \node(B) at (3,1) {$A^{\circ m}$}
      edge [->] node [below right,l] {$\pi_1$} (C);
    \node(AfgB) at (2,2) {$A^{\circ n}_{\pi_n,\pi_1}A^{\circ m}$}
      edge [->] node [above left,l] {$\pi_1$} (A)
      edge [->] node [above right,l] {$\pi_2$} (B);
    \node at (2,1) {$\Rightarrow L$};
    \node (AAAA) at (2,3) {$A^{\circ n}_{f\pi_n,g\pi_1}B^{\circ m}$}
      edge [->, dashed] node [right,l,pos=.9] {$\langle \pi_1, (f^{-1}g)^{\circ m}\pi_2\rangle$}(AfgB)
      edge [->, bend right] node [above left,l] {$\pi_1$} (A)
      edge [->, bend left] node [above right,l] {$(f^{-1}g)^{\circ m}\pi_2$} (B);
    \node at (1.75,2.5) {=};
    \node at (2.25,2.5) {=};
  \end{tikzpicture}
  \begin{tikzpicture}[scale=2.5]
    \node at (0.25,1.5) {$=$};
    \node(C) at (2,0) {$A$};
    \node(A) at (1,1) {$A^{\circ n}$}
      edge [->] node [below left,l] {$\pi_n$} (C);
    \node(B) at (3,1) {$A^{\circ m}$}
      edge [->] node [below right,l] {$\pi_1$} (C);
    \node at (2,1.5) {$\Rightarrow f^{-1}K$};
    \node (AAAA) at (2,3) {$A^{\circ n}_{f\pi_n,g\pi_1}B^{\circ m}$}
      edge [->, bend right] node [above left,l] {$\pi_1$} (A)
      edge [->, bend left] node [above right,l] {$(f^{-1}g)^{\circ m}\pi_2$} (B);
  \end{tikzpicture}
  }
\end{center}

Note that by the dinaturality of $\pi_1,$ $\pi_1 (g^{-1}f)^{\circ m} = g^{-1}f \pi_1,$ so the rightmost morphism on both sides of the top equation is $g g^{-1}f \pi_1 \pi_2 = f \pi_1\pi_2.$  Similarly, the rightmost morphism on both sides of the bottom equation is $\pi_1 (f^{-1}g)^{\circ m} \pi_2 = f^{-1}g \pi_1 \pi_2.$ \hfill $\square$

By the coherence theorem for bicategories, there is a unique isomorphism \[a^{\circ *}\maps A^{\circ n}_{\pi_n,\pi_1}A^{\circ m} \to A^{\circ (n+m)}\] built from associators for composition.  Since we must choose weak pullbacks for each cospan, given a cospan $A^{\circ n} \stackrel{f\pi_n}{\to} C \stackrel{g\pi_1}{\leftarrow} B^{\circ m}$ where $f$ and $g$ are invertible, we choose the weak pullback to be equal to $A^{\circ (n+m)}.$  When $A$ is terminal, for instance, $A$ may not equal $1$ but only be isomorphic.  In that case, the weak pullback of $A \stackrel{!}{\to} 1 \stackrel{!}{\leftarrow} 1$ is
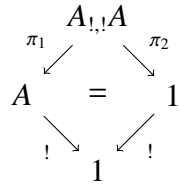
\begin{center}
  \begin{tikzpicture}
    \node (C) at (1,0) {$1$};
    \node (A) at (0,1) {$A$} edge [->] node [below left,l] {$!$} (C);
    \node (B) at (2,1) {$1$} edge [->] node [below right,l] {$!$} (C);
    \node at (1,2) {$A_{!,!}A$} edge [->] node[above left, l] {$\pi_1$} (A) edge [->] node[above right, l] {$\pi_2$} (B);
    \node at (1,1) {$=$};
  \end{tikzpicture}
\end{center}
whereas when $A$ is not terminal, the weak pullback is
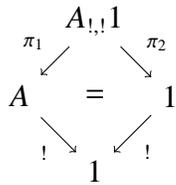
\begin{center}
  \begin{tikzpicture}
    \node (C) at (1,0) {$1$};
    \node (A) at (0,1) {$A$} edge [->] node [below left,l] {$!$} (C);
    \node (B) at (2,1) {$1$} edge [->] node [below right,l] {$!$} (C);
    \node at (1,2) {$A_{!,!}1$} edge [->] node[above left, l] {$\pi_1$} (A) edge [->] node[above right, l] {$\pi_2$} (B);
    \node at (1,1) {$=$};
  \end{tikzpicture}
  .
\end{center}

With that choice, we also have the following useful corollary.
\begin{cor}
  \label{reverse}
  Given an isomorphism $f\maps A \to B$ in $T,$ the composite of the identity span on $B$ and the span $B \stackrel{f}{\leftarrow} A \stackrel{A}{\to} A$ is equal to the composite of the identity span on $B$ and the span $B \stackrel{B}{\leftarrow} B \stackrel{f^{-1}}{\to} A$; both result in the span $B \stackrel{\pi_1}{\leftarrow} B^{\circ 2} \stackrel{f^{-1}\pi_2}{\to} A.$
\end{cor}

Because we mod out by isomorphisms of maps of spans, some spans that at first sight appear different are actually the same.
\begin{lemma}
  \label{swap}
  The braiding $b\maps A^{\circ 2} \to A^{\circ 2}$ in $T$ is 2-isomorphic to the identity.
\end{lemma}

{\em Proof.} The weak pullback of the identity cospan on $A$ is $A^{\circ 2}$ equipped with projections $\pi_1, \pi_2$ and a 2-morphism $L.$  We have $\pi_1 b = \pi_2, \pi_2 b = \pi_1,$ and $L b = L^{-1}.$  The following 2-morphisms are equal:
\begin{center}
  \scalebox{0.9}{
  \begin{tikzpicture}[scale=2]
    \node (C) at (1,0) {$A$};
    \node (A) at (0,1) {$A$}
      edge [->] node [below left,l] {$A$} (C);
    \node (B) at (2,1) {$A$}
      edge [->] node [below right,l] {$A$} (C);
    \node (P2) at (1,2) {$A^{\circ 2}$}
      edge [->] node [above left,l] {$\pi_1$} (A)
      edge [->] node [above right,l] {$\pi_2$} (B);
    \node (P1) at (-1,2) {$A^{\circ 2}$}
      edge [->] node [below left,l] {$\pi_1$} (A);
    \node (X) at (0,3) {$A^{\circ 2}$}
      edge [->] node [above left,l] {$b$} (P1)
      edge [->] node [above right,l] {$A^{\circ 2}$} (P2);
    \node at (1,1) {$\Rightarrow L$};
    \node at (0,2) {$\Rightarrow L^{-1}$};
    \node at (2.25,1.5) {$=$};
  \end{tikzpicture}
  \begin{tikzpicture}[scale=2]
    \node (C) at (1,0) {$A$};
    \node (A) at (0,1) {$A$}
      edge [->] node [below left,l] {$A$} (C);
    \node (B) at (2,1) {$A$}
      edge [->] node [below right,l] {$A$} (C);
    \node (P1) at (1,2) {$A^{\circ 2}$}
      edge [->] node [above left,l] {$\pi_1$} (A)
      edge [->] node [above right,l] {$\pi_2$} (B);
    \node (P2) at (3,2) {$A^{\circ 2}$}
      edge [->] node [below right,l] {$\pi_2$} (B);
    \node (X) at (2,3) {$A^{\circ 2}$}
      edge [->] node [above left,l] {$b$} (P1)
      edge [->] node [above right,l] {$A^{\circ 2}$} (P2);
    \node at (1,1) {$\Rightarrow L$};
    \node at (2,2) {$\Rightarrow L$};
  \end{tikzpicture}
  }
\end{center}
(note that on the right hand side, the lower use of $L$ is whiskered by $b,$ becoming $L^{-1}$), so by the second universal property of the weak pullback, there exists a unique 2-isomorphism $\gamma\maps b \Rightarrow A^{\circ 2}$ such that $L^{-1} = \pi_1 \gamma$ and $L = \pi_2 \gamma.$ \hfill $\square$

\begin{cor}
  The weak pullback of the identity cospan on $A$ is $A^{\circ 2}$ equipped with the projections $\pi_1\maps A^{\circ 2}\to A, \pi_2\maps A^{\circ 2}\to A,$ and a 2-morphism $L\maps \pi_1\Rightarrow \pi_2.$  The map of spans
  \begin{center}
    \begin{tikzpicture}[scale=2]
      \node (C) at (1,0) {$A^{\circ 2}$};
      \node (A) at (0,1) {$A$}
        edge [<-] node [below left,l] {$\pi_1$} (C);
      \node (B) at (2,1) {$A$}
        edge [<-] node [below right,l] {$\pi_2$} (C);
      \node (D) at (1,2) {$A^{\circ 2}$}
        edge [->] node [above left,l] {$\pi_1$} (A)
        edge [->] node [above right,l] {$\pi_2$} (B)
        edge [->] node [right,l] {$b$} (C);
      \node at (.5,1) {$\Downarrow L$};
      \node at (1.5,1) {$\Downarrow L^{-1}$};
    \end{tikzpicture}
  \end{center}
  is in the same equivalence class as the identity map of spans.
\end{cor}

\begin{cor}
  For any permutation $\sigma$ of $n$ elements, the morphism \[\langle \pi_{\sigma(1)}, \pi_{\sigma(2)}, \ldots, \pi_{\sigma(n)}\rangle\maps A^{\circ n} \to A^{\circ n}\] is 2-isomorphic to the identity.
\end{cor}

\begin{cor}
  \label{permute}
  The composite of $n$ identity spans on $A$ has as its apex the weak pullback consisting of the object $A^{\circ n}$ equipped with projections ${\pi_1, \ldots, \pi_n\maps A^{\circ n} \to A}$ and for each ${1 \le i < n}$ an invertible 2-morphism $L_i\maps \pi_i \Rightarrow \pi_{i+1}.$  Let $L'$ be the invertible 2-morphism from $\pi_1$ to $\pi_{\sigma(1)}$ and $L''$ be the invertible 2-morphism from $\pi_n$ to $\pi_{\sigma(n)}$ derived from composing the $L_i$.  The map of spans 
  \begin{center}
    \begin{tikzpicture}[scale=2]
      \node (C) at (1,0) {$A^{\circ n}$};
      \node (A) at (0,1) {$A$}
        edge [<-] node [below left,l] {$\pi_1$} (C);
      \node (B) at (2,1) {$A$}
        edge [<-] node [below right,l] {$\pi_n$} (C);
      \node (D) at (1,2) {$A^{\circ n}$}
        edge [->] node [above left,l] {$\pi_1$} (A)
        edge [->] node [above right,l] {$\pi_n$} (B)
        edge [->] node [right,l, pos=0.6] {$p$} (C);
      \node at (.5,1) {$\Downarrow L'$};
      \node at (1.5,1) {$\Downarrow L''$};
    \end{tikzpicture}
  \end{center}
  where $p = \langle \pi_{\sigma(1)}, \pi_{\sigma(2)}, \ldots, \pi_{\sigma(n)}\rangle$ is in the same equivalence class as the identity map of spans.
\end{cor}  


We are now ready to prove the main theorem.

\begin{thm}
  \label{Span2}
  If $T$ is a 2-category with finite products and weak pullbacks, then $\Span_2(T)$ is a compact closed bicategory.
\end{thm}
{\em Proof.} As noted, Hoffnung \cite{Hoffnung} showed that $\Span_3(T)$
is a monoidal tricategory.  We refer the reader to Hoffnung's paper for the complete definition of a monoidal tricategory, but suffice it to say that it replaces the commuting polyhedra in the above definition of a monoidal bicategory with polyhedra that commute up to a specified 3-morphism, and then adds coherence law polytopes to govern them.  When we mod out by 3-isomorphism classes of maps of spans, these 3-morphisms become trivial, so $\Span_2(T)$ is a monoidal bicategory.

The monoidal associator is the span $$(A\times B)\times C \stackrel{(A\times B)\times C}{\leftarrow} (A\times B)\times C \stackrel{a}{\to} A\times (B\times C).$$  The left and right monoidal unitors are the spans $$1\times A \stackrel{1\times A}{\leftarrow} 1\times A \stackrel{l}{\to} A$$ and $$A\times 1 \stackrel{A\times 1}{\leftarrow} A\times 1 \stackrel{r}{\to} A,$$ respectively.  The monoidal braiding is $$A\times B \stackrel{A\times B}{\leftarrow} A\times B \stackrel{b}{\to} B\times A.$$   The ``bulleted'' morphisms like $a^\bullet$ are the reverse spans.

To define the pentagonator, we start with a ``six-edged'' identity map of spans: each edge is a span whose left leg is the identity and whose right leg is an isomorphism in $T;$ the source and target composite spans are both the composite of three such edges, so by our choice of weak pullbacks, their apexes are equal.
\begin{center}
  \scalebox{0.9}{
  \begin{tikzpicture}[scale=1.5,font=\fontsize{10}{10}\selectfont]
    \filldraw[white,fill=red,fill opacity=0.1](0:2cm)--(60:2cm)--(120:2cm)--(180:2cm)--(240:2cm)--(300:2cm)--cycle;
    \node (A) at (  0:2cm) {$A\tensor (B\tensor (C\tensor D))$};
    \node (B) at ( 60:2cm) {$(A\tensor B)\tensor (C\tensor D)$}
      edge [->] node [l, above right] {$a$} (A);
    \node (C) at (120:2cm) {$((A\tensor B)\tensor C)\tensor D$}
      edge [->] node [l, above] {$a$} (B);
    \node (D) at (180:2cm) {$((A\tensor B)\tensor C)\tensor D$}
      edge [->] node [l, above left] {$((A\tensor B)\tensor C)\tensor D$} (C);
    \node (E) at (240:2cm) {$(A\tensor (B\tensor C))\tensor D$}
      edge [<-] node [l, below left] {$a\tensor D$} (D);
    \node (F) at (300:2cm) {$A\tensor ((B\tensor C)\tensor D)$}
      edge [<-] node [l, below] {$a$} (E)
      edge [->] node [l, below right] {$A\tensor a$} (A);
    \node at (0,0) {$=$};
  \end{tikzpicture}
  \begin{tikzpicture}[scale=2.75,font=\fontsize{10}{10}\selectfont]
    \node (C) at (1,0) {$(((A\times B)\times C)\times D)^{\circ 3}$};
    \node (A) at (0,1) {$((A\times B)\times C)\times D$}
      edge [<-] node [below left,l] {$\pi_1$} (C);
    \node (B) at (2,1) {$A\times (B\times (C\times D))$}
      edge [<-] node [below right,l] {$(A\times a) a (a\times D) \pi_3$} (C);
    \node (D) at (1,2) {$(((A\times B)\times C)\times D)^{\circ 3}$}
      edge [->] node [above left,l] {$\pi_1$} (A)
      edge [->] node [above right,l] {$a a \pi_3$} (B)
      edge [->] node [right,font=\fontsize{6}{6}\selectfont, pos=.6] {$(((A\times B)\times C)\times D)^{\circ 3}$} (C);
    \node at (.75,1.25) {$=$};
    \node at (1.25,1.25) {$=$};
  \end{tikzpicture}
  }
\end{center}
The right-hand 2-morphism in the map of spans is an identity because the pentagon equation holds in the underlying category of $T.$  We define the pentagonator to be the composite of this identity map of spans with the unitor for composition:  

\begin{center}
  \scalebox{0.9}{
  \begin{tikzpicture}[scale=2]
    \filldraw[white,fill=red,fill opacity=0.1](0,3)--(2,4)--(4,2)--(2,0)--(0,1)--cycle;
    \node (ABCD1) at (0,3) {$((A \tensor B)\tensor C)\tensor D$};
    \node (ABCD2) at (2,4) {$(A\tensor B)\tensor (C\tensor D)$}
      edge [<-] node [l, above left] {$a$} (ABCD1);
    \node (ABCD3) at (4,2) {$A\tensor (B\tensor (C\tensor D))$}
      edge [<-] node [l, above right] {$a$} (ABCD2);
    \node (ABCD4) at (2,0) {$A\tensor ((B\tensor C)\tensor D)$}
      edge [->] node [l, below right] {$A\tensor a$} (ABCD3);
    \node (ABCD5) at (0,1) {$(A\tensor (B\tensor C))\tensor D$}
      edge [->] node [l, below left] {$a$} (ABCD4)
      edge [<-] node [l, left] {$a\tensor D$} (ABCD1);
    \node (ABCD6) at (1.5,2.5) {$((A\tensor B)\tensor C)\tensor D$}
      edge [<-] node [l, auto] {$((A\tensor B)\tensor C)\tensor D$} (ABCD1)
      edge [->] node [l, auto] {$a$} (ABCD2);
    \node at (1.25,3) {\tikz\node [rotate=90] {$\Rightarrow$};};
    \node at (1.4,3) {$l^\circ$};
    \node at (2,2) {$=$};
  \end{tikzpicture}
  }
\end{center}

The coherence theorem for bicategories \cite{LeinsterBB} says that any diagram built out of $a^\circ, l^\circ,$ and $r^\circ$ commutes, so any coherence law involving only pentagonators and identity 2-morphisms---such as the associahedron---must hold in $\Span_2(T).$

To define the left 2-unitor for the monoidal product, we start with a ``four-edged'' identity map of spans.  Each edge is a span whose left leg is the identity and whose right leg is an isomorphism in $T;$ the source and target composite spans are both the composite of two such edges, so by our choice of weak pullbacks, their apexes are equal.

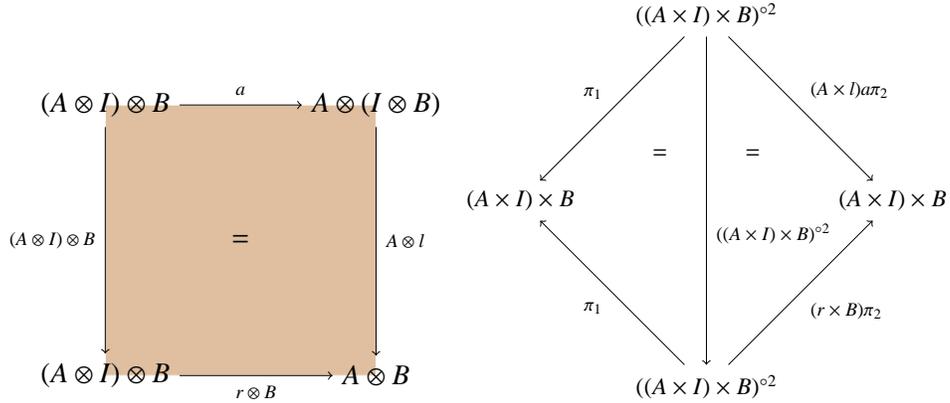
\begin{center}
  \scalebox{0.9}{
  \begin{tikzpicture}[scale=4]
    \filldraw[white,fill=brown,fill opacity=0.5](0,0)--(0,1)--(1,1)--(1,0)--cycle;
    \node (A) at (0,1) {$(A\tensor I) \tensor B$};
    \node (B) at (1,1) {$A\tensor (I \tensor B)$}
      edge [<-] node [above, l] {$a$} (A);
    \node (C) at (0,0) {$(A\tensor I) \tensor B$}
      edge [<-] node [left, l] {$(A\tensor I) \tensor B$} (A);
    \node (D) at (1,0) {$A\tensor B$}
      edge [<-] node [below, l] {$r \tensor B$} (C)
      edge [<-] node [right, l] {$A \tensor l$} (B);
    \node at (.5,.5) {$=$};
  \end{tikzpicture}
  \begin{tikzpicture}[scale=2.75,font=\fontsize{10}{10}\selectfont]
    \node (C) at (1,0) {$((A\times I) \times B)^{\circ 2}$};
    \node (A) at (0,1) {$(A\times I) \times B$}
      edge [<-] node [below left,l] {$\pi_1$} (C);
    \node (B) at (2,1) {$(A\times I) \times B$}
      edge [<-] node [below right,l] {$(r\times B) \pi_2$} (C);
    \node (D) at (1,2) {$((A\times I) \times B)^{\circ 2}$}
      edge [->] node [above left,l] {$\pi_1$} (A)
      edge [->] node [above right,l] {$(A\times l)a\pi_2$} (B)
      edge [->] node [right,l, pos=.6] {$((A\times I) \times B)^{\circ 2}$} (C);
    \node at (.75,1.25) {$=$};
    \node at (1.25,1.25) {$=$};
  \end{tikzpicture}
  }
\end{center}

The right-hand 2-morphism in the map of spans is an identity because the triangle equation holds in the underlying category of $T$.  We define the 2-unitor to be the composite of this identity map of spans with the inverse unitor for composition:

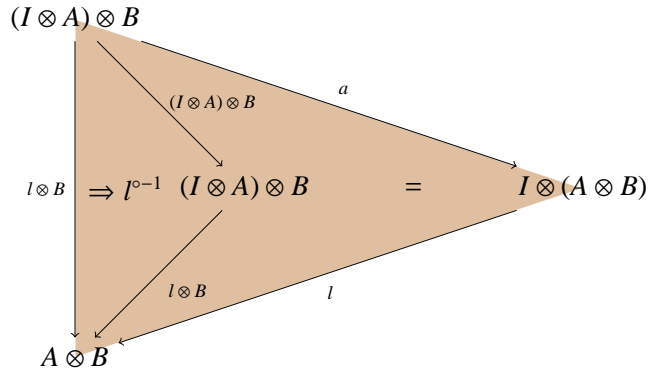
\begin{center}
  \scalebox{0.9}{
  \begin{tikzpicture}[scale=2.5]
    \filldraw[white,fill=brown,fill opacity=0.5](0,2)--(3,1)--(0,0)--cycle;
    \node (IAB1) at (0,2) {$(I\tensor A)\tensor B$};
    \node (IAB2) at (3,1) {$I\tensor (A\tensor B)$}
      edge [<-] node [l, above right] {$a$} (IAB1);
    \node (AB) at (0,0) {$A\tensor B$}
      edge [<-] node [l, below right] {$l$} (IAB2)
      edge [<-] node [l, left] {$l\tensor B$} (IAB1);
    \node at (1,1) {$(I\tensor A)\tensor B$}
      edge [<-] node [l, right] {$(I\tensor A)\tensor B$} (IAB1)
      edge [->] node [l, auto] {$l\tensor B$} (AB);
    \node at (.3,1) {$\Rightarrow l^{\circ-1}$};
    \node at (2,1) {$=$};
  \end{tikzpicture}
  }
\end{center}

The 2-unitors $\mu$ and $\rho$ are also equal to the inverse of the unitor for composition:
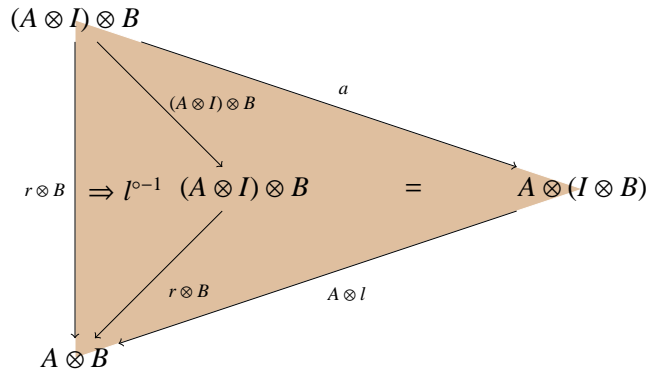
\begin{center}
  \scalebox{0.9}{
  \begin{tikzpicture}[scale=2.5]
    \filldraw[white,fill=brown,fill opacity=0.5](0,2)--(3,1)--(0,0)--cycle;
    \node (IAB1) at (0,2) {$(A\tensor I)\tensor B$};
    \node (IAB2) at (3,1) {$A\tensor (I\tensor B)$}
      edge [<-] node [l, above right] {$a$} (IAB1);
    \node (AB) at (0,0) {$A\tensor B$}
      edge [<-] node [l, below right] {$A\tensor l$} (IAB2)
      edge [<-] node [l, left] {$r\tensor B$} (IAB1);
    \node at (1,1) {$(A\tensor I)\tensor B$}
      edge [<-] node [l, right] {$(A\tensor I)\tensor B$} (IAB1)
      edge [->] node [l, auto] {$r\tensor B$} (AB);
    \node at (.3,1) {$\Rightarrow l^{\circ-1}$};
    \node at (2,1) {$=$};
  \end{tikzpicture}
  }
\end{center}

\begin{center}
  \scalebox{0.9}
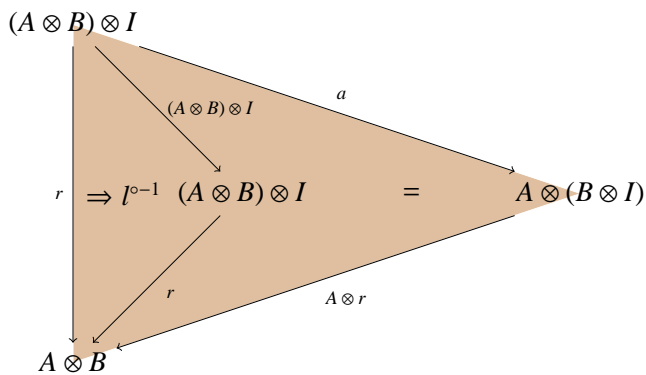
{
  \begin{tikzpicture}[scale=2.5]
    \filldraw[white,fill=brown,fill opacity=0.5](0,2)--(3,1)--(0,0)--cycle;
    \node (IAB1) at (0,2) {$(A\tensor B)\tensor I$};
    \node (IAB2) at (3,1) {$A\tensor (B\tensor I)$}
      edge [<-] node [l, above right] {$a$} (IAB1);
    \node (AB) at (0,0) {$A\tensor B$}
      edge [<-] node [l, below right] {$A\tensor r$} (IAB2)
      edge [<-] node [l, left] {$r$} (IAB1);
    \node at (1,1) {$(A\tensor B)\tensor I$}
      edge [<-] node [l, right] {$(A\tensor B)\tensor I$} (IAB1)
      edge [->] node [l, auto] {$r$} (AB);
    \node at (.3,1) {$\Rightarrow l^{\circ-1}$};
    \node at (2,1) {$=$};
  \end{tikzpicture}
  }
\end{center}

By the coherence theorem for bicategories, any diagram built out of $a^\circ, l^\circ,$ and $r^\circ$ commutes, so any coherence law involving only $\pi, \lambda, \mu, \rho$ and identity 2-morphisms---such as the unitor prisms---must hold in $\Span_2(T).$

The hexagon modification $R$ is a ``six-edged'' identity map of spans: each edge is, again, a span whose left leg is an identity and whose right leg is an isomorphism in $T$.  The source and target spans are the composite of three such edges, so because of our choice of weak pullbacks, the apexes are equal; the right-hand 2-morphism in the map of spans is an equality because the hexagon equations hold in $T.$

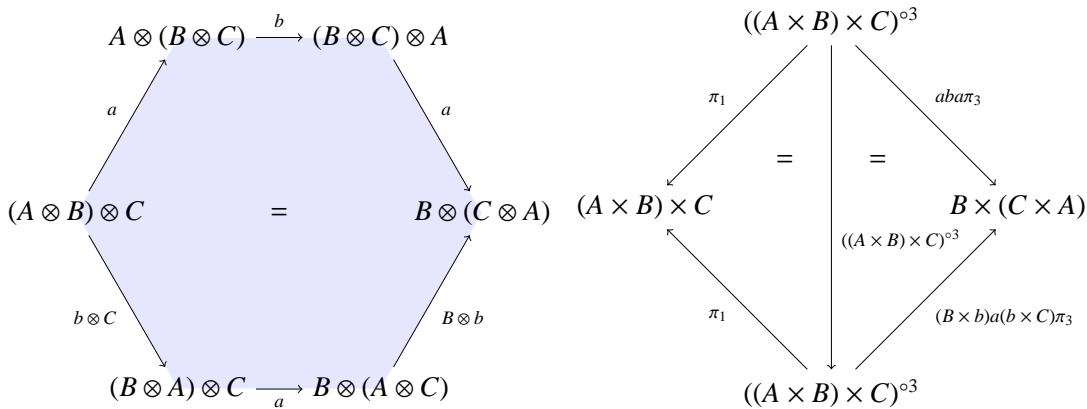
\begin{center}
  \scalebox{0.9}{
  \begin{tikzpicture}[scale=1.5]
    \filldraw[white,fill=blue,fill opacity=0.1](0:2cm)--(60:2cm)--(120:2cm)--(180:2cm)--(240:2cm)--(300:2cm)--cycle;
    \node (A) at (  0:2cm) {$B\tensor (C\tensor A)$};
    \node (B) at ( 60:2cm) {$(B\tensor C)\tensor A$}
      edge [->] node [l, above right] {$a$} (A);
    \node (C) at (120:2cm) {$A\tensor (B\tensor C)$}
      edge [->] node [l, above] {$b$} (B);
    \node (D) at (180:2cm) {$(A\tensor B)\tensor C$}
      edge [->] node [l, above left] {$a$} (C);
    \node (E) at (240:2cm) {$(B\tensor A)\tensor C$}
      edge [<-] node [l, below left] {$b\tensor C$} (D);
    \node (F) at (300:2cm) {$B\tensor (A\tensor C)$}
      edge [<-] node [l, below] {$a$} (E)
      edge [->] node [l, below right] {$B\tensor b$} (A);
    \node at (0,0) {$=$};
  \end{tikzpicture}
  \begin{tikzpicture}[scale=2.75]
    \node (C) at (1,0) {$((A\times B)\times C)^{\circ 3}$};
    \node (A) at (0,1) {$(A\times B)\times C$}
      edge [<-] node [below left,l] {$\pi_1$} (C);
    \node (B) at (2,1) {$B\times (C\times A)$}
      edge [<-] node [below right,l] {$(B\times b)a(b\times C)\pi_3$} (C);
    \node (D) at (1,2) {$((A\times B)\times C)^{\circ 3}$}
      edge [->] node [above left,l] {$\pi_1$} (A)
      edge [->] node [above right,l] {$aba\pi_3$} (B)
      edge [->] node [right,l, pos=.6] {$((A\times B)\times C)^{\circ 3}$} (C);
    \node at (.75,1.25) {$=$};
    \node at (1.25,1.25) {$=$};
  \end{tikzpicture}
  }
\end{center}

The hexagon modification $S$ is more complicated because it has three uses of $a^\bullet.$  To define $S$, we start with a ``ten-edged'' identity map of spans.  The edges are those of $S$ except that instead of using $a^\bullet$ it uses $a^{-1},$ and it also includes four extra identity edges.  Each edge is, again, a span whose left leg is an identity and whose right leg is an isomorphism in $T$.  The source and target spans are the composite of five such edges, so because of our choice of weak pullbacks, the apexes are equal; the right-hand 2-morphism in the map of spans is an equality because the hexagon equations hold in $T.$  By Corollary \ref{reverse}, the composite of an identity span with $a^\bullet$ is equal to the composite of an identity span with $a^{-1},$ so we define $S$ to be the composite of this identity span with four unitors for composition:

\begin{center}
  \begin{tikzpicture}[scale=2]
    \filldraw[white,fill=blue,fill opacity=0.1](0:2cm)--(60:2cm)--(120:2cm)--(180:2cm)--(240:2cm)--(300:2cm)--cycle;
    \node (A) at (  0:2cm) {$(CA)B$};
    \node (B) at ( 60:2cm) {$C(AB)$}
      edge [->] node [l, above right] {$a^\bullet$} (A);
    \node (C) at (120:2cm) {$(AB)C$}
      edge [->] node [l, above] {$b$} (B);
    \node (D) at (180:2cm) {$A(BC)$}
      edge [->] node [l, above left] {$a^\bullet$} (C);
    \node (E) at (240:2cm) {$A(CB)$}
      edge [<-] node [l, below left] {$Ab$} (D);
    \node (F) at (300:2cm) {$(AC)B$}
      edge [<-] node [l, below] {$a^\bullet$} (E)
      edge [->] node [l, below right] {$bB$} (A);
    \node at (150:1cm) {$A(BC)$}
      edge [<-] node [l, below] {$A(BC)$} (D)
      edge [->] node [l, right] {$a^{-1}$} (C);
    \node at (150:1.4cm) {$\Rightarrow l^{\circ -1}$};
    \node at (30:1cm) {$C(AB)$}
      edge [<-] node [l, left] {$C(AB)$} (B)
      edge [->] node [l,below] {$a^{-1}$} (A);
    \node at (30:1.4cm) {$\Leftarrow l^{\circ -1}$};
    \node (G) at (270:1cm) {$A(CB)$}
      edge [<-] node [l, left] {$A(CB)$} (E)
      edge [->] node [l,below] {$a^{-1}$} (F);
    \node at (330:1cm) {$A(CB)$}
      edge [<-] node [l, left] {$A(CB)$} (F)
      edge [->] node [l,above] {$bB$} (A);
    \node at (270:1.5) {$\Downarrow l^\circ$};
    \node at (330:1.5) {$\Rightarrow l^\circ$};
    \node at (0,0) {$=$};
  \end{tikzpicture}
\end{center}

By the coherence theorem for bicategories, any coherence law involving only $\pi, R, S$ and identity 2-morphisms---such as the shuffle and Breen polytopes---must hold in $\Span_2(T).$

To define the syllepsis, we begin with a ``four-edged'' identity span and compose it with two unitors for composition. By the coherence theorem for bicategories, any coherence law involving only $R, S, \nu$ and identity 2-morphisms---such as those governing the syllepsis---must hold in $\Span_2(T).$

\begin{center}
  \begin{tikzpicture}[scale=3, bend angle=270]
    \filldraw[white,fill=salmon,fill opacity=0.1] (1,1) ellipse (30pt and 20pt);
    \node (A) at (0,1) {$AB$};
    \node (B) at (1,1.3) {$AB$}
      edge [<-] node [l,above left] {$AB$} (A);
    \node (C) at (1,.6) {$AB$}
      edge [<-] node [l,below left] {$AB$} (A);
    \node (D) at (2,1) {$BA$}
      edge [<-] node [l,above right] {$b$} (B)
      edge [<-] node [l,below right] {$b^\bullet$} (C)
      edge [<-, bend right] node [l,below] {$b^\bullet$} (A)
      edge [<-, bend left] node [l,above] {$b$} (A);
    \node at (1,1) {$=$};
    \node at (.6,1.4) {$\Downarrow l^{\circ-1}$};
    \node at (1.3,.5) {$\Downarrow l^{\circ}$};
  \end{tikzpicture}
\end{center}

Because all these coherence laws hold in $\Span_2(T),$ it is a symmetric monoidal bicategory.

\newcommand{\pb}[1]{\begin{tikzpicture}[dot/.style={fill,inner sep=1pt,circle}, scale=.2]
  \xdef\lastx{0}
  \foreach \x / \y [count=\n] in {#1} {
    \draw[draw=\ifnum\y>0 black\else white\fi] (\lastx, 0) to[out=90,in=90,looseness=1.7] (\x, 0);
    \xdef\lastx{\x}
  }
  \foreach \x in {0,...,\n} \draw[black] (\x,0) node[dot]{};
\end{tikzpicture}}

In order to prove that the swallowtail coherence law holds, we have to demonstrate an equation between two maps of spans for every object $A$ in $T$.  These maps go between spans whose legs are not necessarily isomorphisms, so the approach taken above will not work to prove that the swallowtail coherence law holds.  Each leg is, however, a natural transformation: either a unitor, an associator, duplication, deletion, a projection, or some product of these.  The feet and apexes of the spans are cartesian products involving only copies of $A$ and the terminal object $1$.

As a calculational aid, we introduce some topological notation for weak pullbacks.  We use one dot for each of the projections from the weak pullback to $A$ that it comes equipped with, and we use an arc for each 2-isomorphism between two projections.  We will denote the terminal object by $1$.

Some examples, assuming $A$ is not terminal:
\begin{enumerate}
  \item We denote $1$ by $1$.
  \item We denote $A$ by $\begin{tikzpicture}[dot/.style={fill,inner sep=1pt,circle}, scale=.2]\draw[black](0,0)node[dot]{};\end{tikzpicture}$.
  \item We denote $A \times 1$ by $\begin{tikzpicture}[dot/.style={fill,inner sep=1pt,circle}, scale=.2]\draw[black](0,0)node[dot]{};\end{tikzpicture} 1$.
  \item \label{product} We denote $A \times A$ by $\pb{1/0}$.
  \item \label{pullback} We denote $A^{\circ 2}$ by $\pb{1/1}$.
  \item \label{firstthree}The weak pullback of the cospan $A^{\circ 2} \stackrel{\pi_2}{\to} A \stackrel{\pi_2}{\leftarrow} A\times A$ is the object $(A^{\circ_2})_{\pi_2, \pi_2}(A\times A)$ equipped with morphisms $$\pi_1, \pi_2, \pi_3, \pi_4\maps (A^{\circ_2})_{\pi_2, \pi_2}(A\times A) \to A$$ and 2-isomorphisms $$K_1\maps \pi_1 \Rightarrow \pi_2$$ and $$K_2\maps \pi_2 \Rightarrow \pi_4.$$  We denote the object $(A^{\circ_2})_{\pi_2, \pi_2}(A\times A)$ by $$\pb{1/1,3/1,2/0}.$$  Note that we form this diagram by juxtaposing examples \ref{pullback} and \ref{product} and adding an arc between the second dot in each pair.
  \item \label{four} We denote $A^{\circ 4}$ by $\pb{1/1,2/1,3/1}$.
  \item The weak pullback of the cospan $A\times A \stackrel{A \times \Delta}{\to} A \times (A \times A) \stackrel{a \circ (\Delta \times A)}{\leftarrow} A\times A$ is the object ${(A \times A)_{A \times \Delta, a \circ (\Delta \times A)}(A\times A)}$ equipped with morphisms $$\pi_1, \pi_2, \pi_3, \pi_4\maps {(A \times A)_{A \times \Delta, a \circ (\Delta \times A)}(A\times A)} \to A$$ and 2-morphisms $$K_1\maps \pi_1 \Rightarrow \pi_3,$$ $$K_2\maps \pi_2 \Rightarrow \pi_3,$$ and $$K_3\maps \pi_2 \Rightarrow \pi_4.$$  We denote the object $(A \times A)_{\pi_2, \pi_2}(A\times A)$ by $$\pb{2/1,1/1,3/1}.$$  Note that we form this diagram by juxtaposing two copies of example \ref{product} and adding three arcs.  This object is isomorphic to example \ref{four} by ${A \times b \times A.}$
\end{enumerate}

To show that $\Span_2(T)$ is compact closed, we have to show the existence of the 1-morphisms $i$ and $e$, the existence of the 2-morphisms $\zeta$ and $\theta$, and show that $\zeta$ and $\theta$ satisfy the swallowtail coherence law.  The real meat of the proof will be in showing that $\zeta$ (and therefore $\theta$) can be defined in terms of an identity span much like the pentagonator and other 2-morphisms above; the ``dressing'' of this span with unitors for composition follows very much as above.

The cap $i\maps I \to A \tensor A^*$ is the span $$1 \stackrel{!}{\leftarrow} A \stackrel{\Delta}{\to} A\times A;$$ the cup $e\maps A^* \tensor A \to I$ is its reverse $i^\bullet$, $$A\times A \stackrel{\Delta}{\leftarrow} A \stackrel{!}{\to} 1.$$  When $A$ is terminal, we define $\zeta_A$ and $\theta_A$ to be the unique 2-morphism on the unique morphism from $A$ to itself.

To define $\zeta_A$ when $A$ is not terminal, we start with an identity map of spans.  The source span is $$A \stackrel{\pi_1}{\leftarrow} A^{\circ 10} \stackrel{\pi_{10}}{\to} A.$$  The target span is the composite $$(r^{-1})^\bullet \circ (A\tensor e)\circ a\circ (i\tensor A)\circ l^{-1} \circ A,$$
where by $r^{-1}$ we mean the span $A \leftarrow A \stackrel{r^{-1}}{\to} A \times 1,$ and similarly for $l^{-1}.$  To see that it is, in fact, an identity map of spans, consider the target span.  In the diagrams below, we elide the 2-isomorphisms for clarity; we also denote the morphism $\langle\pi_i, \ldots, \pi_j\rangle$ out of a weak pullback by $\pi_{i-j}$.

We start building the composite span by composing the spans $A$ and $l^{-1}$.  The cospan in the composite is the identity on $A$, so the apex is $A^{\circ 2}:$
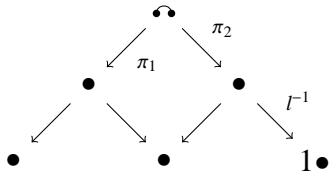
\begin{center}
  \begin{tikzpicture}
    \node (Al) at (0,0) {$\bullet$};
    \node (Ar) at (2,0) {$\bullet$};
    \node (Aa) at (1,1) {$\bullet$}
      edge [->] (Al)
      edge [->] (Ar);
    \node (lr) at (4,0) {$1 \bullet$};
    \node (la) at (3,1) {$\bullet$}
      edge [->] (Ar)
      edge [->] node[auto,l] {$l^{-1}$}(lr);
    \node (Ala) at (2,2) {$\pb{1/1}$}
      edge [->] node[auto,l] {$\pi_1$} (Aa)
      edge [->] node[auto,l] {$\pi_2$} (la);
  \end{tikzpicture}
\end{center}
Next, we compose with $i \tensor A;$ this cospan is the same as example \ref{firstthree} except for the addition of an irrelevant terminal object at the nadir of the cospan:
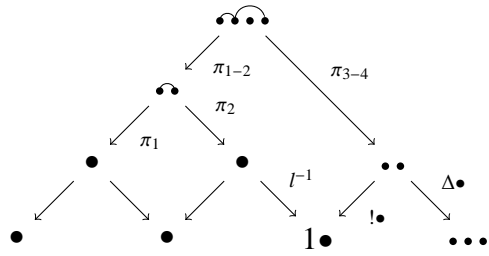
\begin{center}
  \begin{tikzpicture}
    \node (Al) at (0,0) {$\bullet$};
    \node (Ar) at (2,0) {$\bullet$};
    \node (Aa) at (1,1) {$\bullet$}
      edge [->] (Al)
      edge [->] (Ar);
    \node (lr) at (4,0) {$1 \bullet$};
    \node (la) at (3,1) {$\bullet$}
      edge [->] (Ar)
      edge [->] node[auto,l] {$l^{-1}$}(lr);
    \node (Ala) at (2,2) {$\pb{1/1}$}
      edge [->] node[auto,l] {$\pi_1$} (Aa)
      edge [->] node[auto,l] {$\pi_2$} (la);
    \node (ir) at (6,0) {$\pb{1/0,2/0}$};
    \node (ia) at (5,1) {$\pb{1/0}$}
      edge [->] node [auto,l] {$! \bullet$} (lr)
      edge [->] node [auto,l] {$\Delta  \bullet$} (ir);
    \node (Alia) at (3,3) {$\pb{1/1,3/1,2/0}$}
      edge [->] node [auto,l] {$\pi_{1-2}$} (Ala)
      edge [->] node [auto,l] {$\pi_{3-4}$} (ia);
  \end{tikzpicture}
\end{center}
Next, we compose with the associator:
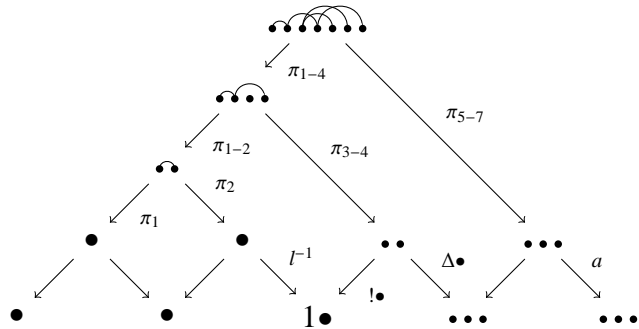
\begin{center}
  \begin{tikzpicture}
    \node (Al) at (0,0) {$\bullet$};
    \node (Ar) at (2,0) {$\bullet$};
    \node (Aa) at (1,1) {$\bullet$}
      edge [->] (Al)
      edge [->] (Ar);
    \node (lr) at (4,0) {$1 \bullet$};
    \node (la) at (3,1) {$\bullet$}
      edge [->] (Ar)
      edge [->] node[auto,l] {$l^{-1}$}(lr);
    \node (Ala) at (2,2) {$\pb{1/1}$}
      edge [->] node[auto,l] {$\pi_1$} (Aa)
      edge [->] node[auto,l] {$\pi_2$} (la);
    \node (ir) at (6,0) {$\pb{1/0,2/0}$};
    \node (ia) at (5,1) {$\pb{1/0}$}
      edge [->] node [auto,l] {$! \bullet$} (lr)
      edge [->] node [auto,l] {$\Delta  \bullet$} (ir);
    \node (Alia) at (3,3) {$\pb{1/1,3/1,2/0}$}
      edge [->] node [auto,l] {$\pi_{1-2}$} (Ala)
      edge [->] node [auto,l] {$\pi_{3-4}$} (ia);
    \node (ir') at (8,0) {$\pb{1/0,2/0}$};
    \node (aa) at (7,1) {$\pb{1/0,2/0}$}
      edge [->] (ir)
      edge [->] node [auto,l] {$a$} (ir');
    \node (Aliaa) at (4,4) {$\pb{1/1,3/1,6/1,5/0,2/1,4/1}$}
      edge [->] node [auto,l] {$\pi_{1-4}$} (Alia)
      edge [->] node [auto,l] {$\pi_{5-7}$} (aa);
  \end{tikzpicture}
\end{center}
Next, we compose with $A \tensor e$: 
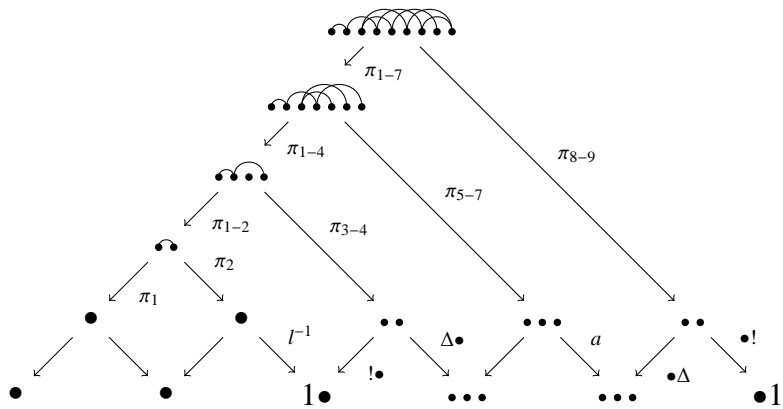
\begin{center}
  \begin{tikzpicture}
    \node (Al) at (0,0) {$\bullet$};
    \node (Ar) at (2,0) {$\bullet$};
    \node (Aa) at (1,1) {$\bullet$}
      edge [->] (Al)
      edge [->] (Ar);
    \node (lr) at (4,0) {$1 \bullet$};
    \node (la) at (3,1) {$\bullet$}
      edge [->] (Ar)
      edge [->] node[auto,l] {$l^{-1}$}(lr);
    \node (Ala) at (2,2) {$\pb{1/1}$}
      edge [->] node[auto,l] {$\pi_1$} (Aa)
      edge [->] node[auto,l] {$\pi_2$} (la);
    \node (ir) at (6,0) {$\pb{1/0,2/0}$};
    \node (ia) at (5,1) {$\pb{1/0}$}
      edge [->] node [auto,l] {$! \bullet$} (lr)
      edge [->] node [auto,l] {$\Delta  \bullet$} (ir);
    \node (Alia) at (3,3) {$\pb{1/1,3/1,2/0}$}
      edge [->] node [auto,l] {$\pi_{1-2}$} (Ala)
      edge [->] node [auto,l] {$\pi_{3-4}$} (ia);
    \node (ir') at (8,0) {$\pb{1/0,2/0}$};
    \node (aa) at (7,1) {$\pb{1/0,2/0}$}
      edge [->] (ir)
      edge [->] node [auto,l] {$a$} (ir');
    \node (Aliaa) at (4,4) {$\pb{1/1,3/1,6/1,5/0,2/1,4/1}$}
      edge [->] node [auto,l] {$\pi_{1-4}$} (Alia)
      edge [->] node [auto,l] {$\pi_{5-7}$} (aa);
    \node (er) at (10,0) {$\bullet 1$};
    \node (ea) at (9,1) {$\pb{1/0}$}
      edge [->] node [auto,l] {$\bullet \Delta$} (ir')
      edge [->] node [auto,l] {$\bullet !$} (er);
    \node (Aliaea) at (5,5) {$\pb{1/1,3/1,6/1,8/1,5/1,2/1,4/1,7/1}$}
      edge [->] node [auto,l] {$\pi_{1-7}$} (Aliaa)
      edge [->] node [auto,l] {$\pi_{8-9}$} (ea);
  \end{tikzpicture}
\end{center}

Finally, we compose with $(r^{-1})^\bullet$:
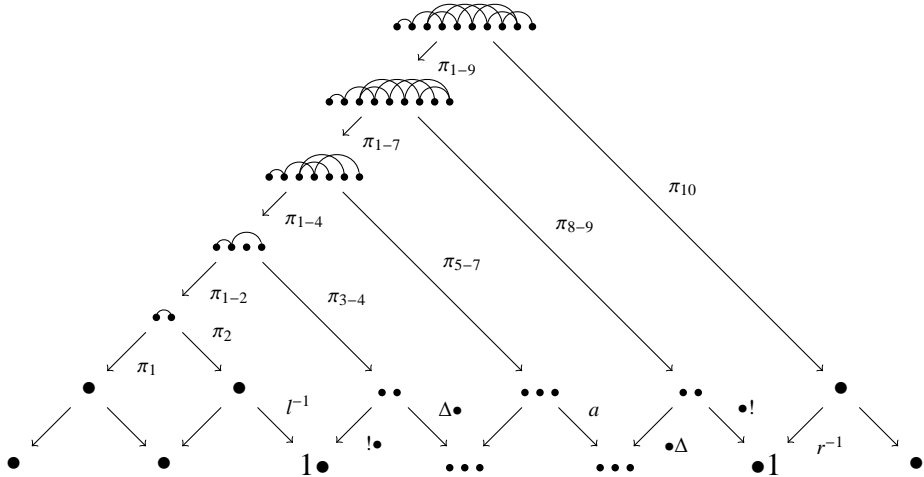
\begin{center}
  \begin{tikzpicture}
    \node (Al) at (0,0) {$\bullet$};
    \node (Ar) at (2,0) {$\bullet$};
    \node (Aa) at (1,1) {$\bullet$}
      edge [->] (Al)
      edge [->] (Ar);
    \node (lr) at (4,0) {$1 \bullet$};
    \node (la) at (3,1) {$\bullet$}
      edge [->] (Ar)
      edge [->] node[auto,l] {$l^{-1}$}(lr);
    \node (Ala) at (2,2) {$\pb{1/1}$}
      edge [->] node[auto,l] {$\pi_1$} (Aa)
      edge [->] node[auto,l] {$\pi_2$} (la);
    \node (ir) at (6,0) {$\pb{1/0,2/0}$};
    \node (ia) at (5,1) {$\pb{1/0}$}
      edge [->] node [auto,l] {$! \bullet$} (lr)
      edge [->] node [auto,l] {$\Delta  \bullet$} (ir);
    \node (Alia) at (3,3) {$\pb{1/1,3/1,2/0}$}
      edge [->] node [auto,l] {$\pi_{1-2}$} (Ala)
      edge [->] node [auto,l] {$\pi_{3-4}$} (ia);
    \node (ir') at (8,0) {$\pb{1/0,2/0}$};
    \node (aa) at (7,1) {$\pb{1/0,2/0}$}
      edge [->] (ir)
      edge [->] node [auto,l] {$a$} (ir');
    \node (Aliaa) at (4,4) {$\pb{1/1,3/1,6/1,5/0,2/1,4/1}$}
      edge [->] node [auto,l] {$\pi_{1-4}$} (Alia)
      edge [->] node [auto,l] {$\pi_{5-7}$} (aa);
    \node (er) at (10,0) {$\bullet 1$};
    \node (ea) at (9,1) {$\pb{1/0}$}
      edge [->] node [auto,l] {$\bullet \Delta$} (ir')
      edge [->] node [auto,l] {$\bullet !$} (er);
    \node (Aliaea) at (5,5) {$\pb{1/1,3/1,6/1,8/1,5/1,2/1,4/1,7/1}$}
      edge [->] node [auto,l] {$\pi_{1-7}$} (Aliaa)
      edge [->] node [auto,l] {$\pi_{8-9}$} (ea);
    \node (rnbr) at (12,0) {$\bullet$};
    \node (rnba) at (11,1) {$\bullet$}
      edge [->] node [auto,l] {$r^{-1}$} (er)
      edge [->] (rnbr);
    \node (Aliaernba) at (6,6) {$\pb{1/1,3/1,6/1,8/1,5/1,2/1,4/1,7/1,9/1}$}
      edge [->] node [auto,l] {$\pi_{1-9}$} (Aliaea)
      edge [->] node [auto,l] {$\pi_{10}$} (rnba);
  \end{tikzpicture}
\end{center}

Inspection of the composite span above shows that none of the cospans involving $\Delta$ are of the form $A^{\circ n} \stackrel{f\pi_n}{\to} C \stackrel{g\pi_1}{\leftarrow} B^{\circ m}$ where $f$ and $g$ are isomorphisms, so the choice of weak pullback for those cospans does not matter there.  The apex $\pb{1/1,3/1,6/1,8/1,5/1,2/1,4/1,7/1,9/1}$ of this composite is made up of ten dots connected by nine arcs in a single chain.  It is evident that $\pb{1/1,3/1,6/1,8/1,5/1,2/1,4/1,7/1,9/1}$ can be permuted to $A^{\circ 10} = \pb{1/1,2/1,3/1,4/1,5/1,6/1,7/1,8/1,9/1}$.  By Corollary \ref{permute}, there is a map of spans
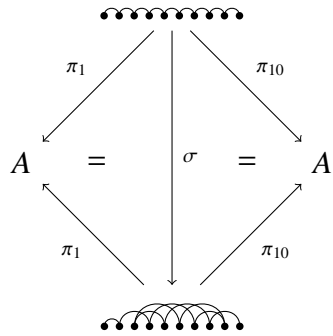
\begin{center}
  \begin{tikzpicture}[scale=2]
    \node (A1) at (0,1) {$A$};
    \node (A10) at (2,1) {$A$};
    \node (BOT) at (1,0) {$\pb{1/1,3/1,6/1,8/1,5/1,2/1,4/1,7/1,9/1}$}
      edge [->] node [below left, l] {$\pi_1$} (A1)
      edge [->] node [below right, l] {$\pi_{10}$} (A10);
    \node (TOP) at (1,2) {$\pb{1/1,2/1,3/1,4/1,5/1,6/1,7/1,8/1,9/1}$}
      edge [->] node [above left, l] {$\pi_1$} (A1)
      edge [->] node [above right, l] {$\pi_{10}$} (A10)
      edge [->] node [right,l] {$\sigma$} (BOT);
    \node at (.5,1) {$=$};
    \node at (1.5,1) {$=$};
  \end{tikzpicture}
\end{center}
in the same equivalence class as the identity.

We define the 2-morphism $\zeta$ to be the composite of this identity map of spans with
\begin{enumerate}
  \item inverse unitors for composition on the source morphism mapping from the identity span on $A$ to the span $A \stackrel{\pi_1}{\leftarrow} A^{\circ 10} \stackrel{\pi_{10}}{\to} A,$
  \item a unitor for composition on the target morphism mapping $(l^{-1}\circ A)$ to $l^\bullet$, similar to what we did when defining the 2-morphism $S$, and
  \item \label{switch} the isomorphism of spans
  \begin{center}
    \begin{tikzpicture}[scale=2]
      \node (A1) at (0,1) {$A\times 1$};
      \node (A10) at (2,1) {$A$};
      \node (BOT) at (1,0) {$A\times 1$}
        edge [->] (A1)
        edge [->] node [below right, l] {$r$} (A10);
      \node (TOP) at (1,2) {$A$}
        edge [->] node [above left, l] {$r^{-1}$} (A1)
        edge [->] (A10)
        edge [->] node [right,l] {$r^{-1}$} (BOT);
      \node at (.5,1) {$=$};
      \node at (1.5,1) {$=$};
    \end{tikzpicture}
  \end{center}
  on the target.
\end{enumerate}
The 2-morphism $\theta_A$ follows {\em mutatis mutandis}.

In the left hand-side of the swallowtail coherence law, the only parts not accounted for by the coherence theorem for bicategories are the two uses of the isomorphism of spans in item \ref{switch} above: once in $\zeta\tensor A$ and once in $A^*\tensor \theta^{-1}$.  The composite isomorphism of spans is
\begin{center}
  \begin{tikzpicture}[scale=2]
    \node (LEFT) at (0,1.5) {$(A\times 1) \times A$};
    \node (RIGHT) at (2,1.5) {$A\times A$};
    \node (BOT) at (1,0) {$A\times A$}
      edge [->] node [below left, l] {$r^{-1}\times A$} (LEFT)
      edge [->] (RIGHT);
    \node (MID1) at (1,1) {$A\times (1 \times A)$}
      edge [->] node [right, l,pos=.2] {$A\times l$} (BOT);
    \node (MID2) at (1,2) {$(A\times 1) \times A$}
      edge [->] node [right, l] {$a$} (MID1);
    \node (TOP) at (1,3) {$A\times A$}
      edge [->] node [above left, l] {$r^{-1}\times A$} (LEFT)
      edge [->] (RIGHT)
      edge [->] node [right,l,pos=.9] {$r^{-1}\times A$} (MID2);
    \node at (.6,1.5) {$=$};
    \node at (1.5,1.5) {$=$};
  \end{tikzpicture}
\end{center}
Because the triangle laws hold in $T$, the composite isomorphism is the identity.  Therefore the swallowtail coherence law holds in $\Span_2(T)$ and $\Span_2(T)$ is compact closed. \hfill $\square$

\begin{cor}
  \label{span}
  When $C$ is a category with finite products and pullbacks, the bicategory Span($C$) of objects of $C$, spans in $C,$ and maps of spans is compact closed.
\end{cor}
{\em Proof.} When $C$ is a category with finite products and pullbacks,
$\Span(C)$ is a special case of Theorem \ref{Span2} where all the
2-morphisms in the weak pullbacks are identities. \hfill $\square$

\begin{cor}
  The bicategories Cospan(ResNet) and Circ are compact closed.
\end{cor}
{\em Proof.} The coproduct of two resistor networks is given by juxtaposition; the pushout of a cospan $S \hookleftarrow R
\hookrightarrow T$ of resistor networks is given by juxtaposition followed by identifying the images of $R$ in $S$ and $T$.  Cospans in ResNet are spans in ResNet${}^{\op},$ where the coproduct and pushout become product and pullback, so Cospan(ResNet) is compact closed by the previous corollary.  Since every object is self-dual in Cospan(ResNet), the subcategory Circ whose objects are resistor networks with no edges is also compact closed. \hfill $\square$


\begin{thebibliography}{99}

\bibitem{Asada} Kazuyuki Asada, Arrows are strong monads, 
\textsl{MSFP '10} (2002) 33--42.  Also available at
\href{http://www.ipl.t.u-tokyo.ac.jp/~asada/papers/arrStrMnd.pdf}
{http://www.ipl.t.u-tokyo.ac.jp/$\sim$asada/papers/arrStrMnd.pdf}.

\bibitem{AtiyahTQFT} Michael Atiyah, Topological quantum field theories,
\textsl{Publications Math\'ematiques de l'Institut des Hautes \'Etudes
Scientifiques.} \textbf{68:1} (1988) 175--186.

\bibitem{HDATQFT} John C. Baez and James Dolan, Higher-dimensional algebra and
topological quantum field theory, \textsl{Journal of Mathematical Physics} \textbf{36}
(1995) 6073--6105.  Also available at \href
{http://arxiv.org/abs/q-alg/9503002v2}{arXiv:q-alg/9503002v2}.

\bibitem{HDA1} John C. Baez and Martin Neuchl, Higher-dimensional algebra I:
Braided monoidal 2-categories, \textsl{Advances in Mathematics} \textbf{121} (1996)
196--244. Also available at \href{http://arxiv.org/abs/q-alg/9511013}
{arXiv:q-alg/9511013}.

\bibitem{HDA3} John C. Baez and James Dolan, Higher-dimensional algebra III:
n-categories and the algebra of opetopes, \textsl{Advances in Mathematics} \textbf{135}
(1998) 145--206.  Also available at \href
{http://arxiv.org/abs/q-alg/9702014}{arXiv:q-alg/9702014}.

\bibitem{HDA4} John C. Baez and Laurel Langford, Higher-dimensional algebra IV:
2-Tangles, \textsl{Advances in Mathematics} \textbf{180} (2003) 705--764. Also
available at \href{http://arxiv.org/abs/math.QA/9811139}
{arXiv:math.QA/9811139}.

\bibitem{BaezWiseCrans}John C. Baez, Derek K. Wise, and Alissa S. Crans, Exotic
statistics for strings in 4d BF theory, 
\textsl{Advances in Theoretical and Mathematical Physics}
\textbf{11} (2007) 707--749.  Also available at \href
{http://arxiv.org/abs/gr-qc/0603085}{arXiv:gr-qc/0603085}.

\bibitem{Bartlett} Bruce Bartlett, Quasistrict symmetric monoidal 2-categories via wire diagrams.  
Available at \href{http://arxiv.org/abs/1409.2148}{arXiv:1409.2148}.

\bibitem{Batanin}Michael A. Batanin, The Eckmann-Hilton argument and higher operads.  Available at
\href{http://arxiv.org/abs/math/0207281}{arXiv:math/0207281}.

\bibitem{Ben00} Jean B{\'e}nabou, Distributors at work.
Available at \\
\href{http://www.mathematik.tu-darmstadt.de/~streicher/FIBR/DiWo.pdf}
{http://www.mathematik.tu-darmstadt.de/$\sim$streicher/FIBR/DiWo.pdf}.

\bibitem{Ben67} Jean B{\'e}nabou, Introduction to bicategories, 
\textsl{Reports of the Midwest Category Seminar} (1967) 1--77.

\bibitem{Ben73} Jean B{\'e}nabou, Les distributeurs, 
\textsl{S\'eminaires de Math\'ematique Pure Universit\'e Catholique de Louvain} 
\textbf{33} (1973).

\bibitem{BD86} Francis Borceux and Dominique Dejean, 
Cauchy completion in category theory, 
\textsl{Cahiers de Topologie et G\'eom\'etrie Diff\'erentielle Cat\'egoriques} 
\textbf{27} (1986) 133--146.

\bibitem{Breen} Lawrence Breen, On the classification of 2-gerbes and
2-stacks, \textsl{Ast\'erisque} \textbf{225} (1994) 1--160.

\bibitem{CKWW08} Aurelio Carboni, G.\ Maxwell Kelly, Robert F.\ C.\ Walters, and
Richard J.\ Wood, Cartesian bicategories II, \textsl{Theory and Applications
of Categories} {\bf 19:6} (2008) 93--124.  Available at \href
{http://www.emis.de/journals/TAC/volumes/19/6/19-06abs.html}
{http://www.emis.de/journals/TAC/volumes/19/6/19-06abs.html}.

\bibitem{CW87} Aurelio Carboni and Robert F.\ C.\ Walters, Cartesian
bicategories I, \textsl{Journal of Pure and Applied Algebra} 
\textbf{49} (1987) 11--32.

\bibitem{CW04} Gian Luca Cattani and Glynn Winskel, Profunctors, open maps
and bisimulation, \textsl{Mathematical Structures in Computer Science}
\textbf{15:3} (2005) 553--614.

\bibitem{ChengGurskinCob} Eugenia Cheng and Nick Gurski, Towards an
$n$-category of cobordisms, \textsl{Theory and Applications of
Categories} \textbf{18:10} (2007) 274--302.

\bibitem{Crans} Sjoerd E.\ Crans, Generalized centers of braided and sylleptic
monoidal 2-categories, \textsl{Advances in Mathematics} 
\textbf{136} (1998) 183--223.

\bibitem{Day} Brian J.\ Day. On closed categories of functors. 
\textsl{Lecture Notes in Mathematics} \textbf{137} (1970) 1--38.

\bibitem{DS97} Brian J.\ Day and Ross Street, Monoidal bicategories and Hopf
algebroids, \textsl{Advances in Mathematics} \textbf{129} (1997) 99-157.

\bibitem{FLObCatComplex} Marcelo Fiore, Tom Leinster, Objects of categories
as complex numbers, \textsl{Advances in Mathematics} \textbf{190} (2005) 264--277. 
Available at \href{http://arxiv.org/abs/math/0212377}
{arXiv:math/0212377}.

\bibitem{FKLW} Michael H.\ Freedman, Alexei Kitaev, Michael J.\ Larsen, 
and Zhenghan Wang, Topological quantum computation.  Available at \href
{http://arxiv.org/abs/quant-ph/0101025}{arXiv:quant-ph/0101025v2}.

\bibitem{GK} Nicola Gambino and Joachim Kock, 
Polynomial functors and polynomial monads,
\textsl{Mathematical Proceedings of the Cambridge Philosophical Society} 
\textbf{154} (2013) 153--192. 
Also available at \href{http://arxiv.org/abs/0906.4931}{arXiv:0906.4931}.

\bibitem{Gray} John W.\ Gray, \textsl{Formal category theory: Adjointness for
2-Categories,} Lecture Notes in Mathematics \textbf{391} (1974).

\bibitem{GPS} Robert Gordon, A.\ John Power, and Ross Street, Coherence for
tricategories, \textsl{Memoirs of the American Mathematical Society} \textbf{117}
(1995).

\bibitem{Permutohedra} Georges-th\'eodule Guilbaud and Pierre Rosenstiehl, 
Analyse alg\'ebrique d'un scrutin, \textsl{Math\'ematiques et Sciences Humaines}
\textbf{4} (1963) 9--33.

\bibitem{Gurski} Nick Gurski, \textsl{An algebraic theory of tricategories},
PhD thesis, University of Chicago (2007).  Available at
\href{http://gauss.math.yale.edu/~mg622/tricats.pdf}
{http://gauss.math.yale.edu/$\sim$mg622/tricats.pdf}

\bibitem{GO} Nick Gurski and Ang\'elica M.\ Osorno, 
Infinite loop spaces, and coherence for symmetric 
monoidal bicategories, \textsl{Advances in Mathematics} 
\textbf{246} (2013) 1--32.  Also available at 
\href{http://arxiv.org/abs/1210.1174v1}{arXiv:1210.1174v1}.

\bibitem{Hoffnung} Alex Hoffnung, Spans in 2-categories: A monoidal
tricategory. \textsl{Theory and Applications of Categories,} to appear.  
Also available at \href{http://arxiv.org/abs/1112.0560}{arXiv:1112.0560}

\bibitem{Houston} Robin Houston, \textsl{Linear Logic without Units},
PhD Thesis, University of Manchester (2007).  Available at \href
{http://arxiv.org/abs/1305.2231}{arXiv:1305.2231}

\bibitem{Hughes} John Hughes, Generalising Monads to Arrows, 
\textsl{Science of Computer Programming} \textbf{37} (2000) 67--111. 
Also available at \href
{http://www.cse.chalmers.se/~rjmh/Papers/arrows.pdf}
{http://www.cse.chalmers.se/$\sim$rjmh/Papers/arrows.pdf}.

\bibitem{JHH} Bart Jacobs, Chris Heunen, and Ichiro Hasuo, 
Categorical semantics for arrows, 
\textsl{Journal of Functional Programming} 
\textbf{19:3-4} (2009) 403--438.  The results appeared in an 
earlier conference paper C.\ Heunen, B.\ Jacobs, Arrows, like 
Monads, are Monoids, {\sl ENTCS} \textbf{158} (2006) 219--236.

\bibitem{JSV} Andr\'e Joyal, Ross Street, and Dominic Verity, 
Traced monoidal categories, 
\textsl{Mathematical Proceedings of the Cambridge Philosophical Society} 
\textbf{119:3} (1996) 447--468.

\bibitem{KV94} Mikhail Kapranov and Vladimir Voevodsky, 2-categories and
Zamolodchikov tetrahedra equations, 
\textsl{Proceedings of Symposia in Pure Mathematics}
\textbf{56} (1994) 177--259.

\bibitem{KSW98} Piergiulio Katis, Nicoletta Sabadini, 
Robert F.\ C.\ Walters, On partita doppia, 
\textsl{Theory and Applications of Categories.} (Written in 1998,
to appear.)  Also available at \href{http://arxiv.org/abs/0803.2429}
{arXiv:0803.2429v1}.

\bibitem{KW99} Piergiulio Katis, Robert F.\ C.\ Walters, 
The compact closed bicategory of left adjoints, 
\textsl{Mathematical Proceedings of the Cambridge Philosophical Society} 
\textbf{130:1} (2001) 77--87.  Also available at \href
{http://citeseerx.ist.psu.edu/viewdoc/summary?doi=10.1.1.23.4825}
{http://citeseerx.ist.psu.edu/viewdoc/summary?doi=10.1.1.23.4825}.

\bibitem{Kelly1982} G.\ Maxwell Kelly, Basic concepts of enriched category
theory, \textsl{Lecture Notes in Mathematics} \textbf{64} (1982). 
Also available at \\
\href{http://www.tac.mta.ca/tac/reprints/articles/10/tr10abs.html}
{http://www.tac.mta.ca/tac/reprints/articles/10/tr10abs.html}.

\bibitem{Kelly1974} G.\ Maxwell Kelly, Coherence theorems for lax algebras and
distributive laws, \textsl{Lecture Notes in Mathematics} \textbf{420}
(1974) 281--375.

\bibitem{KellyCC} G.\ Maxwell Kelly, Many variable functorial calculus I,
\textsl{Lecture Notes in Mathematics} \textbf{281} (1972) 66--105.

\bibitem{KellyLaplaza} G.\ Maxwell Kelly and M.\ Laplaza, 
Coherence for compact closed categories, 
\textsl{Journal of Pure and Applied Algebra} \textbf{19} (1980) 193--213.

\bibitem{Lack} Stephen Lack, A 2-categories companion, \textsl
{Towards Higher Categories} (2009).  Also available at 
\href{http://arxiv.org/abs/math/0702535}{arXiv:math/0702535}.

\bibitem{Laplaza} M.\ Laplaza, Coherence for distributivity, 
\textsl{Lecture Notes in Mathematics} \textbf{281} (1972) 29--72.

\bibitem{LeinsterBB} Tom Leinster, Basic bicategories. Available
at \href{http://arxiv.org/abs/math/9810017}{arXiv:math/9810017}.

\bibitem{McCrudden} Paddy McCrudden, Balanced coalgebroids, \textsl{Theory
and Applications of Categories}, \textbf{7} 6 (2000) 71--147. 
Available at \\
\href{http://www.emis.de/journals/TAC/volumes/7/n6/7-06abs.html}
{http://www.emis.de/journals/TAC/volumes/7/n6/7-06abs.html}.

\bibitem{Morton} Jeffrey Morton, Double bicategories and double cospans,
\textsl{Journal of Homotopy and Related Structures} \textbf{4:1} (2009)
 389--428.  Also available at \href{http://arxiv.org/abs/math/0611930}
{arXiv:math/0611930}.

\bibitem{MortonVicary} Jeffrey Morton and Jamie Vicary, The Categorified
Heisenberg Algebra I: A Combinatorial Representation.  Available at
\href{http://arxiv.org/abs/1207.2054}{arXiv:1207.2054}.

\bibitem{nlab2pullbacks} \textsl{2-pullbacks}, 
\href{http://ncatlab.org/nlab/show/2-pullback#strict_weighted_limits}
{http://ncatlab.org/nlab/show/2-pullback\#strict\_weighted\_limits}

\bibitem{PL07} Anne Preller and Joachim Lambek, Free compact 2-categories,
\textsl{Mathematical Structures in Computer Science} \textbf{17} 2
(2007) 309--340.  Also available at \\ \href
{ftp://ftp.math.mcgill.ca/pub/lambek/freecompact.pdf}
{ftp://ftp.math.mcgill.ca/pub/lambek/freecompact.pdf}.

\bibitem{SPT} Chris Schommer-Pries, \textsl{The classification of two-dimensional
extended topological field theories}, PhD thesis, UC Berkeley (2003). 
Available at \\
\href{http://sites.google.com/site/chrisschommerpriesmath/Home/Schommer-Pries-Thesis.pdf}
{http://sites.google.com/site/chrisschommerpriesmath/Home/Schommer-Pries-Thesis.pdf}.

\bibitem{Shulman} M.\ Shulman, Constructing symmetric monoidal
bicategories.  Avaliable at \\ \href{http://arxiv.org/abs/1004.0993}
{arXiv:1004.0993}.

\bibitem{ShulCafe} M.\ Shulman, Re: Inevitability in Mathematics,
\textsl{$n$-Category Caf\'e} (2010).  Available at
\href{https://golem.ph.utexas.edu/category/2010/06/inevitability_in_mathematics.html#c033741}
{https://golem.ph.utexas.edu/category/2010/06/inevitability\_in\_mathematics.html\#c033741}

\bibitem{Stasheff} J.\ Stasheff, Homotopy associativity of $H$-spaces I,
II, \textsl{Transactions of the American Mathematical Society} \textbf{108} (1963) 275--312.

\bibitem{Street82} Ross Street, Fibrations in bicategories, \textsl
{Cahiers de Topologie et G\'eom\'etrie Diff\'erentielle Cat\'egoriques}, 
\textbf{21:2} (1980) 111--160.  Also available at \\
\href{http://www.numdam.org/item?id=CTGDC_1980__21_2_111_0}
{http://www.numdam.org/item?id=CTGDC\_1980\_\_21\_2\_111\_0}

\bibitem{Street03} Ross Street, Functorial calculus in monoidal
bicategories, \textsl{Applied Categorical Structures} \textbf{11} (2003)
219--227.  Also available at 
\href{http://www.maths.mq.edu.au/~street/Extraord.pdf}
{http://www.maths.mq.edu.au/$\sim$street/Extraord.pdf}.

\bibitem{Street04} Ross Street, Frobenius monads and pseudomonoids,
\textsl{Journal of Mathematical Physics} \textbf{45:10} (2004) 3930--3948. 
Also available at
\href{http://www.maths.mq.edu.au/~street/Frob.pdf}
{http://www.maths.mq.edu.au/$\sim$street/Frob.pdf}.

\end{thebibliography}
\end{document}